\documentclass[11pt,fleqn]{article}
\usepackage{amsfonts, epsfig, amssymb, amsmath}
\usepackage{color}
\newcommand{\ol}{\setlength{\itemsep}{0pt.}\begin{enumerate}}
\newcommand{\eol}{\end{enumerate}\setlength{\itemsep}{-\parsep}}
\newcommand{\ignore}[1]{}
\title{A moment ratio bound for polynomials and some extremal properties of Krawchouk polynomials and Hamming spheres}
\author{Naomi Kirshner and Alex Samorodnitsky}

\begin{document}
\date{}
\maketitle


\newtheorem{THEOREM}{Theorem}[section]
\newenvironment{theorem}{\begin{THEOREM} \hspace{-.85em} {\bf :}
}%
                        {\end{THEOREM}}
\newtheorem{LEMMA}[THEOREM]{Lemma}
\newenvironment{lemma}{\begin{LEMMA} \hspace{-.85em} {\bf :} }%
                      {\end{LEMMA}}
\newtheorem{COROLLARY}[THEOREM]{Corollary}
\newenvironment{corollary}{\begin{COROLLARY} \hspace{-.85em} {\bf
:} }%
                          {\end{COROLLARY}}
\newtheorem{PROPOSITION}[THEOREM]{Proposition}
\newenvironment{proposition}{\begin{PROPOSITION} \hspace{-.85em}
{\bf :} }%
                            {\end{PROPOSITION}}
\newtheorem{DEFINITION}[THEOREM]{Definition}
\newenvironment{definition}{\begin{DEFINITION} \hspace{-.85em} {\bf
:} \rm}%
                            {\end{DEFINITION}}
\newtheorem{EXAMPLE}[THEOREM]{Example}
\newenvironment{example}{\begin{EXAMPLE} \hspace{-.85em} {\bf :}
\rm}%
                            {\end{EXAMPLE}}
\newtheorem{CONJECTURE}[THEOREM]{Conjecture}
\newenvironment{conjecture}{\begin{CONJECTURE} \hspace{-.85em}
{\bf :} \rm}%
                            {\end{CONJECTURE}}
\newtheorem{MAINCONJECTURE}[THEOREM]{Main Conjecture}
\newenvironment{mainconjecture}{\begin{MAINCONJECTURE} \hspace{-.85em}
{\bf :} \rm}%
                            {\end{MAINCONJECTURE}}
\newtheorem{PROBLEM}[THEOREM]{Problem}
\newenvironment{problem}{\begin{PROBLEM} \hspace{-.85em} {\bf :}
\rm}%
                            {\end{PROBLEM}}
\newtheorem{QUESTION}[THEOREM]{Question}
\newenvironment{question}{\begin{QUESTION} \hspace{-.85em} {\bf :}
\rm}%
                            {\end{QUESTION}}
\newtheorem{REMARK}[THEOREM]{Remark}
\newenvironment{remark}{\begin{REMARK} \hspace{-.85em} {\bf :}
\rm}%
                            {\end{REMARK}}

\newcommand{\thm}{\begin{theorem}}
\newcommand{\lem}{\begin{lemma}}
\newcommand{\pro}{\begin{proposition}}
\newcommand{\dfn}{\begin{definition}}
\newcommand{\rem}{\begin{remark}}
\newcommand{\xam}{\begin{example}}
\newcommand{\cnj}{\begin{conjecture}}
\newcommand{\mcnj}{\begin{mainconjecture}}
\newcommand{\prb}{\begin{problem}}
\newcommand{\que}{\begin{question}}
\newcommand{\cor}{\begin{corollary}}
\newcommand{\prf}{\noindent{\bf Proof:} }
\newcommand{\ethm}{\end{theorem}}
\newcommand{\elem}{\end{lemma}}
\newcommand{\epro}{\end{proposition}}
\newcommand{\edfn}{\bbox\end{definition}}
\newcommand{\erem}{\bbox\end{remark}}
\newcommand{\exam}{\bbox\end{example}}
\newcommand{\ecnj}{\bbox\end{conjecture}}
\newcommand{\emcnj}{\bbox\end{mainconjecture}}
\newcommand{\eprb}{\bbox\end{problem}}
\newcommand{\eque}{\bbox\end{question}}
\newcommand{\ecor}{\end{corollary}}
\newcommand{\eprf}{\bbox}
\newcommand{\beqn}{\begin{equation}}
\newcommand{\eeqn}{\end{equation}}
\newcommand{\wbox}{\mbox{$\sqcap$\llap{$\sqcup$}}}
\newcommand{\bbox}{\vrule height7pt width4pt depth1pt}
\newcommand{\qed}{\bbox}
\def\sup{^}

\def\H{\{0,1\}^n}

\def\S{S(n,w)}

\def\g{g_{\ast}}
\def\xop{x_{\ast}}
\def\y{y_{\ast}}
\def\z{z_{\ast}}

\def\f{\tilde f}

\def\n{\lfloor \frac n2 \rfloor}

\def \E{\mathop{{}\mathbb E}}
\def \R{\mathbb R}
\def \Z{\mathbb Z}
\def \F{\mathbb F}
\def \T{\mathbb T}

\def \x{\textcolor{red}{x}}
\def \r{\textcolor{red}{r}}
\def \Rc{\textcolor{red}{R}}

\def \noi{{\noindent}}

\def \iff{~~~~\Leftrightarrow~~~~}

\def\<{\left<}
\def\>{\right>}
\def\({\left(}
\def\){\right)}

\def\myblt{\noi --\, }

\def \e{\epsilon}
\def \l{\lambda}

\def \I{{\cal I}}

\def \queq {\quad = \quad}

\def\Tp{Tchebyshef polynomial}
\def\Tps{TchebysDeto be the maximafine $A(n,d)$ l size of a code with distance $d$hef polynomials}
\newcommand{\rarrow}{\rightarrow}

\newcommand{\larrow}{\leftarrow}

\overfullrule=0pt
\def\setof#1{\lbrace #1 \rbrace}

\begin{abstract}
Let $p \ge 2$. We improve the bound $\frac{\|f\|_p}{\|f\|_2} \le (p-1)^{s/2}$ for a polynomial $f$ of degree $s$ on the boolean cube $\H$, which comes from hypercontractivity, replacing the right hand side of this inequality by an explicit bivariate function of $p$ and $s$, which is smaller than $(p-1)^{s/2}$ for any $p > 2$ and $s > 0$. We show the new bound to be tight, within a smaller order factor, for the {\it Krawchouk polynomial} of degree $s$.

This implies several nearly-extremal properties of Krawchouk polynomials and Hamming spheres (equivalently, Hamming balls). In particular, Krawchouk polynomials have (almost) the heaviest tails among all polynomials of the same degree and $\ell_2$ norm\footnote{This has to be interpreted with some care.}. The Hamming spheres have the following approximate edge-isoperimetric property: For all $1 \le s \le \frac{n}{2}$, and for all even distances $0 \le i \le \frac{2s(n-s)}{n}$, the Hamming sphere of radius $s$ contains, up to a multiplicative factor of $O(i)$, as many pairs of points at distance $i$ as possible, among sets of the same size\footnote{There is a similar, but slightly weaker and somewhat more complicated claim for general distances.}. This also implies that Hamming spheres are (almost) stablest with respect to noise among sets of the same size. In coding theory terms this means that a Hamming sphere (equivalently a Hamming ball) has the maximal probability of undetected error, among all binary codes of the same rate.

We also describe a family of hypercontractive inequalities for functions on $\H$, which improve on the `usual' ``$q \rarrow 2$" inequality by taking into account the concentration of a function (expressed as the ratio between its $\ell_r$ norms), and which are nearly tight for characteristic functions of Hamming spheres.
\end{abstract}

\section{Introduction}

\noi We prove upper bounds on the moments of polynomials on the discrete cube $\{0,1\}^n$ endowed with uniform measure. Let $H$ be the binary entropy function, and let $\psi(p,x)$ be a function on $[2, \infty) \times \left[0,1/2\right]$, defined by
\[
\psi(p,x) \queq (p-1) + \log_2\Big((1-\delta)^p + \delta^p\Big) - \frac{p}{2} H(x) - px\log_2(1-2\delta),
\]
where $\delta$ is determined by $x = \(\frac12 - \delta\) \cdot \frac{(1-\delta)^{p-1} - \delta^{p-1}}{(1-\delta)^p + \delta^p}$. (It will be shown that $\delta$ is well-defined.)

\noi Then, for $p \ge 2$, $0 \le s \le \frac n2$,\footnote{This is the interesting range of parameters in terms of $s$, since the spaces of homogeneous polynomials of degree $s$ and $n-s$ on $\H$ are isometric, for any $\ell_p$ norm, see Section~\ref{subsubsec:Fourier}.} and for a homogeneous polynomial $f$ of degree $s$ on $\H$,  holds
\beqn
\label{ineq:moments}
\frac{\E |f|^p}{\(\E f^2\)^{\frac p2}} ~\le~ 2^{\psi\(p,\frac sn\) \cdot n}.
\eeqn

\noi We will show this to be an improvement over the well-known bound
\beqn
\label{ineq:moments-HC}
\frac{\E |f|^p}{\(\E f^2\)^{\frac p2}} ~\le~ (p-1)^{\frac{ps}{2}},
\eeqn
which follows from the hypercontractive inequality (\ref{ineq:HC}) below (see e.g. \cite{O'Donnel}). Let $\psi_1(p,x)  = \frac{p \log_2(p-1)}{2} \cdot x$, so that $\psi_1\(p, \frac sn\) = \frac1n \log_2\((p-1)^{\frac{ps}{2}}\)$. We will show that for any fixed $p > 2$ the functions $\psi$ and $\psi_1$ and their first derivatives coincide at $x = 0$ and, moreover, that the function $\psi$ is strongly concave in $x$. This will imply $\psi(p,x) < \psi_1(p,x)$ for any $p > 2$ and $x > 0$.

\noi For a fixed $p > 2$ and for $s \ll n$, the bounds in (\ref{ineq:moments}) and in (\ref{ineq:moments-HC}) are very close to each other, in accord with the fact (\cite{Larsson-Cohn}) that if $s$ is a slowly growing function of $n$, the RHS in (\ref{ineq:moments-HC}) cannot be replaced by $C^s$ with $C < (p-1)^{\frac p2}$. However, if we allow $p$ to grow with $n$, the two bounds can be significantly different, even for small $s$. This will be important in estimates which take into account higher moments of polynomials, as is the cases we discuss below.

\noi Let us also observe that both bounds hold in somewhat higher generality - for all polynomials of degree {\it at most} $s$ on $\H$ (see Corollary~\ref{cor:also balls} below).

\noi We proceed with an informal description of several applications of (\ref{ineq:moments}). The formal statements and a more extensive discussion of these results will be given below, in Section~\ref{subsec:results}. First, it will be shown that (\ref{ineq:moments}) is "nearly tight" (in the sense that will be clarified below) if $f$ is the {\it Krawchouk polynomial} $K_s$ defined by
\[
K_s(x) ~=~ \sum_{S \subseteq [n],|S| = s} (-1)^{\sum_{i \in S} x_i}, \quad \mathrm{for} \quad x = \(x_1,x_2,...x_n\) \in \H.
\]
Recalling that $K_s$ is proportional to the Fourier transform of the characteristic function of the Hamming sphere of radius $s$ around zero, this says, alternatively, that Fourier transforms of Hamming spheres are nearly extremal with respect to (\ref{ineq:moments}). This will be shown to imply that Krawchouk polynomials and Hamming spheres have certain nearly extremal properties, compared to other objects with similar characteristics. Specifically, we will show that, up to at most polynomial in $n$ error, the following facts hold for functions on $\H$:

\begin{itemize}

\item

\noi Krawchouk polynomials  have (almost) the heaviest tails among all polynomials of the same degree and $\ell_2$ norm. That is, for a polynomial $f$ of degree $s$ with $\|f\|_2 = \|K_s\|_2$, and for a threshold $T > 0$ holds
\[
\mathrm{Pr}\{|K_s| \ge T\} ~\gtrsim~ \mathrm{Pr}\{|f| \ge T'\},
\]
where $T'$ is not much larger than $T$. For the exact formulation see Theorem~\ref{thm:tails}.

\item

\noi For any $0 \le s \le \frac n2$ and any even $0 \le i \le \frac{2s(n-s)}{n}$, the Hamming sphere of radius $s$ around $0$ contains (almost) the ``maximal" number of pairs of points at distance $i$, among all sets of the same size. For a general distance $i$, the same holds for the union of two Hamming spheres of consecutive radii.

\noi For the exact formulation see Theorem~\ref{thm:edge-isop}.

\item

\noi For any $p \ge 2$, characteristic functions of Hamming spheres are (almost) stablest with respect to noise among all functions with the same $\ell_1$ and $\ell_p$ norms. That is, let $0 \le s \le \frac n2$, let $f_s$ be the characteristic function of the Hamming sphere of radius $s$ around $0$, and let $f$ be a function with $\|f\|_1 = \|f_s\|_1$ and $\|f\|_p = \|f_s\|_p$. Then,
for the noise operator $T_{\e}$, $0 \le \e \le \frac12$, holds
\[
\<f_s, T_{\e} f_s\> ~\gtrsim~ \<f, T_{\e} f\>.
\]
For the exact formulation see Corollary~\ref{cor:sphere-stable} and the discussion after it.

\item

\noi For any $p \ge 2$, characteristic functions of Hamming spheres have (almost) the largest spectral projections among all functions with the same $\ell_1$ and $\ell_p$ norms. That is, in the notation of the previous item, for any $0 \le s \le \frac n2$ and `many' $0 \le k \le \frac n2$ holds
\[
\|\Pi_k f_s\|_2  ~\gtrsim~ \|\Pi_k f\|_2.
\]
Here $\Pi_k f$ is the orthogonal projection on the span of Walsh-Fourier characters of weight $k$. For the exact formulation see Theorem~\ref{thm:max-proj}.

\end{itemize}

\noi Let us make several comments about these results.

\noi -- \,In all the statements above 'homogeneous polynomials of degree $s$' can be replaced with 'polynomials of degree $s$', and 'Hamming spheres of radius $s$' with 'Hamming balls of radius $s$' (we do not go into details for lack of space, but see Corollary~\ref{cor:also balls}.)

\noi -- \,It can be seen that the last three of the claims above are essentially equivalent to each other.

\noi -- \,The exact formulations of the claims above will be in terms of functional inequalities (for functions on $\H$) involving certain explicit, but rather complicated, functions of two variables. These bivariate functions describe the relevant aspects of behavior of Hamming spheres or of Krawchouk polynomials. For instance, consider the function $\psi$ defined above. As will be seen, $\psi\(p,\frac sn\)$ is the right constant in the exponent of the ratio between the $p^{\mathrm{\small th}}$ moment of the Krawchouk polynomial $K_s$ and the $p/2$-power of its second moment. We point out that the appearance of these functions in the statements of the results indicates that Hamming spheres / Krawchouk polynomials are indeed (almost) extremal objects for these results.

\noi -- \,Continuing from the preceding comment, we observe that while these bivariate functions describe the correct exponential behavior of Hamming spheres or Krawchouk polynomials, they do introduce error, which is polynomial in the dimension $n$ of the discrete cube. This is the cause of imprecision in all of the results above. Let us provide some details. Krawchouk polynomials on $\H$ and $n$-dimensional Hamming spheres are discrete objects (if we view a polynomial as a vector of its coefficients), whose behavior is described by expressions involving binomial coefficients. Hence it cannot be reduced to a simply exponential expression without incurring a certain loss. In our case this (lossy) reduction is achieved by replacing the binomial coefficient ${b \choose a}$ by a larger exponential expression $2^{H\(\frac ab\) \cdot b}$ (see (\ref{binomial-H}) below). This is the main source of loss we incur. For an illustration see Example~\ref{xam:isop} below and observe that the gap between the upper bound and the lower bound given by a Hamming sphere is due solely to replacing two binomial coefficients by corresponding exponential expressions.

\noi -- \,Finally, we observe that a polynomial error will typically be much smaller than the main term in the estimates we discuss, since the approximation of ${b \choose a}$ by $2^{H\(\frac ab\) \cdot b}$ is usually a very good one. However, this fact has to be interpreted with some care, since the significance of an inaccuracy depends on the context. Consider the following two examples.

\xam
\label{xam:isop}

\noi We will show in Theorem~\ref{thm:edge-isop} that if $A$ a subset of $\H$ with $|A| \le {n \choose s}$, and if $0 \le i \le \frac{2s(n-s)}{n}$, then the number of pairs of points at distance $i$ in $A$ is bounded from above by $|A| \cdot 2^{H\(\frac{i}{2s}\) \cdot s + H\(\frac{i}{2(n-s)}\) \cdot (n-s)}$. On the other hand, if $A$ is a Hamming sphere of radius $s$, and if $i$ is even, this number is $|A| \cdot {s \choose {i/2}} {{n-s} \choose {i/2}}$. This, by (\ref{binomial-H}), is at least $\Omega\(\frac 1i\) \cdot |A| \cdot 2^{H\(\frac{i}{2s}\) \cdot s + H\(\frac{i}{2(n-s)}\) \cdot (n-s)}$.

\noi So here the error is of order $i$, which is significant if we view this as an isoperimetric-type result, since in this context one is typically interested in almost tight results. (With that, to the best of our knowledge, the bounds we obtain are new. In particular, for $i = 2$ we seem to obtain some new estimates related to the Kleitman-West problem. See the discussion in Section~\ref{subsubsec:results:edge-isop}.)
\exam

\xam
\label{xam:undetect}

\noi We will show (as a corollary of Theorem~\ref{thm:edge-isop}) that if $A$ is a binary code of length $n$ used over a binary symmetric channel, then the undetected error probability of $A$ is at most $O\(n^2\)$ times that of the union of two Hamming spheres of adjacent radii, whose size is roughly that of $A$. So here the error is of order $n^2$. However, in this type of coding estimates sub-exponential errors are ignored. Hence this result implies that unions of Hamming spheres (one can also take a Hamming sphere or a Hamming ball of an appropriate size) have, asymptotically, the largest  undetected error probability over the binary symmetric channel. (See (\ref{worst error exponent}) and the discussion preceding it, and also Section~\ref{subsubsec:8}.)
\exam

\myblt Finally, let us draw attention to the special case of the third of the claims above (it is also closely related to the second example above) in which $f$ is a characteristic function of a set. The claim then is that characteristic functions of Hamming spheres (or Hamming balls) are almost stablest with respect to noise among all sets of the same cardinality. To say this differently, consider the following  probabilistic experiment. Given a subset $A$ of $\H$, choose uniformly at random a point $x$ in $A$. Flip each coordinate in $x$ independently with probability $\e$ and check whether the obtained point is also in $A$. Then, the probability of this event is maximized (up to a sublinear in $n$ factor) if $A$ is a Hamming sphere (ball).

\noi Let us say a few words about the proofs, focusing on the proof of (\ref{ineq:moments}), since the  applications described above follow from it in a more or less standard manner. We prove (\ref{ineq:moments}) in Theorem~\ref{thm:norms} by a comparison argument, showing by induction on the dimension that for a homogeneous polynomial $f$ of degree $s$, the ratio $\frac{\E |f|^p}{\(\E f^2\)^{p/2}}$ cannot be much larger than that for the Krawchouk polynomial $K_s$. The error we obtain in this part of the argument is subexponential in the dimension. It is then reduced to a polynomial error by a tensorization argument (see Subsection~\ref{subsubsec:tensorization} below), applying the claim proved in the first step to tensor powers $f^{\otimes m}$ and passing to the limit as $m \rarrow \infty$. In this limit argument, the behavior of discrete objects such as Krawchouk polynomials is smoothened out, leading to a simply exponential expression in (\ref{ineq:moments}), and incurring a polynomial loss (see also the discussion above).

\noi A key element in controlling the growth of $\frac{\E |f|^p}{\(\E f^2\)^{p/2}}$ with dimension in the induction part of this argument is Hanner's inequality \cite{Lieb-Loss}: for $p \ge 2$ and for any two functions $g_0, g_1$ holds
\[
\|g_0+g_1\|_p^p + \|g_0-g_1\|_p^p ~\le~ \(\|g_0\|_p + \|g_1\|_p\)^p + \Big | \|g_0\|_p - \|g_1\|_p \Big |^p.
\]
An important part of our argument is showing the following intriguing fact: for any fixed $p \ge 2$ and for  sufficiently large $n$ and $s$, Krawchouk polynomials $K_{s-1}$ and $K_s$ on $\H$ satisfy Hanner's inequality almost with equality. To show this we rely on many known properties of Krawchouk polynomials (see Section~\ref{subsec:krawchouk}) and also prove some seemingly new ones: In particular, we provide a rather tight estimate for the $\ell_p$ norms of Krawchouk polynomials; and show their $\ell_2$ norm to be attained with only polynomial loss between any two of their roots, and also before their first and after their last roots. An additional implication of our results is that the above mentioned bivariate functions provide an accurate description of the behavior of Krawchouk polynomials $K_s$ for any sufficiently large $s$ (even a large constant $s$). Previously this seems to have been known mostly for $s$ growing linearly with dimension $n$ (see also \cite{Krasikov:Q} where the behavior of the magnitude of $|K_s|$ was analyzed for any $s$).

\subsubsection*{Related work}

\myblt A special case of (\ref{ineq:moments}), for $p = 4$, was shown in \cite{KS}, where it was also conjectured that the Krawchouk polynomials actually attain the maximum for $\frac{\E f^4}{\(\E f^2\)^2}$ among all homogeneous polynomials of the same degree. This conjecture has been recently proved in \cite{Aaronson}, by a short and a very elegant argument (using compression). It seems possible to extend the argument in \cite{Aaronson} to work for any even integer $p$. However, since this argument is essentially combinatorial in nature, it is not immediately obvious how to extend it to general $p > 2$.

\myblt After completing this paper, we have learned \cite{P-private-19} that a generalization of Theorem~\ref{thm:NHC} and Corollary~\ref{cor:sphere-stable} was proved in a concurrent work \cite{P-concurrent}. More specifically \cite{P-concurrent} proves the conjecture of \cite{P-private} (see the discussion following Corollary~\ref{cor:sphere-stable} in Section~\ref{subsubsec:hypercontractive-ineq}).

\myblt It was shown in \cite{Bogdanov-Mossel} that characteristic functions of Hamming spheres (or Hamming balls) of cardinality $2^{n - \alpha(n)}$, where $\alpha(n)$ is a slowly growing function of $n$, are (almost) stablest with respect to noise among all sets of the same cardinality. In \cite{OPS} Hamming spheres (or Hamming balls) of any cardinality are shown to be nearly stablest if the noise is very small, and it is conjectured that this should hold for any level of noise.

\myblt The hypercontractive inequality (\ref{ineq:HC}) was used  in \cite{ACKL} to obtain bounds on the distance components and other parameters of binary codes. We follow the approach of \cite{ACKL} in deriving some of our results, such as the second of the claims above, but replacing (\ref{ineq:HC}) with a (stronger) inequality (\ref{ineq:NHC}). We remark that the idea of using (\ref{ineq:HC}) to study the distance distribution of binary codes was introduced already in \cite{KL}.

\subsubsection*{Organization of the paper}

\noi The remainder of this paper is organized as follows. We describe the relevant notions and provide some background in the next subsection. Our results are stated formally and discussed in Section~\ref{subsec:results}. Somewhat unfortunately, the statements of the results involve certain functions of two variables, which will be defined later on in Section~\ref{subsec:bivariate}. This is done in order not to interrupt the flow of presentation.

\noi We define several bivariate functions which play an important role in this paper and describe their pertinent properties in Section~\ref{subsec:bivariate}. Some properties of Krawchouk polynomials and Hamming spheres are described in Sections~\ref{subsec:krawchouk}~-~\ref{subsec:hamming}. These subsections also clarify the relevance of some of the bivariate functions defined in Section~\ref{subsec:bivariate}, by showing them to describe aspects of behavior of Krawchouk polynomials or of Hamming spheres.

\noi Theorems~\ref{thm:tails}~to~\ref{thm:max-proj} and some related results are derived from Theorem~\ref{thm:norms} in Section~\ref{sec:thms}. Theorem~\ref{thm:norms} itself is proved in Section~\ref{sec:norms}. This paper contains many auxiliary results describing the behaviour of various univariate and bivariate functions. The proofs of these results are relegated to the Appendix.

\noi Let us suggest that (most of) Section~\ref{sec:technical} and the Appendix are better viewed as reference sections, written as laundry lists of results, and suitable for lookup, rather than for continuous reading.

\subsection{Background, definitions, and notation}

\noi We view $\H$ as a metric space, with the Hamming distance between $x, y \in \H$ given by $|x - y| = |\{i: x_i \not = y_i\}|$. The {\it Hamming sphere} of radius $r$ centered at $x$ is the set $S(x,r) = \left\{y \in \H:~|x-y| = r\right\}$. The {\it Hamming ball} of radius $r$ centered at $x$ is the set $B(x,r) = \left\{y \in \H:~|x-y| \le r\right\}$. Clearly, for any $x \in \H$ and $0 \le r \le n$ holds $|S(x,r)| = {n \choose r}$ and $|B(x,r) = \sum_{k=0}^r {n \choose k}$.

\noi Let $H(t) = t \log_2\(\frac 1t\) + (1-t) \log_2\(\frac{1}{1-t}\)$ be the binary entropy function. We will make repeated use of the following sequence of estimates (the first estimate follows from the Stirling formula, for the second estimate see e.g., Theorem~1.4.5. in \cite{van Lint}): For $x \in \H$ and $0 < r \le \frac n2$ holds
\beqn
\label{binomial-H}
\Theta\(\sqrt{\frac{n}{r(n-r)}}\) \cdot 2^{H\(\frac{r}{n}\) \cdot n} ~=~  |S(x,r)| ~\le~ |B(x,r)| ~\le~ 2^{H\(\frac{r}{n}\) \cdot n}.
\eeqn

\noi The asymptotic notation will always hide absolute constants (unless specifically stated otherwise).

\subsubsection{Distance distribution, edge-isoperimetry, binary codes}
\label{subsubsec:distance distribution}

\noi The distance distribution of a subset $A$ of $\H$ is given by the vector $\(a_0,a_1,...a_n\)$ with $a_i = |\{(x,y) \in A \times A, |x-y| = i\}|$. That is, $a_i = a_i(A)$ counts the pairs of points at distance $i$ in $A$. The distance distribution captures many important properties of a subset.

\noi {\it Edge-Isoperimetry}. For $1 \le i \le n$, let $G_i$ be the graph with vertices indexed by $\H$, in which two vertices are connected by an edge iff the Hamming distance between them is $i$. In particular, $G_1$ is the usual graph of the boolean cube. The {\it edge-isoperimetric problem} (see \cite{Bezrukov} for a survey on discrete isoperimetry) in a graph $G$ asks for a subset of vertices of a given cardinality, such that the number of edges crossing from this subset to its complement is as small as possible. If $G$ is regular, this is the same as maximizing the number of edges in an induced subgraph of $G$ with a given number of vertices. Note that a subset $A$ of vertices of $G_i$, this number is given by $a_i(A)$. The edge-isoperimetric problem has been resolved for $i=1$ \cite{Harper, Hart}, in which case the solution to the problem is the initial segment of the lexicographic ordering on the cube. To the best of our knowledge, the problem is still open for any $i > 1$.

\noi {\it Undetected error probability}. A {\it binary symmetric channel} (see e.g., \cite{Cover-Thomas}), with crossover probability $0 \le \e \le 1/2$ is a communication channel which on input $x \in \H$ outputs a random vector $y \in \H$ obtained by flipping each bit of $x$ independently, with probability $\e$. Given a binary code $C \subseteq \H$, the undetected error probability \cite{Klove-Korzhik} of $C$ is the average probability (over the codewords) that a codeword transmitted over a binary symmetric channel is distorted in such a way that the received word, though different from the transmitted one, also belongs to the code. It is easy to see that this can be expressed in terms of the distance distribution of $C$:
\[
P_{\mathrm{ue}}(C,\e) ~=~ \frac{1}{|C|} \cdot \sum_{i=1}^n a_i(C) \e^i (1-\e)^{n-i}.
\]
\noi The {\it worst asymptotic undetected error exponent} for codes of rate $0 \le R \le 1$ and crossover probability $\e$ was defined in \cite{ACKL} as
\[
P_{\mathrm{ue}}(R,\e) ~=~ \mathrm{lim sup}_{n \rarrow \infty} \(\frac 1n \max_C \log_2\(P_{\mathrm{ue}}(C,\e)\)\),
\]
where the maximum is taken over all codes $C \subseteq \H$ of cardinality at most $2^{Rn}$.

\noi {\it Binary error-correcting codes}. A binary error-correcting code $C$ of length $n$ and minimal distance $d$ is a subset of $\H$ such that the Hamming distance between any two distinct points in $C$ is at least $d$. This is clearly equivalent to $a_1= ... = a_{d-1} = 0$. The problem of finding the largest possible code with a given minimal distance is open. In \cite{dels} a family of linear inequalities holding for the distance distribution vector of any binary code were obtained. These inequalities play a key role in the linear programming relaxation of this problem \cite{dels}, which led to the best known upper bounds \cite{MRRW} on the cardinality of a code with a given minimal distance.

\noi The {\it largest asymptotic distance component rate} (see e.g., \cite{ABL, ACKL}) of a code with given rate and minimal distance is defined for $0 \le \mu, \delta \le \frac12$ and $0 \le R \le 1$ as
\[
b_{\mu}(R, \delta) ~=~ \mathrm{lim sup}_{n \rarrow \infty} \(\frac 1n \max_{C} \log_2\(a_{\lfloor \mu n \rfloor} (C)\)\),
\]
where the maximum is taken over all codes $C \subseteq \H$ of cardinality at most $2^{Rn}$ and minimal distance at least $\delta n$.

\subsubsection{Fourier analysis, polynomials, noise operators, and spectral projections}
\label{subsubsec:Fourier}

\noi We recall some basic notions in Fourier analysis on the boolean cube (see \cite{O'Donnel}). For $\alpha \in \H$, define the Walsh-Fourier character $W_{\alpha}$ on $\H$ by setting $W_{\alpha}(y) = (-1)^{\sum \alpha_i y_i}$, for all $y \in \H$. The {\it weight} of the character $W_{\alpha}$ is the Hamming weight $|\alpha|$ of $\alpha$.  The characters $\{W_{\alpha}\}_{\alpha \in \H}$ form an orthonormal basis in the space of real-valued functions on $\H$, under the inner product $\<f, g\> = \frac{1}{2^n} \sum_{x \in \H} f(x) g(x)$. The expansion $f = \sum_{\alpha \in \H} \widehat{f}(\alpha) W_{\alpha}$ defines the Fourier transform $\widehat{f}$ of $f$. We also have the Parseval identity, $\|f\|^2_2 = \sum_{\alpha \in \H} {\widehat{f}}^2(\alpha)$. One additional simple fact will be used several times in this paper: Let $g(x) = (-1)^{|x|} f(x)$. Then, writing $\bar{\alpha} = \alpha \oplus 1$ for the complement of a vector $\alpha \in \H$, for all $\alpha \in \H$ holds $\widehat{f}(\alpha) = \widehat{g}\(\bar{\alpha}\)$.

\noi {\it Polynomials on $\H$}. A function $f$ on $\H$ is a called a polynomial of degree $s$, for some $0 \le s \le n$, if $f$ belongs to the span of Walsh-Fourier characters of weight at most $s$. Alternatively, for $1 \le i \le n$, let $r_i$ be the Walsh-Fourier character of weight $1$, corresponding to $\alpha = \{i\}$. The functions $r_i$ are known as the {\it Rademacher functions} on $\H$. Then $f$ is a polynomial of degree $s$ if and only if $f$ is a multilinear polynomial of degree $s$ in $r_1,...,r_n$.

\noi A function $f$ is a homogeneous polynomial of degree $s$ on $\H$ if $f$ is a homogeneous multilinear polynomial of degree $s$ in $r_1,...,r_n$. Such functions are also called {\it Rademacher chaos of order $s$}. Note that if $f$ is a homogeneous polynomial of degree $s$ and $g(x) = (-1)^{|x|} f(x)$, then $g$ is a homogeneous polynomial of degree $n-s$. In particular, the spaces of homogeneous polynomials of degrees $s$ and $n-s$ are isometric for any $\ell_p$ norm on $\H$.

\noi {\it Krawchouk polynomials}. For $0 \le s \le n$, let $F_s$ be the sum of all Walsh-Fourier characters of weight $s$, that is $F_s = \sum_{|\alpha| = s} W_{\alpha}$. Note that $F_s$ is the Fourier transform of $2^n \cdot 1_S$, where $S$ is the Hamming sphere of radius $s$ around $0$. It is easy to see that $F_s(x)$ depends only on the Hamming weight $|x|$ of $x$, and it can be viewed as a univariate function on the integer points $0,...,n$, given by the restriction to $\{0,...,n\}$ of the univariate polynomial $K_s = \sum_{k=0}^s (-1)^k {x \choose k} {{n-x} \choose {s-k}}$ of degree $s$. That is, $F_s(x) = K_s(|x|)$. The polynomial $K_s$ is the $s^{th}$ {\it Krawchouk polynomial}. Abusing notation, we will also call $F_s$ the $s^{th}$ Krawchouk polynomial, and write $K_s$ for $F_s$ when the context is clear.

\noi {\it Spectral projections}. For $0 \le k \le n$ we define $\Pi_k$ to be the orthogonal projection to the subspace spanned by Walsh-Fourier characters of weight $k$. (This is the eigenspace of the Laplacian of the discrete cube corresponding to eigenvaleue $2k$.) That is, for a function $f$ on $\H$, and $0 \le k \le n$, we have $\Pi_k f = \sum_{|\alpha| = k} \widehat{f}(\alpha) W_{\alpha}$. We will also write $f_k$ for $\Pi_k f$ for ease of notation.

\noi {\it The noise operator}. Given a noise parameter $0 \le \e \le 1/2$, the noise operator $T_{\e}$ is a linear operator acting on functions on the boolean cube as follows: for $f:\H \rarrow \R$, $T_{\e} f$ at a point $x$ is the expected value of $f$ at $y$, where $y$ is "$(1-\e)$-correlated" with $x$. That is, $y$ is a random binary vector whose $i^{\small{th}}$ coordinate is $x_i$ with probability $1-\e$ and $1 - x_i$ with probability $\e$, independently for different coordinates. In other words, $\(T_{\e} f\)(x) = \E_y f(y)$, where $y$ is the output of the binary symmetric channel on input $x$. Writing this out explicitly, we have $\(T_{\e} f\)(x) =  \sum_{y \in \H} \e^{|y - x|}  (1-\e)^{n - |y-x|}  f(y)$. The noise operators form a semigroup: $T_{\e_1} T_{\e_2} = T_{\e_1 + \e_2 - 2 \e_1 \e_2}$. We will also write $f_{\e}$ for $T_{\e} f$, for brevity.
It is easy to see that $\widehat{f_{\e}}(\alpha) = (1-2\e)^{|\alpha|} \widehat{f}(\alpha)$, which means that $T_{\e} = \sum_{k=0}^n (1-2\e)^k \Pi_k$.

\subsubsection{Hypercontractive inequalities}
\label{subsubsec:hypercontractive}

\noi Hypercontractive inequalities \cite{Bonami,Gross,Beckner} form a family of analytic inequalities for functions on $\H$, with many applications in discrete mathematics, information theory, and theoretical computer science, see e.g., \cite{KKL, MOO, KKMPSU}, and also the monograph \cite{O'Donnel} and the references there.

\noi Let the $\ell_p$ norm of a real-valued function $f$ on $\H$ be given by $\|f\|_p = \(\frac{1}{2^n} \sum_{x \in \H} |f(x)|^p\)^{1/p}$. Hypercontractive inequalities  assert that applying noise to a function flattens it in a well-defined sense: a higher norm of the noisy function is upperbounded by a lower norm of the original function. A useful special case is the one involving the $\ell_2$ norm, (since this norm is easy to work with in applications). The inequality in this case is
\beqn
\label{ineq:HC}
\|f_{\e} \|_2 \le \|f\|_{1 + (1-2\e)^2}.
\eeqn

\noi It is easy to see that if $f$ is a characteristic function of a subset $A \subseteq \H$, then $\|f_{\e} \|^2_2 = \frac{1}{2^n} \sum_{i=0}^n a_i(A) \delta^i (1-\delta)^{n-i}$, for $\delta = 2\e(1-\e)$. The relevance of (\ref{ineq:HC}) to the study of distance distributions of binary codes has been pointed out in \cite{KL}. In \cite{ACKL} this inequality was used to obtain new bounds on the distance distribution, the undetected error probability, and other related parameters of binary codes of a given cardinality and minimal distance.

\noi {\it Stronger hypercontractive inequalities for highly concentrated functions.} While (\ref{ineq:HC}) is known to be essentially tight for functions which are almost flat to begin with, stronger hypercontractive inequalities were proved in \cite{PS} for functions $f$ on $\H$ for which the ratio $\frac{\|f\|_p}{\|f\|_1}$, for some $p > 1$, is exponentially large in $n$.

\noi {\it An uncertainty theorem}. Strong hypercontractive inequalities for highly concentrated functions were used in \cite{PS} to obtain a tight uncertainty-type result for $\H$. Let a non-zero function $f$ be supported on a set $A \subseteq \H$, with cardinality of $A$ being at most that of a Hamming ball of radius $\rho_1 n$ (for some $0 < \rho_1 < \frac12$). In fact, it suffices to assume, more generally, that the ratio $\frac{\|f\|_2}{\|f\|_1}$ is lower-bounded by $2^{\frac{1 - H\(\rho_1\)}{2} \cdot n}$. Then $\widehat{f}$ attains only an exponentially small fraction of its $\ell_2$ norm on any Hamming ball of radius $\rho_2 n$, provided $\rho_2 < \frac12 - \sqrt{\rho_1\(1-\rho_1\)}$.

\subsubsection{Tensorization}
\label{subsubsec:tensorization}

\noi We describe a useful and well-known tool in analysis which will be used several times in this paper. Let $f$ be a function on $\H$. For an integer $m \ge 1$, the $m^{\mathrm{th}}$ tensor power $F_m := f^{\otimes m}$ is a function on $nm$ boolean variables defined for $x_1,...,x_m \in \H$ by $F_m\(x_1,...,x_m\) = \prod_{i=1}^m f\(x_i\)$. We recall some useful properties of tensor powers:

\noi -- \,For any $\alpha_1,...,\alpha_m \in \H$ holds $\widehat{F_m}\(\alpha_1,...,\alpha_m\) = \prod_{i=1}^m \widehat{f}\(\alpha_i\)$. In particular, if $f$ is a homogeneous polynomial of degree $s$, then $F_m$ is a homogeneous polynomial of degree $sm$; and if $f$ is a (not necessarily homogeneous) polynomial of degree $s$, then $F_m$ is a polynomial of degree $sm$;

\noi -- \,For any $q \in \R$ holds $\E |F_m|^q = \(\E |f|^q\)^m$.

\subsection{Our results}
\label{subsec:results}

\subsubsection{Upper bounds for moments of polynomials}
\label{subsec:moments}
\noi We show that for any $p \ge 2$ and for any  $1 \le s \le \frac n2$, the $s^{\tiny th}$ Krawchouk polynomial $K_s = \sum_{|\alpha| = s} W_{\alpha}$ attains, within a relatively small error, the maximal ratio of $\frac{\|f\|_p}{\|f\|_2}$ among all homogeneous polynomials of degree $s$. Let $\psi(p,x)$ be the function defined in Subsection~\ref{subsubsec:ref:psi}.

\thm
\label{thm:norms}
For any $p \ge 2$, $0 \le s \le \frac n2$, and for any homogeneous polynomial $f$ of degree $s$ on $\H$ holds
\beqn
\label{ineq:norms}
\frac{\E |f|^p}{\(\E f^2\)^{\frac p2}} ~\le~ 2^{\psi\(p,\frac sn\) \cdot n}.
\eeqn

\noi There is an absolute constant $C > 0$ such that for any $p \ge 2$ and $0 \le s \le \frac n2$ holds
\[
2^{\psi\(p,\frac sn\) \cdot n} ~\le~ n \cdot C^p \cdot s^{\frac p4} \cdot \frac{\E |K_s|^p}{\(\E K_s^2\)^{\frac p2}}.
\]
\ethm

\noi {\bf Discussion}.

\noi --\, The assumption that $f$ is homogeneous is not necessary. In fact, we have, as a simple corollary of (\ref{ineq:norms}):
\cor
\label{cor:also balls}
The upper bound (\ref{ineq:norms}) holds for general polynomials of degree {\it as most} $s$.
\ecor

\myblt As mentioned above, if $p$ is an even integer, it seems possible to extend the argument given in \cite{Aaronson} for $p=4$ and to show that Krawchouk polynomials actually attain the maximum for $\frac{\E |f|^p}{\(\E f^2\)^{\frac p2}}$ among all homogeneous polynomials of the same degree.

\noi --\, The inequality (\ref{ineq:norms}) is a {\it Khintchine-type inequality}. Recall that Khintchine-type inequalities establish an upper bound on the ratio of two $\ell_p$ norms for functions coming from a certain restricted domain, typically a space of multivariate polynomials of a specified degree over a given product space. In particular, the prototypical Khintchine inequality \cite{Khintchine} states that the ratio of $\ell_2$ and $\ell_1$ norms of linear polynomials over the boolean cube $\H$ is bounded by an absolute constant. See \cite{IT} for a recent discussion and references. Viewed in this context, Theorem~\ref{thm:norms} states that for any $p > 2$ the ``Khintchine ratio" $\frac{\|f\|_p}{\|f\|_2}$ for polynomials of a given degree on the boolean cube is maximized, up to a small error, on the Krawchouk polynomial of this degree.

\noi --\, It is easy to see that Theorem~\ref{thm:norms} essentially determines the $\|\cdot\|_{2 \rarrow p}$ norm of the spectral projection operator $\Pi_k$ (see \cite{P2} where the norms of these operators are investigated).

\subsubsection{Tail bounds for polynomials}

\noi We show that Krawchouk polynomials have (almost) the heaviest tails among all polynomials of same degree and $\ell_2$ norm. Let $\tau$ be the function defined in Subsection~\ref{subsubsec:ref:tau}. Let $H$ be the binary entropy function.

\noi To make the statement of the second part of the following theorem more legible, recall (see Section~\ref{subsec:krawchouk}) that the Krawchouk polynomial $K_s$ on $\H$ has all its roots in the interval \\ $\left[\frac n2 - \sqrt{s(n-s)}, \frac n2 +\sqrt{s(n-s)}\right]$, and that the distance between any two consecutive roots is $o(n)$.

\thm
\label{thm:tails}

\noi Let $f$ be a polynomial of degree $s \le \frac n2$ on $\H$. Then for all $0 \le i \le \frac n2$ holds
\beqn
\label{ineq:tails}
\mathrm{Pr}\left\{|f| \ge \|f\|_2 \cdot 2^{\(\tau\(\frac sn, \frac in\) - \frac12 H\(\frac sn\)\) \cdot n}\right\} ~\le~ 2^{\(H\(\frac in\) - 1\) \cdot n}.
\eeqn

\noi Moreover, for $f = K_s$ we have:

\begin{itemize}

\item

\noi For any $0 \le i \le \frac n2 - \sqrt{s(n-s)}$ holds
\[
\mathrm{Pr}\left\{|K_s| \ge \sqrt{\frac{{n \choose s}}{2^{H\(\frac sn\) \cdot n}}} \cdot \|K_s\|_2 \cdot 2^{\(\tau\(\frac sn, \frac in\) - \frac12 H\(\frac sn\)\) \cdot n}\right\} ~\ge~ \Omega\(\frac{1}{\sqrt{i}}\) \cdot 2^{\(H\(\frac in\) - 1\) \cdot n}.
\]

\item

\noi Between any two consecutive roots of $K_s$ there is a point $i$ for which
\[
\mathrm{Pr}\left\{|K_s| \ge \Omega\(\frac{1}{n^{5/2}}\) \cdot \|K_s\|_2 \cdot  2^{\(\tau\(\frac sn, \frac in\) - \frac12 H\(\frac sn\)\) \cdot n}\right\} ~\ge~ \Omega\(\frac{1}{\sqrt{i}}\) \cdot 2^{\(H\(\frac in\) - 1\) \cdot n}.
\]

\end{itemize}
\ethm

\noi {\bf Discussion}

\noi --\, Note that, by (\ref{binomial-H}), the correction factor $\sqrt{\frac{{n \choose s}}{2^{H\(\frac sn\) \cdot n}}}$ is $\Theta\(s^{\frac 14}\)$.

\myblt A polynomial $f$ of degree $s$ is a linear combination of $\sum_{k = 0}^s {n \choose k}$ Walsh-Fourier characters, which are orthonormal, and all of which evaluate to $1$ at $0$. It is easy to see that this implies that $\frac{\|f\|_{\infty}}{\|f\|_2} \le \sqrt{\sum_{k = 0}^s {n \choose k}} \le 2^{\frac12 H\(\frac sn\) \cdot n}$, with equality attained for $f = \sum_{|\alpha| \le s} W_s = \sum_{k=0}^s K_k$. On the other hand, the function $\tau\(\frac sn, \frac in\) - \frac12 H\(\frac sn\)$ decreases from $\frac12 H\(\frac sn\)$ to $0$ as $i$ does from $0$ to $\frac n2$. Hence (\ref{ineq:tails}) provides tail estimates for the whole range of values of $|f|$.

\myblt The bound (\ref{ineq:tails}) is a pointwise improvement over the estimate $\mathrm{Pr}\{|f| \ge \|f\|_2 \cdot t\} \le e^2 \cdot e^{-\(\frac te\)^{\frac 1k}}$, due to \cite{BKS, Talagrand-inc-pos}. To see this, note that the latter bound was obtained by applying the inequality $\mathrm{Pr}\{|f| \ge T\} \le \min_{p \ge 2} \left\{\frac{\E |f|^p}{T^p}\right\}$, using (\ref{ineq:moments-HC}) to bound the RHS, and choosing a suitable value of $p$. The proof of (\ref{ineq:tails}) uses the same approach, while replacing (\ref{ineq:moments-HC}) with a stronger bound (\ref{ineq:moments}) and choosing the {\it optimal} value of $p$.

\myblt There is a gap between the upper bound (\ref{ineq:tails}) and the lower bound provided by Krawchouk polynomials, due to the correction factor $\sqrt{\frac{{n \choose s}}{2^{H\(\frac sn\) \cdot n}}}$ in the value of the threshold. One can ask how accurate (\ref{ineq:tails}) is for Krawchouk polynomials, if there is no correction factor. It's not hard to see, using the properties of the function $\tau$, that e.g., for $i$ growing linearly with $n$, that is for $t$ bounded away from $\frac{\|K_s\|_{\infty}}{\|K_s\|_2}$, the error in the estimate of (\ref{ineq:tails}) for the probability of $|K_s| \ge \|K_s\|_2 \cdot t$ is of order at most $2^{O\(\frac{\log(s)}{s}\) \cdot n}$. In particular, if $s$ is an increasing function of $n$, (\ref{ineq:tails}) provides the right constant in the large deviation inequalities for homogeneous polynomials of degree $s$ (that is, Rademacher chaos of order $s$).

\subsubsection{An isoperimetric-type inequality}
\label{subsubsec:results:edge-isop}

\noi Recall that $G_i$ is a graph with vertices indexed by $\H$, in which two vertices are connected by an edge iff the Hamming distance between them is $i$. We prove an edge-isoperimetric inequality for the graphs $\{G_i\}_i$, and show this inequality to be somewhat tight for the Hamming sphere, or a union of two adjacent spheres, depending on the parity of $i$.

\noi Let $A \subseteq \H$. Recall that $a_i(A) = |\{(x,y) \in A \times A, |x-y| = i\}|$ is the $i^{\tiny{th}}$ distance component of $A$.

\thm
\label{thm:edge-isop}

\noi Let $A \subseteq \H$, with $|A| \le 2^{H(\sigma) n}$, for some $0 \le \sigma \le \frac12$. Then for $1 \le i \le 2\sigma(1-\sigma) n$ holds
\beqn
\label{ineq:edge-isop}
a_i(A) \quad \le \quad |A| \cdot 2^{\(\sigma H\(\frac{i}{2\sigma n}\) + (1-\sigma)H\(\frac{i}{2(n-\sigma n)}\)\) \cdot n}.
\eeqn

\noi Let $s = \sigma n$, and assume $s$ to be integer. If $i$ is even, this inequality is tight, up to a factor of $O(i)$, if $A$ a Hamming sphere of radius $s$. For an arbitrary $i$, this is tight, up to a factor of $O\(\sqrt{\frac{n-s}{s}} \cdot i\)$, if $A$ is the union of two adjacent spheres of dimension $n-1$ and radii $s-1$ and $s$.

\ethm

\noi {\bf Discussion}

\myblt For $2\sigma(1-\sigma) n \le i \le \frac n2$, the distance component $a_i(A)$ could be (essentially) as large as $|A|^2$. For $i > \frac n2$, the bounds on $a_i$ reduce to these on $a_{n-i}$. See Remark~\ref{rem:edge-isop}.

\myblt The proof of (\ref{ineq:edge-isop}) follows the argument in \cite{ACKL}, replacing the hypercontractive inequality (\ref{ineq:HC}) used in \cite{ACKL} with a stronger inequality proved in Corollary~\ref{cor:hypercontractive} (which is a special case of (\ref{ineq:NHC})).

\myblt Choosing $i=2$ in (\ref{ineq:edge-isop}), we get the following claim:
\cor
\label{cor:Kleitman-West}
For $A \subseteq \H$ with $|A| \le {n \choose s}$ holds $a_2(A) \le e^2 s(n-s) \cdot |A|$.
\ecor

\noi This is tight, within a factor of $e^2$, if $A$ is a Hamming sphere of radius $s$.

\noi Let us also consider this bound in the context of the Kleitman-West problem (see e.g., \cite{Harper-KW}). This is the edge-isoperimetric problem for the Hamming sphere (see Section~\ref{subsubsec:distance distribution}). One way to pose this problem is as follows. Given the dimension $n$ and the radius $0 \le r \le \frac n2$ of the sphere $S = S(0,r) \subseteq \H$, determine how large can $a_2(A)$ be for a subset $A$ of $S$ of a given size.

\noi Let $s \le r$ and let $|A| = {n \choose s}$. Let $A$ be an $\big(n-(r-s)\big)$-dimensional Hamming sphere of radius $s$ embedded in $S$, by concatenating $r-s$ coordinates to points in $A$ and setting them to be $1$. Then $a_2(A) = s(n-r) \cdot |A|$, so the bound in Corollary~\ref{cor:Kleitman-West} is tight for $A$ within a factor of $2e^2$. For $s < \frac{r}{C}$, where $C$ is a sufficiently large constant, this bound seems to improve the best known upper bounds on $a_2(A)$ for $A \subseteq S$ (which, to the best of our knowledge, come from the logarithmic Sobolev inequality for the Hamming sphere \cite{Lee-Yau}).

\myblt Recall (see Section~\ref{subsubsec:distance distribution}) that an undetected error probability of a binary code $C \subseteq \H$ is given by $P_{\mathrm{ue}}(C,\e) = \frac{1}{|C|} \cdot \sum_{i=1}^n a_i(C) \e^i (1-\e)^{n-i}$. Theorem~\ref{thm:edge-isop} implies (as shown in Section~\ref{subsubsec:8}) that a union of two adjacent spheres of dimension $n-1$ and radii $s-1$ and $s$ maximizes the undetected error probability for all codes of cardinality at most $2^{H\(\frac sn\) \cdot n}$, up to at most a polynomial in $n$ factor. A simple consequence of this fact is the following expression for the worst asymptotic undetected error exponent: For all $0 < R \le 1$ and $0 < \e \le \frac12$ holds
\beqn
\label{worst error exponent}
P_{\mathrm{ue}}(R,\e) ~=~ \sigma H\(\frac{x}{\sigma}\) + (1-\sigma) H\(\frac{x}{1-\sigma}\) + 2x\log_2(\e) + (1-2x) \log_2(1-\e),
\eeqn
where $\sigma = H^{-1}(R)$ and $x = x(\sigma, \e) = \frac{-\e^2 + \e \sqrt{\e^2 + 4(1-2\e) \sigma(1-\sigma)}}{2(1-2\e)}$.

\subsubsection{A hypercontractive inequality}
\label{subsubsec:hypercontractive-ineq}

\noi We prove a nearly tight hypercontractive inequality for functions on $\H$, which takes into account the distribution of a function, specifically the ratio between its $\ell_p$ and $\ell_1$ norms. (See \cite{PS} for a different family of hypercontractive inequalities taking into account the ratio between $\ell_p$ and $\ell_1$ norms of a function.)

\noi Let $\eta$ be the function defined in Section~\ref{subsubsec:ref:eta}. Recall that $\eta(x,\e) < 0$ for all $0 < \e < \frac12$ and $0 < x \le \frac{(1-2\e)^2}{1 + (1-2\e)^2}$. Moreover, $\eta$ is concave and decreasing in $x$ for any $0 \le \e \le \frac12$.

\noi For a function $f$ on $\H$ and $1 \le p \le \infty$ let $r(p) = r_f(p) = \frac{1}{n} \log_2\(\frac{\|f\|_p}{\|f\|_1}\)$. Note that $0 \le r(p) \le \frac{p-1}{p}$.

\thm
\label{thm:NHC}

\noi Let $f$ be a function on $\H$, and let $0 \le \e \le 1/2$. Then
\beqn
\label{ineq:NHC}
\|f_{\e}\|_2 ~~\le~~ 2^{\eta\(r\(1 + (1-2\e)^2\),~\e\) \cdot n} \cdot \|f\|_{1 + (1-2\e)^2}.
\eeqn

\noi This is tight up to a factor of $O\(s^{3/4}\)$ if $f$ is proportional to a characteristic function of a  Hamming sphere of radius $s$.
\ethm

\noi {\bf Discussion}.

\myblt (\ref{ineq:NHC}) is a strengthening of the hypercontractive inequality (\ref{ineq:HC}). However, since (as is easily seen) $\frac{\partial \eta(x,\e)}{\partial x}_{|x = 0} = 0$, this improvement is significant only if $\frac{\|f\|_{1 + (1-2\e)^2}}{\|f\|_1} \ge 2^{\Omega\(\sqrt{n}\)}$.

\myblt We have the following corollary of (\ref{ineq:NHC}), extending it to other $\ell_p$ norms. For $p \ge 1$, let $\eta_p$ be the function defined in Section\ref{subsubsec:ref:eta}.
\cor
\label{cor:sphere-stable}
Let $f$ be a function on $\H$, let $0 \le \e \le 1/2$ and let $p \ge  1 + (1-2\e)^2$. Then
\[
\|f_{\e}\|_2 ~~\le~~ 2^{\eta_p\(r(p),\e\)\cdot n} \cdot \|f\|_p.
\]
\ecor

\noi As will be seen in the proof of this result, the RHS of the inequality above increases in $p$, so it is in general weaker than (\ref{ineq:NHC}). However, it is still tight up to a factor of $O\(s^{3/4}\)$ if $f$ is proportional to the characteristic function of a  Hamming sphere of radius $s$. To see this, note that the RHS of this inequality does not depend on $p$, if $f$ is a characteristic function of a set. Hence, this corollary can be rephrased as follows: For any $0 \le \delta \le \frac12$ and $p \ge 2 - 2\delta$, a characteristic function of a  Hamming sphere of radius $s$ maximizes, within a factor of $O\(s^{3/4}\)$, the inner product $\<f_{\delta}, f\>$ among all functions with the same $\ell_1$ and $\ell_p$ norms\footnote{Recall that $\|f_{\e}\|^2_2 = \<f_{\delta},f\>$, for $\delta = 2\e(1-\e)$.}.

\myblt Let us mention a more general conjecture \cite{P-private}: For any $q > 1$ and a threshold value $t = t(q,\e)$ there exists $p_0 = p_0(q,\e,t) \le 1 + (1-2\e)^2 (q-1)$  such that for any $p \ge p_0$ the maximum of the ratio $\frac{\|f_{\e}\|_q}{\|f\|_p}$ over all functions $f$ on $\H$ with $r(p) \ge t$ is essentially attained at the characteristic function of a Hamming sphere of an appropriate radius.

\myblt Upper bounds on the asymptotic distance component rates $b_{\mu}(R,\delta)$ of binary codes with given rate and minimal distance (see Section~\ref{subsubsec:distance distribution}) were obtained in \cite{ACKL} using the hypercontractive inequality (\ref{ineq:HC}). These bounds can be improved by using (\ref{ineq:NHC}) instead of (\ref{ineq:HC}), similarly to the improvement obtained in (\ref{ineq:edge-isop}) over the bounds of \cite{ACKL} for $b_{\mu}(R,0)$. We do not go into details, since the bounds, both in \cite{ACKL} and here, are not explicit, but rather given, for each $R$ and $\delta$, as the minimal value of a certain explicit function in a constant number of variables (three in \cite{ACKL} and two in our case) over its domain.

\subsubsection{An uncertainty theorem}
\label{subsubsec:proj}

\noi We give an extension of an uncertainty-type result from \cite{PS} (see Section~\ref{subsubsec:hypercontractive}). Let $\pi$ be the function defined in Section~\ref{subsubsec:ref:pi}.
Recall that for a function $f$ on $\H$ and $0 \le k \le n$ we write $f_k$ for the orthogonal projection of $f$ on the space of Walsh-Fourier characters of weight $k$. We write $x \wedge y$ for the minimum of $x$ and $y$ and, as above, given a function $f$ on $\H$, write $r(p)$ for $\frac{1}{n} \log_2\(\frac{\|f\|_p}{\|f\|_1}\)$.

\thm
\label{thm:max-proj}

\noi Let $f$ be a function on $\H$. Then for any $p \ge 2$ and $0 \le k \le n$ holds
\beqn
\label{ineq:max-proj}
\|f_k\|_2 ~~\le~~ 2^{\(\pi\(\frac kn \wedge \frac{n-k}{n},H^{-1}\(1 - \frac{p}{p-1} \cdot r(p)\)\) - \frac{p-2}{2p - 2} \cdot r(p)\) \cdot n} \cdot \|f\|_p.
\eeqn

\noi Let $f$ be proportional to a characteristic function of a Hamming sphere of radius $s$. Then this inequality is tight up to a factor of $O\((ks)^{1/4}\)$ for $0 \le k \le \frac n2 - \sqrt{s(n-s)}$ and $\frac n2 +\sqrt{s(n-s)} \le k \le n$. In addition, between any two consecutive roots of the Krawchouk polynomial $K_s$ there is a point $k$ for which this is tight up to a factor of $O\(n^{5/2}\)$.

\ethm

\noi {\bf Discussion}.

\myblt An alternative (somewhat imprecise) way to phrase this result is as follows: For any $p \ge 2$, characteristic functions of Hamming spheres have (almost) the largest spectral projections among all functions with the same $\ell_1$ and $\ell_p$ norms.

\myblt Let $\frac{\|f\|_2}{\|f\|_1} = 2^{\frac{1 - H\(\rho\)}{2} \cdot n}$, for some $0 < \rho < \frac12$. Then the theorem with $p=2$ implies $\|f_k\|_2 \le 2^{\pi\(\rho, \frac kn\) \cdot n} \cdot \|f\|_2$. In particular, for $\frac{k}{n}$ bounded away from below from $\frac12 - \sqrt{\rho\(1-\rho\)}$, this implies that $\|f_k\|_2$ is exponentially smaller than $\|f\|_2$ (since $\pi\(\rho, \frac kn\)$ is negative in this range, see Lemma~\ref{lem:ref:pi}), recovering the result in \cite{PS}. Furthermore, we get a quantitative upper bound on the exponent of the ratio $\frac{\|f_k\|_2}{\|f\|_2}$ for $0 < \frac kn \le \frac12 - \sqrt{\rho\(1-\rho\)}$.

\section{Bivariate functions, Krawchouk polynomials, and Hamming spheres}
\label{sec:technical}

\subsection{Some bivariate functions}
\label{subsec:bivariate}

\noi Section~\ref{subsec:results} describes some functional inequalities on the Hamming cube. These inequalities involve certain functions of two variables. A good way to come to terms with these functions is to realize that they describe various aspects of the behavior of Hamming spheres or of Krawchouk polynomials (see Sections~\ref{subsec:krawchouk}~and~\ref{subsec:hamming}). In this subsection we define these functions and list their relevant properties.

\subsubsection{The function $I$}
\label{subsubsec:ref:I}

\noi For $0 \le x \le \frac12$ and $0 < y \le \frac12 - \sqrt{x(1-x)}$, let
\[
I(x,y) ~=~ \log_2(1-x) + \frac a2 \log_2\(1 - 2x - b\) + x \log_2\(\frac{a + b}{2(1-x)}\) - \frac12 \log_2\(2(1-x) - a^2 - ab\),
\]
where $a = 1 - 2y$ and $b = \sqrt{a^2 - 4x(1-x)}$. Extend this by continuity to $y = 0$ by setting $I(x,0) = -1$ for all $0 \le x \le \frac12$.

\noi Let $r(x,y) = \frac{\(1-2x\) + \sqrt{\(1-2x\)^2 - 4y(1-y)}}{2(1-y)}$ for For $0 \le x \le \frac12$ and $0 \le y \le \frac12 - \sqrt{x(1-x)}$. Then (\cite{KL}, with a correction in \cite{Kr-Lit}) $I$ is an indefinite integral of $r$, that is
$\int_{0}^y \log_2\(r(x,z)\) dz = I\(y, x\) - I\(0, x\)$.

\lem
\label{lem:ref:r}
For a fixed $0 \le x \le \frac12$, the function $r(x,y)$ decreases in $y$. In particular, for any $y \ge 0$ holds $r(x,y) \le r(x,0) = 1-2x$.
\elem

\cor
For a fixed $0 \le x \le \frac12$, the function $I(y,x)$ is decreasing and concave in $y$.
\ecor

\subsubsection{The function $\tau$}
\label{subsubsec:ref:tau}

\noi For $0 \le x, y \le \frac12$, let
\[
\tau(x,y) ~=~ \left\{\begin{array}{ccc}  H(x) + I(y,x) - I(0,x) & \mathrm{if} & y \le \frac12 - \sqrt{x(1-x)} \\ \frac{1 + H(x) - H(y)}{2} & \mathrm{otherwise} \end{array}\right.
\]
It is easy to verify that $\tau$ is continuous in both variables. Using the results in Subsection~\ref{subsubsec:ref:I}, we see that $\frac{\partial \tau(x,y)}{\partial y} = \left\{\begin{array}{ccc}  \log_2(r(x,y)) & \mathrm{if} & y < \frac12 - \sqrt{x(1-x)} \\ \frac12 \log_2\(\frac{y}{1-y}\) & \mathrm{if} & \frac12 - \sqrt{x(1-x)} < y < \frac12 \end{array}\right.$. This means that $\tau(x,y)$ is decreasing and concave in $y$ on $0 \le y \le \frac12 - \sqrt{x(1-x)}$, and is decreasing and convex in $y$ afterwards.

\lem
\label{lem:ref:tau}

\noi For all $0 \le x, y \le \frac12$ holds
\[
H(y) + \tau(x,y) ~=~ H(x) + \tau(y,x).
\]
\elem

\subsubsection{The function $h$}
\label{subsubsec:ref:h}

\noi For $2 \le p < \infty$ and $0 \le x \le \frac12$, let
\[
h(p,x) ~=~  x^{\frac{1}{p}} (1-x)^{\frac{p-1}{p}} + x^{\frac{p-1}{p}} (1-x)^{\frac{1}{p}}.
\]

\lem
\label{lem:ref:h}

\begin{enumerate}

\item For any $p \ge 2$ he function $h(p,x)$ increases from $0$ to $1$, as $x$ goes from $0$ to $\frac 12$.

\item For any $0 < x < \frac 12$,  the function $h(p,x)$ increases from $2\sqrt{x(1-x)}$ to $1$, as $p$ goes from $2$ to $\infty$.

\item Let $0 < x < \frac12$. If $p > 2$ and $h(p,y) = 1 - 2x$, then $y < \frac12 - \sqrt{x(1-x)}$.

\end{enumerate}
\elem

\subsubsection{The function $\psi$}
\label{subsubsec:ref:psi}

\noi For $2 \le p < \infty$ and $0 \le x \le \frac12$, let
\[
\psi(p,x) ~=~  H(y) - 1 + p\tau(x,y) - \frac p2 H(x),
\]
where $y$ is determined by $h(p,y) = 1 - 2x$.

\noi The function $\psi$ has another useful representation. First, we state an auxiliary lemma.
\lem
\label{lem:ref:psi-aux}
For any $p \ge 2$ the function $a(p, \delta) = \(\frac12 - \delta\) \cdot \frac{(1-\delta)^{p-1} - \delta^{p-1}}{(1-\delta)^p + \delta^p}$ decreases from $\frac12$ to $0$, as $\delta$ goes from $0$ to $\frac12$.
\elem

\pro
\label{pro:psi-second rep}
We also have
\[
\psi(p,x) ~=~ (p-1) + \log_2\Big((1-\delta)^p + \delta^p\Big) - \frac{p}{2} H(x) - px\log_2(1-2\delta),
\]
where $\delta$ is determined by $x = a(p,\delta)$. (Note that $\delta$ is well-defined by Lemma~\ref{lem:ref:psi-aux}.)
\epro

\pro
\label{pro:psi-concave convex}
\begin{enumerate}

\item

\noi For any $p \ge 2$ and for any $x \ge 0$ holds $\psi(p,0) = \psi(2,x) = 0$.

\item

\noi The function $\psi(p,x)$ is increasing and (strongly) convex in $p$ for any $x > 0$.

\item

\noi The function $\psi(p,x)$ is increasing and (strongly) concave in $x$ for any $p > 2$. We also have $\frac{\partial \psi(p,x)}{\partial x}_{| x = 0} = \frac{p \log_2(p-1)}{2}$.

\end{enumerate}

\epro

\subsubsection{The function $\pi$}
\label{subsubsec:ref:pi}

\noi For $0 \le x, y \le \frac12$, let
\[
\pi(x,y) ~~=~~ \left\{\begin{array}{ccc} I(y,x) - I(0,x) + \frac{H(x) + H(y) - 1}{2}  & \mathrm{if} & y \le \frac12 - \sqrt{x(1-x)} \\
0 & \mathrm{otherwise} \end{array} \right.
\]

\noi Note that $\pi(x,y) = \tau(x,y) - \frac{1 + H(x) - H(y)}{2}$. In particular, $\pi$ is continuous in both variables.

\lem
\label{lem:ref:pi}
The function $\pi$ is symmetric, that is $\pi(x,y) = \pi(y,x)$ for all $0 \le x,y \le \frac12$.
It is strictly negative for $y < \frac12 - \sqrt{x(1-x)}$ and increasing in both arguments.
\elem

\lem
\label{lem:ref:pi-min}
For any $0 \le \kappa, \sigma \le \frac12$ holds
\[
\pi(\sigma,\kappa) ~=~ \frac12 \min_{0 \le \delta \le \frac12} \left\{\sigma H\(\frac{x}{\sigma}\) + \(1-\sigma\) H\(\frac{x}{1-\sigma}\) + 2x\log_2(\delta) + (1-2x) \log_2(1-\delta) - \kappa \log_2(1-2\delta)\right\}   ,
\]
where $x = x(\sigma,\delta) = \frac{-\delta^2 + \delta \sqrt{\delta^2 + 4(1-2\delta) \sigma(1-\sigma)}}{2(1-2\delta)}$.
\elem

\subsubsection{The functions $\phi$ and $\tilde{\phi}$}
\label{subsubsec:ref:phi}

\noi For $0 \le \sigma, \e \le \frac12$, let
\[
\alpha_{\sigma,\e}(x) ~=~ \sigma H\(\frac{x}{\sigma}\) + (1-\sigma) H\(\frac{x}{1-\sigma}\) + 2x\log_2(\e) + (1-2x) \log_2(1-\e).
\]

\noi Let
\[
\phi(\sigma, \e) ~=~
H(\sigma) - 1 + \max_{0 \le x \le \sigma} \alpha_{\sigma,\e}(x).
\]
The value of $x$ for which the maximum is attained is $x = x(\sigma, \e) = \frac{-\e^2 + \e \sqrt{\e^2 + 4(1-2\e) \sigma(1-\sigma)}}{2(1-2\e)}$ (see e.g., \cite{ACKL})

\noi For $0 \le y \le 1$ and $0 \le \e \le \frac12$, let
\[
\tilde{\phi}(y,\e) ~=~ \phi\(H^{-1}(y),\e\).
\]

\noi We list some relevant properties of $\phi$.

\lem
\label{lem:ref:phi}

\noi For all $0 \le \sigma, \e \le \frac12$ holds
\[
\phi(\sigma, \e) ~=~ \max_{0 \le y \le \frac12} \Big\{y \log_2(1 - 2\e) + H(y) + 2\tau(\sigma,y) \Big\} - 2.
\]
\elem

\lem
\label{lem:ref:phi:edge-isop}
Let $0 \le \sigma \le \frac12$ and let $0 \le y \le 2\sigma(1-\sigma)$. Then
\[
\min_{0 < \e \le \frac12} \Big\{\phi\(\sigma,\e\) + 1 - H(\sigma) - y \log_2(\e) - \(1 - y\) \log_2(1-\e) \Big\}  =  \sigma H\(\frac{y}{2\sigma}\) + (1-\sigma) H\(\frac{y}{2(1 - \sigma)}\).
\]
\elem

\lem
\label{lem:ref:phi:disc-cont}
Let $\sigma = \frac sn$, $s$ an integer between $0$ and $\frac n2$. Let $0 \le \e \le \frac12$. Let $A = \max_{0 \le x \le \sigma} \alpha_{\sigma,\e}(x)$ and $B = \max_{0 \le i \le s} \alpha_{\sigma,\e}\(\frac in\)$. Then $B \le A  \le B + O\(\frac 1n\)$, where the constant in the asymptotic notation is absolute.
\elem

\noi We also list some relevant properties of $\tilde{\phi}$.

\lem
\label{lem:ref:tilde-phi}

\noi The function $\tilde{\phi}(y,\e)$ is (strictly) increasing and concave in $y$, for any fixed $0 \le \e \le \frac12$. Moreover, $\tilde{\phi}(1, \e) = 0$, and
the one-sided derivatives with respect to $y$ of $\tilde{\phi}$ at the endpoints of the interval are $\frac{\partial \tilde{\phi}}{\partial y}(0,\e) = 2$ and $\frac{\partial \tilde{\phi}}{\partial y}(1,\e) = \frac{1}{1-\e}$.
\elem

\subsubsection{The functions $\eta$ and $\eta_p$}
\label{subsubsec:ref:eta}

\noi For $0 \le \e \le \frac12$, and for $0 \le x \le \frac{(1-2\e)^2}{1 + (1-2\e)^2}$, let
\[
\eta(x,\e) ~=~ \frac12 \tilde{\phi}\(1 - \frac{1 + (1-2\e)^2}{(1-2\e)^2} \cdot x, ~2\e(1-\e)\) + \frac{1}{(1-2\e)^2} \cdot x,
\]
where the function $\tilde{\phi}$ is defined in Section~\ref{subsubsec:ref:phi}.

\noi More generally, for $p > 1$, for $0 \le \e \le \frac12$, and for $0 \le x \le \frac{p-1}{p}$, let
\[
\eta_p(x,\e) ~~=~~ \frac12 \tilde{\phi}\(1 - \frac{p}{p-1} \cdot x, ~2\e(1-\e)\) + \frac{1}{p-1} \cdot x.
\]
Note that $\eta = \eta_{1 + (1-2\e)^2}$.

\lem
\label{lem:ref:eta-p}
The function $\eta_p(x,\e)$ is concave and decreasing in $x$ for fixed $p$ and $\e$ satisfying $p \ge 1 + (1 - 2\e)^2$. Moreover, if we also assume $0 < \e < \frac12$ then it is strictly negative for $0 < x \le \frac{p-1}{p}$.
\elem

\subsection{Krawchouk polynomials}
\label{subsec:krawchouk}

\noi Krawchouk polynomials were defined in Section~\ref{subsubsec:Fourier}. In this subsection we list some of their properties. We refer to \cite{Krasikov:Q, KZ, Lev, Lev-chapter} for many of the facts stated below. Some of the properties we describe, in particular Proposition~\ref{pro:norm-any-two-roots} and Corollary~\ref{cor:norm-any-two-roots} seem to be new and might be of independent interest.

\noi {\it Notation:} Here and below we will write $a \in b \pm \e$ as a shorthand for $b - \e \le a \le b + \e$.

\begin{enumerate}

\item {\it Value at $0$}. For all $0 \le s \le n$ holds $K_s(0) = {n \choose s}$.

\item {\it Symmetry}. For all $0 \le i, s \le n$ holds $K_s(i) = (-1)^s \cdot K_s(n-i)$.

\item {\it Reciprocity}. For all $0 \le i, s \le n$ holds ${n \choose i} K_s(i) = {n \choose s} K_i(s)$.

\item {\it $\ell_2$ norm}. Viewing $K_s$ as a function on $\H$ or, equivalently, as a univariate real polynomial, endowing $\R$ with the binomial measure $\mu(i) = \frac{{n \choose i}}{2^n}$, for $0 \le i \le n$, we have $\|K_s\|_2 = \sqrt{{n \choose s}}$.

\item {\it Roots}. The polynomial $K_s$ (viewed as a univariate polynomial) has $s$ distinct real roots, which lie in the interval $\frac{n}{2} \pm \sqrt{s(n-s)} $. For $1 \le s \le \frac n2$, the distance between any two consecutive roots is at least $2$ and at most $o(n)$.

\item {\it Magnitude outside the root region}. Let $I$ and $\tau$ be the functions defined in Section~\ref{subsec:bivariate}. As shown in \cite{Krasikov:Q} (a more precise version of a result in \cite{KL}), for any  $0 \le  s \le \frac n2$ and $0 \le i \le \frac n2 - \sqrt{s(n-s)}$ holds
\beqn
\label{Krawchouk-tail}
\frac{{n \choose s}}{2^{H\(\frac sn\) \cdot n}} \cdot 2^{\tau\(\frac sn, \frac in\) \cdot n}~=~ {n \choose s} \cdot 2^{\(I\(\frac in, \frac sn\) - I\(0,\frac sn\)\) \cdot n} ~\le~
K_s(i) ~\le~ 2^{\(\tau\(\frac sn, \frac in\) + o(1)\) \cdot n}.
\eeqn
In particular, using (\ref{binomial-H}),
\[
2^{\(\tau\(\frac sn, \frac in\) - o(1)\) \cdot n} ~\le~ K_s(i) ~\le~ 2^{\(\tau\(\frac sn, \frac in\) + o(1)\) \cdot n}.
\]

\item {\it Magnitude in the root region}.
By Corollary~\ref{cor:norm-any-two-roots}, the polynomial $K_s$, for $1 \le s \le \frac n2$, attains its $\ell_2$ norm, up to an error of $O\(n^{5/2}\)$, between any two consecutive roots. That is, there are at least $s-1$ points $i$ between the minimal and the maximal roots of $K_s$ so that $\Omega\(\frac{{n \choose s}}{n^5}\) \le \frac{{n \choose i}}{2^n} K_s^2(i) \le {n \choose s}$. By (\ref{binomial-H}), and by the definition of $\tau$, for any such $i$ holds
\beqn
\label{Krawchouk-between-roots}
|K_s(i)| ~\in~ 2^{\(\tau\(\frac sn, \frac in\) \pm o(1)\) \cdot n}.
\eeqn

\item {\it Higher norms}. Let $h$ be the function defined in Section~\ref{subsubsec:ref:h}. The following estimate seems to be new\footnote{It might be known to experts, but we are unaware of its appearance in the literature.}. Let $p > 2$ be fixed. Let $n$ be sufficiently large. Let $\frac{n}{\ln n} < s < \frac n2 - \frac{n}{\ln n}$. Let $0 < y < \frac12$ be such that $h(p,y) = 1 - \frac{2s}{n}$. By Proposition~\ref{pro:2k-norm} the $\ell_p$ norm of $K_s$ is attained, up to a small factor, in the vicinity of $yn$. More precisely, for some (in fact for all) $i \in yn \pm o(n)$ holds $\|K_s\|^p_p \le 2^{o(n)} \cdot \frac{{n \choose i}}{2^n} \cdot |K_s(i)|^p$.

\noi By Lemma~\ref{lem:ref:h}, $yn < \frac n2 - \sqrt{s(n-s)}$. Hence, using (\ref{Krawchouk-tail}), we get the following estimate for $\|K_s\|^p_p$:
    \[
    \|K_s\|^p_p ~~\in~~ 2^{\(H(y) - 1 + p\tau\(\frac sn,y\) \pm o(1)\) \cdot n}.
    \]

    \noi Since the $\ell_2$ norm of $K_s$ is $\sqrt{{n \choose s}}$ this implies that for $\frac{n}{\ln n} < s < \frac n2 - \frac{n}{\ln n}$ holds
    \[
     \frac{||K_s||^p_p}{||K_s||^p_2} ~~\in~~ 2^{\(\psi\(p,\frac sn\) \pm o(1)\) \cdot n},
     \]
     where $\psi$ is the function defined in Section~\ref{subsubsec:ref:psi}.

    \noi We remark that Theorem~\ref{thm:norms} will imply that these estimates are valid for all $0 \le s \le n$.

\end{enumerate}

\subsection{Attaining norms between consecutive roots}

\noi The goal of this subsection is to show that Krawchouk polynomials attain their $\ell_2$ norm, within a polynomial factor, between any two of their consecutive roots, and also in the intervals below their first and above their last roots. We prove this property, in a somewhat higher degree of generality, for any family of polynomials orthogonal with respect to a discrete measure supported on $\{0,...,n\}$.

\noi Let $\mu$ be a positive measure on $\{0,...,n\}$, and let $\{P_0,...,P_n\}$ be the family of polynomials orthogonal with respect to the inner product $\<f,g\> = \sum_{i=0}^n \mu(i) f(i) g(i)$ induced by $\mu$, and normalized\footnote{Note that we need to choose a normalization to make the polynomials $\{P_s\}$ well-defined, however the specific choice of a normalization is immaterial for the discussion below.} so that $P_s(0) = 1$ for all $0 \le s \le n$.

\noi We will need two properties of orthogonal polynomials. First (see e.g., \cite{Szego}), for any $0 \le s \le  n$, the roots of $P_s$ are real and distinct, and lie in the interval $(0,n)$; and second that its $\ell_2$ norm $\|P_s\|_2 = \sqrt{\sum_{i=0}^n \mu(i) P_s^2(i)}$ is minimal among all polynomials of degree $s$ with the same leading coefficient. This is a simple and a well-known fact, but we provide an argument for completeness. Let $P$ be a polynomial of degree $s$ with the same leading coefficient as $P_s$. Then $Q = P - P_s$ is a polynomial of a smaller degree and hence is orthogonal to $P_s$. This means that $||P||^2_2 = ||P_s + Q||^2_2 = ||P_s||^2_2 + ||Q||^2_2 \ge ||P_s||^2_2$.

\noi We can now state our claim.

\pro
\label{pro:norm-any-two-roots}
Let $s > 0$. Let the roots of $P_s$ be $y_1 < y_2 < ... < y_s$. Assume that $y_1 \ge 1$ and that $y_s \le n-1$, and that the distance between any two consecutive roots is at least $2$. Assume also that the ratios $\frac{\mu(j)}{\mu(j+1)}$ and their inverses are uniformly bounded by some $R > 0$.

\begin{enumerate}

\item
The $\ell_2$ norm of $P_s$ is attained on the intervals $\left[0,y_1\right]$ and $\left[y_s,n\right]$ up to a factor of at most $O\(\sqrt{R n}\)$.

\item
For any $1 \le k \le s-1$ the $\ell_2$ norm of $P_s$ is attained between $y_k$ and $y_{k+1}$, up to a factor of at most $O\(\sqrt{R} n^2\)$.

\end{enumerate}

\epro

\prf

\noi We start with the first claim. We will prove it for the interval $\left[0,y_1\right]$, the proof for $\left[y_s,n\right]$ is similar. We will assume that the claim does not hold, and reach contradiction by constructing a polynomial $P$ as above with $||P||^2_2 < ||P_s||^2_2$. There are two cases to consider: $y_1$ is non-integer, and $y_1$ is an integer. We will deal only with the first case, the second case is similar (and easier).

\noi Assume then that the claim does not hold and that $y_1$ is not integer. Let $i_m = \lfloor y_1 \rfloor$. Let $a_s$ be the leading coefficient of $P_s$. That is $P_s(y) = a_s \cdot \prod_{j=1}^s \(y - y_j\)$. For a (small) parameter $\tau$ define the polynomial $P_{\tau}$ as follows: $P_{\tau}(y) = a_s \cdot \(y - y_1 - \tau\) \cdot \prod_{j=2}^s \(y - y_j\)$. Note that the roots of $P_{\tau}$, except for the first root, are those of $P_s$, and the first root is shifted downwards by $\tau$. In particular, $P_0 = P_s$. We will show that $\frac{d}{d \tau}_{| \tau = 0} ||P_{\tau}||^2_2 < 0$. This would mean that for some $\tau > 0$ we have $\|P_{\tau}\|^2_2 < \|P\|^2_2$, reaching a contradiction.

\noi A simple computation shows that $\frac{d}{d \tau}_{| \tau = 0} ||P_{\tau}||^2_2$ is proportional to $\sum_{i=0}^n \frac{\mu(i) P_s^2(i)}{{y_1 - i}}$. We write this expression as follows:
\beqn
\label{P-der-1}
\sum_{i=0}^{i_m-1} \frac{\mu(i) P_s^2(i)}{y_1 - i} ~+~ \frac{\mu\(i_m\) P_s^2\(i_m\)}{y_1 - i_m} ~+~ \sum_{i > y_1} \frac{\mu(i) P_s^2(i)}{y_1 - i}.
\eeqn
Since the denominators in the summands in the first sum are at least one, the first sum is bounded from above by $\sum_{i=0}^{i_m-1} \mu(i) P_s^2(i)$, which, by assumption, is at most $\frac{\|P_s\|_2^2}{CR n}$, for some constant $C$ which we may assume to be large. The last sum is negative. Since the denominators in its summands are at most $n$ in absolute value, its absolute value is bounded from below by $\frac{1}{n} \cdot \sum_{i > y_1} \mu(i) P_s^2(i)$, which, by assumption, is at least $\frac{1}{n} \cdot \(1-\frac{1}{CRn}\) \cdot \|P_s\|_2^2$.

\noi Finally, we need to bound the second summand. Note that both $P_s\(i_m - 1\)$ and $P_s\(i_m\)$ are positive, since $P_s$ is positive at $0$, and both points lie below the first root of $P_s$. Note also that
$\frac{P_s\(i_m - 1\)}{P_s\(i_m\)} = \prod_{i=1}^n \frac{y_i - \(i_m - 1\)}{y_i - i_m}
\ge \frac{y_1 - i_m + 1}{y_1 - i_m} > \frac{1}{y_1 - i_m}$. Hence
\[
\frac{\mu\(i_m\) P_s^2\(i_m\)}{y_1 - i_m} ~<~ R \cdot \(y_1 - i_m\) \mu\(i_m - 1\) P_s^2\(i_m-1\) ~<~ R \cdot \mu\(i_m - 1\) P_s^2\(i_m-1\) ~\le~ \frac{\|P_s\|^2_2}{Cn},
\]
Summing up, we see that for a sufficiently large constant $C$, the derivative  $\frac{d}{d \tau}_{| \tau = 0} ||P_{\tau}||^2_2$ is negative, proving the claim.

\noi We pass to the second claim of the proposition, proceeding via a similar line of argument. We will assume that both $y_k$, $y_{k+1}$ are non-integer. The other cases are similar (and simpler).

\noi For a parameter $\tau$ define $P_{\tau}(y) = a_s \cdot \prod_{j \not = k, k+1} \(y - y_j\) \cdot \(y - y_k + \tau\) \(y - y_{k+1} - \tau\)$. That is, we move the two roots in question {\it outwards} by $\tau$. A simple computation gives that $\frac{\partial}{\partial \tau}_{|\tau = 0} \|P_{\tau}\|^2_2 = \(y_k - y_{k+1}\) \cdot \sum_{i=0}^n \frac{\mu(i) P^2_s(i)}{\(i - y_k\)\(i - y_{k+1}\)}$. Hence the contribution of all integer points outside the region between the two roots is negative, and inside positive. We want to argue that if the norm inside is smaller than the total $\ell_2$ norm by a factor of more than $CRn^4$, for some sufficiently large constant $C$, then $\frac{d}{d \tau}_{| \tau = 0} ||P_{\tau}||^2_2$ is negative, reaching a contradiction.

\noi Dividing out by $|y_k - y_{k+1}|$, the outside contributes in absolute value at least $\frac{1}{n^2} \cdot \(1-\frac{1}{CR^2 n^2}\) \cdot \|P_s\|_2^2$. All the terms on the inside, for which the distance from both roots is at least $\frac{1}{2R}$, contribute together (note that the larger of these two distances is always at least $1$) at most $\frac{2}{C n^4} \cdot \|P_s\|^2_2$. It remains to deal with the inside terms which are close to one of the roots. Since the distance between the roots is at least $2$, there could be only one such term at the most. Say, $i$ is close to $y_k$ from the inside. But then the contribution of $i+1$ would be at least $\frac{R}{4n^2}$ that of $i$, by an argument similar to the argument above. Since $i+1$ contributes at most $\frac{2}{C n^4} \cdot \|P_s\|^2_2$, we have that $i$ contributes at most $\frac{8}{C R n^2} \cdot \|P_s\|^2_2$, and the total contribution of the inside is bounded by $\frac{10}{C n^2} \cdot \|P_s\|^2_2$. This means that for a sufficiently large constant $C$, the derivative  $\frac{d}{d \tau}_{| \tau = 0} ||P_{\tau}||^2_2$ is negative, proving the second claim, and completing the proof of the proposition.
\eprf

\cor
\label{cor:norm-any-two-roots}
Let $1 \le s \le \frac n2$, and let $x_s$ be the minimal root of the Krawchouk polynomial $K_s$. Then $K_s$ attains its $\ell_2$ norm within a factor of $O\(n\)$ on $\left[0,x_s\right]$ and within a factor of $O\(n^{5/2}\)$ between any two consecutive roots.
\ecor

\prf
Recall that for $1 \le s \le \frac n2$ the distance between any two consecutive roots of $K_s$ is at least $2$. We also use one additional facts about the Krawchouk polynomials: for $1 \le s \le \frac n2$, the first root of $K_s$ is at least $1$ (see \cite{Lev-chapter}). Hence we may apply the previous proposition with $\mu$ being the binomial measure on $\{0,...,n\}$. Note that the value of $R$ in this case is $n = \frac{\mu(1)}{\mu(0)}$. The claim of the corollary follows.
\eprf

\subsection{Hamming spheres}
\label{subsec:hamming}

\noi Let $f = 1_S$, where $S$ is the Hamming sphere of radius $s \le \frac n2$ around zero. Let $\phi$ be the function defined in Subsection~\ref{subsubsec:ref:phi}. Let $\sigma = \frac sn$.
Then (see e.g., \cite{ACKL}):
\[
\<T_{\e} f, f\> ~=~ \frac{1}{2^n} {n \choose s} \sum_{i=0}^{s} {s \choose i} {{n - s} \choose i} \e^{2i} (1-\e)^{n - 2i} ~\in~2^{\(\phi(\sigma, \e) \pm o(1)\) \cdot n}.
\]
The second step is by (\ref{binomial-H}) and Lemma~\ref{lem:ref:phi:disc-cont}.

\noi Since $\widehat{f} = \frac{1}{2^n} K_s$, and since $T_{\e} = \sum_{k=0}^n (1-2\e)^k \Pi_k$, we have, by Parseval's identity, that $\<T_{\e} f, f\> = \sum_{k=0}^n (1-2\e)^k \<f_k, f_k\> = \frac{1}{2^{2n}} \sum_{k=0}^n (1-2\e)^k {n \choose k} K_s^2(k)$. Using (\ref{Krawchouk-tail}) and (\ref{Krawchouk-between-roots}) it can bee seen that the last expression is in $2^{\(\max_{0 \le y \le \frac12} \left\{y \log_2(1 - 2\e) + H(y) + 2\tau(\sigma,y)\right\}-2 \pm o(1) \) \cdot n}$.

\noi Comparing the two expressions for $\<T_{\e} f, f\>$, the following identity should hold:
\beqn
\phi(\sigma, \e) ~~=~~ \max_{0 \le y \le \frac12} \Big\{y \log_2(1 - 2\e) + H(y) + 2\tau(\sigma,y) \Big\}- 2.
\eeqn
This is verified directly in Lemma~\ref{lem:ref:phi}. We remark that this identity, which shows that $\phi$ is, in an appropriate sense, a transform of $\tau$, might be considered as a step towards understanding of the somewhat 'arbitrary looking' functions $\tau$ and $I$.

\section{Some Proofs}
\label{sec:thms}

\noi In this section we prove theorems~\ref{thm:tails},~\ref{thm:edge-isop},~\ref{thm:NHC}~and~\ref{thm:max-proj} and some related results.

\noi Note that the arguments in this and in the following sections will rely, without further justification, on the properties of the bivariate functions detailed in Section~\ref{subsec:bivariate}.

\subsection{Proof of Theorem~\ref{thm:tails}}

\noi We start with the proof of (\ref{ineq:tails}), distinguishing two cases: $0 \le i < \frac n2 - \sqrt{s(n-s)}$, and $\frac n2 - \sqrt{s(n-s)} \le i \le \frac n2$.

\noi Consider first the case $\frac n2 - \sqrt{s(n-s)} \le i \le \frac n2$. By the definition of $\tau$, for $i$ in this range we have $\tau\(\frac sn, \frac in\) = \frac{1 + H\(\frac sn\) - H\(\frac in\)}{2}$. Therefore (\ref{ineq:tails}) reduces to
\[
\mathrm{Pr} \left\{|f| \ge \|f\|_2 \cdot 2^{\frac{1 - H\(\frac in\)}{2} \cdot n} \right\} ~\le~ 2^{\(H\(\frac in\) - 1\) \cdot n}.
\]
Set $t = \|f\|_2 \cdot 2^{\frac{1 - H\(\frac in\)}{2} \cdot n}$. Then, by Markov's inequality, $\mathrm{Pr} \left\{|f| \ge t \right\} \le \frac{\E f^2}{t^2} = 2^{\(H\(\frac in\) - 1\) \cdot n}$, completing the argument in this case.

\noi For $0 \le i < \frac n2 - \sqrt{s(n-s)}$, set $t = \|f\|_2 \cdot 2^{\(\tau\(\frac sn, \frac in\) - \frac12 H\(\frac sn\)\) \cdot n}$ and use Markov's inequality and Theorem~\ref{thm:norms} to obtain
\[
\mathrm{Pr} \left\{|f| \ge t \right\} ~\le~ \min_{p \ge 2} \left\{\frac{\E |f|^p}{t^p}\right\} ~=~ \min_{p \ge 2} \left\{\frac{\E |f|^p}{\(\E f^2\)^{\frac p2}} \cdot 2^{\(-p\tau\(\frac sn, \frac in\) + \frac p2 H\(\frac sn\)\) \cdot n}\right\} ~\le~
\]
\[
2^{n \cdot \min_{p \ge 2} \left\{\psi\(p,\frac sn\) - p\tau\(\frac sn, \frac in\) + \frac p2 H\(\frac sn\)\right\}}.
\]

\noi Consider the function $g(p) = \psi\(p,\frac sn\) - p\tau\(\frac sn, \frac in\) + \frac p2 H\(\frac sn\)$. We claim that its minimum on $[2, \infty)$ is attained at $p^{\ast}$, which is defined by $h\(p^{\ast}, \frac in\) = 1 - \frac{2s}{n}$. (Note that $p^{\ast}$ is well-defined in this range of $i$.) Since $\psi\(p,x\)$ is convex in $p$, the function $g$ is convex, and it will suffice to verify that $g'\(p^{\ast}\) = 0$. Proceeding as in the proof of Proposition~\ref{pro:psi-concave convex} below, we have that
\[
g'\(p^{\ast}\) ~=~ \frac{\partial \psi\(p, \frac sn\)}{\partial p}_{|p = p^{\ast}} - \tau\(\frac sn, \frac in\) + \frac12 H\(\frac sn\) ~=~ \tau\(\frac sn, y\) - \tau\(\frac sn, \frac in\),
\]
where $y$ is determined by $h\(p^{\ast},y\) = 1 - \frac{2s}{n}$. By the definition of $p^{\ast}$ we have $y = \frac in$, and therefore $g'\(p^{\ast}\) = 0$, as claimed.

\noi We now compute $g\(p^{\ast}\)$. Recall that $\psi(p,x) = H(y) - 1 + p\tau(x,y) - \frac p2 H(x)$, where $y$ is determined by $h(p,y) = 1 - 2x$. In our case $x = \frac sn$ and $y = \frac in$. Substituting, we get $g\(p^{\ast}\) = H\(\frac in\) - 1$, completing the proof of (\ref{ineq:tails}).

\noi We proceed to the second part of the theorem. Let $0 \le i \le \frac n2 - \sqrt{s(n-s)}$, and let $0 \le j \le i$. Recall that by (\ref{Krawchouk-tail}) we have for $0 \le j \le \frac n2 - \sqrt{s(n-s)}$ that $K_s(j) \ge {n \choose s} \cdot 2^{\(I\(\frac jn, \frac sn\) - I\(0,\frac sn\)\) \cdot n} = \frac{{n \choose s}}{2^{H\(\frac sn\) \cdot n}} \cdot 2^{\tau\(\frac sn, \frac jn\)} \ge \frac{{n \choose s}}{2^{H\(\frac sn\) \cdot n}} \cdot 2^{\tau\(\frac sn, \frac in\)}$.
Recall also that $\|K_s\|_2 = \sqrt{{n \choose s}}$.  Hence, using (\ref{binomial-H}) in the last inequality,
\[
\mathrm{Pr}\left\{|K_s| \ge \|K_s\|_2 \cdot \sqrt{\frac{{n \choose s}}{2^{H\(\frac sn\) \cdot n}}} \cdot 2^{\(\tau\(\frac sn, \frac in\) - \frac12 H\(\frac sn\)\) \cdot n}\right\} ~\ge~ \mathrm{Pr}\left\{K_s \ge \frac{{n \choose s}}{2^{H\(\frac sn\) \cdot n}} \cdot 2^{\tau\(\frac sn, \frac in\) \cdot n} \right\} ~\ge~
\]
\[
\frac{1}{2^n} \cdot \sum_{j=0}^i {n \choose j} ~\ge~ \Omega\(\frac{1}{\sqrt{i}}\) \cdot 2^{\(H\(\frac in\) - 1\) \cdot n}.
\]

\noi Next, by Corollary~\ref{cor:norm-any-two-roots}, between any two consecutive roots of $K_s$ there is a point $i$ for which $\frac{{n \choose i}}{2^n} \cdot K^2_s(i) \ge \Omega\(\frac{1}{n^5}\) \cdot \|K_s\|^2_2$. This means that
\[
\mathrm{Pr}\left\{|K_s| \ge \Omega\(\frac{1}{n^{5/2}}\) \cdot  \|K_s\|_2 \cdot \sqrt{\frac{2^n}{{n \choose i}}}\right\} ~\ge~ \frac{{n \choose i}}{2^n} ~\ge~ \Omega\(\frac{1}{\sqrt{i}}\) \cdot 2^{\(H\(\frac in\) - 1\) \cdot n}.
\]
Since in this range of $i$ we have $\sqrt{\frac{2^n}{{n \choose i}}} \ge 2^{\frac{1- H\(\frac in\)}{2} \cdot n} = 2^{\(\tau\(\frac sn, \frac in\) - \frac12 H\(\frac sn\)\) \cdot n}$, this proves the last claim of the theorem.

\eprf

\subsection{Two auxiliary claims}

\noi The following claim provides a key to all the remaining results in this section. Note that this is a special case of (\ref{ineq:max-proj}). Let $\pi$ be the function defined in Section~\ref{subsubsec:ref:pi}.
We write $x \wedge y$ for the minimum of $x$ and $y$.

\pro
\label{pro:max-proj}
Let $0 \le \sigma \le \frac 12$. Let $f$ be a function on $\H$ supported on a set of cardinality at most $2^{H(\sigma) n}$. Then, for any $0 \le k \le n$ holds
\[
\|f_k\|_2 \quad \le \quad 2^{\pi\(\sigma,\frac kn \wedge \frac{n-k}{n}\) \cdot n} \cdot \|f\|_2.
\]
\epro

\prf

\noi Given a function $f$, let $g$ be defined by $g(x) = (-1)^{|x|} f(x)$. Then $f$ and $g$ have supports of the same cardinality and (see Section~\ref{subsubsec:Fourier}) for any $0 \le k \le n$ holds $g_k = f_{n-k}$. Hence we may and will assume in the following argument that $0 \le k \le \frac n2$. Next, recall that $\pi\(\sigma, \frac kn\) = 0$ for $\frac kn \ge \frac12 - \sqrt{\sigma(1-\sigma)}$, reducing the claim of the proposition to the trivial inequality $\|f_k\|_2 \le \|f\|_2$. So, we may assume $0 \le \frac kn < \frac12 - \sqrt{\sigma(1-\sigma)}$.

\noi Let $f$ be a function on $\H$, supported on a subset $A \subseteq \H$. Let $0 \le k \le \frac n2$. Using the fact that $f_k$ is an orthogonal projection of $f$ in the first step and the Cauchy-Schwarz inequality in the last step, we have
\[
\<f_k,f_k\> ~=~ \<f,f_k\> ~=~ \< f \cdot 1_A, f_k\> ~=~ \<f,f_k \cdot 1_A\> ~\le~ \|f\|_2 \cdot \|f_k \cdot 1_A\|_2,
\]
implying that $\frac{\<f_k,f_k\>}{\|f\|_2} \le \|f_k \cdot 1_A\|_2$. On the other hand, for $p \ge 2$, we can apply H\"older's inequality to obtain $\|f_k \cdot 1_A\|^2_2 = \<f^2_k, 1_A\> \le \|1_A\|_{\frac{p}{p-2}} \cdot \|f^2_k\|_{\frac p2} = \(\frac{|A|}{2^n}\)^{\frac{p-2}{p}} \cdot \|f_k\|^2_p$. Combining the two estimates gives
\[
\frac{\|f_k\|_2}{\|f\|_2} ~=~ \frac{1}{\|f_k\|_2} \cdot \frac{\<f_k,f_k\>}{\|f\|_2} ~\le~ \(\frac{|A|}{2^n}\)^{\frac{p-2}{2p}} \cdot \frac{\|f_k\|_p}{\|f_k\|_2} ~\le~ \(\frac{|A|}{2^n}\)^{\frac{p-2}{2p}} \cdot 2^{\frac 1p \psi\(p,k/n\) \cdot n},
\]
where in the last step we have applied  Theorem~\ref{thm:norms}, using the fact that $f_k$ is a homogeneous polynomial of degree $k$.

\noi Since, by assumption, $|A| \le 2^{H\(\sigma\) n}$, we have that for any $p \ge 2$ holds
\[
\|f_k\|_2 ~~\le~~ 2^{\(\frac{p-2}{2p} \cdot (H(\sigma)-1) + \frac 1p \psi\(p,\frac kn\)\) \cdot n} \cdot \|f\|_2.
\]

\noi Since $\frac kn < \frac12 - \sqrt{\sigma(1-\sigma)}$, there is a unique $p > 2$ such that $h(p,\sigma) = 1 - 2\frac kn$. Fix this $p$. We then have, by the first definition of $\psi$, that $\psi\(p,\frac kn\) = H(\sigma) - 1 + p\tau\(\frac kn,\sigma\) - \frac p2 H\(\frac kn\)$. And hence $\frac{p-2}{2p} \cdot (H(\sigma)-1) + \frac 1p \psi\(p,\frac kn\) = \tau\(\frac kn, \sigma\) - \frac{1 + H\(\frac kn\) - H(\sigma)}{2} = \pi\(\frac kn, \sigma\) = \pi\(\sigma,\frac kn\)$, completing the proof of the proposition.
\eprf

\noi As a corollary we prove the following special case of (\ref{ineq:NHC}). Let $T_{\e}$ be the noise operator corresponding to a noise parameter $0 \le \e \le \frac12$ (see Section~\ref{subsubsec:Fourier}). Let $\phi$ be the function defined in Section~\ref{subsubsec:ref:phi}.

\cor
\label{cor:hypercontractive}

\noi Let $0 \le \sigma \le \frac 12$. Let $f$ be a function on $\H$ supported on a set of cardinality at most $2^{H(\sigma) n}$. Then, for any $0 \le \e \le \frac12$ holds
\[
\<T_{\e} f,f\> \quad \le \quad 2^{\(\phi\(\sigma,\e\) + 1 - H(\sigma)\) \cdot n} \cdot \|f\|^2_2.
\]
\ecor

\prf

\noi We have, using Proposition~\ref{pro:max-proj} in the first inequality, that
\[
\<T_{\e} f, f\> ~~=~~ \sum_{k=0}^n (1-2\e)^k \<f_k, f_k\> ~~\le~~ \|f\|^2_2 \cdot \sum_{k=0}^n (1-2\e)^k 2^{2\pi\(\sigma,\frac kn \wedge \frac{n-k}{n}\) \cdot n} ~~\le~~
\]
\[
n \|f\|^2_2 \cdot \max_{0 \le k \le \frac n2} \left\{(1-2\e)^k 2^{2\pi\(\sigma,\frac kn\) \cdot n}\right\} ~~\le~~ n \|f\|^2_2 \cdot 2^{\(\max_{0 \le y \le \frac12} \left\{y \log_2(1-2\e) + 2\pi(\sigma,y)\right\}\) \cdot n}.
\]
Since $\pi(x,y) = \tau(x,y) - \frac{1 + H(x) - H(y)}{2}$ we have
\[
\max_{0 \le y \le \frac12} \left\{y \log_2(1-2\e) + 2\pi(\sigma,y)\right\} ~=~ \max_{0 \le y \le \frac12} \left\{y \log_2(1-2\e) + H(y) + 2\tau(\sigma,y)\right\} - (1 + H(\sigma)) ~=~
\]
\[
\phi(\sigma,\e) + 1 - H(\sigma),
\]
where the last equality is by Lemma~\ref{lem:ref:phi}. So, $\<T_{\e} f,f\> \le n \cdot 2^{\(\phi\(\sigma,\e\) + 1 - H(\sigma)\) \cdot n} \cdot \|f\|^2_2$.

\noi Finally, we remove the extra $n$-factor by a tensorization argument (see Section~\ref{subsubsec:tensorization}). For an integer $m \ge 1$, let $F_m = f^{\otimes m}$. Observe that $F_m$ is supported on a subset of $\{0,1\}^{nm}$ of cardinality at most $2^{H(\sigma) nm}$.
In addition, $\<F_m, F_m\> = \<f,f\>^m$ and  $\<T_{\e} F_m, F_m\> = \<T_{\e} f,f\>^m$. Hence, by the above argument, we have
\[
\<T_{\e} f,f\> ~=~ \<T_{\e} F_m, F_m\>^{\frac 1m} ~\le~ (nm)^{\frac 1m} \cdot \(2^{\(\phi\(\sigma,\e\) + 1 - H(\sigma)\) \cdot nm}\)^{\frac 1m} \cdot \(\|F_m\|^2_2\)^{\frac 1m} ~=
\]
\[
(nm)^{\frac 1m} \cdot 2^{\(\phi\(\sigma,\e\) + 1 - H(\sigma)\) \cdot n} \cdot \|f\|^2_2.
\]
Taking $m$ to infinity, gives the claim of the corollary.
\eprf

\noi We can now prove Theorems~\ref{thm:edge-isop}~to~\ref{thm:max-proj}.

\subsection*{Proof of Theorem~\ref{thm:edge-isop} and related statements}

\rem
\label{rem:edge-isop}
Let us first briefly explain why $1 \le i \le 2\sigma(1-\sigma)n$ is the relevant range of parameters. (See also the discussion following the proof of Theorem~5 in \cite{ACKL}.) Let $s = \sigma n$ and assume $s$ to be integer, and $i = 2\sigma(1-\sigma)n = \frac{2s(n-s)}{n}$ to be an even integer. If $A$ is the Hamming sphere of radius $s$ around zero then $i$ is the expected distance between two points chosen uniformly at random from $A$, and it is easy to see that, up to at most an $O\(\frac{1}{\sqrt{n}}\)$-factor, we have $a_i(A) = |A|^2$. Hence in this case we cannot expect a non-trivial upper bound on $a_i(A)$.

\noi For a larger $i$, that is $2\sigma(1-\sigma)n \le i \le \frac n2$, write $i = \frac{2t(n-t)}{n}$, and, assuming $t > s$ to be integer, choose $A$ to be a random subset of cardinality ${n \choose s}$ of the sphere of radius $t$ around zero, to a similar effect. For $\frac n2 < i \le n$, let $j = n - i$, and let $B = A \cup \bar{A}$, where $\bar{A}$ is the shift of $A$ by an all-$1$ vector. Then $0 \le j < \frac n2$, $|A| \le |B| \le 2|A|$, $a_i(B) = a_j(B)$ and  $a_i(A) \le a_j(B) \le 2 \cdot \(a_i(A) + a_j(A)\)$.
\erem

\subsubsection*{Proof of (\ref{ineq:edge-isop})}

\noi We follow the argument in the proof of Theorem~5 from \cite{ACKL} replacing the hypercontractive inequality (\ref{ineq:HC}) used in \cite{ACKL} with Corollary~\ref{cor:hypercontractive}.

\noi Let $\e$ be between $0$ and $\frac 12$. Note that for any function $f$ on $\H$ holds $\<f_{\e},f\> = \frac{1}{2^n} \sum_{x,y} f(x) f(y) \e^{|x+y|} (1-\e)^{n - |x+y|}$. Substituting $f = 1_A$ gives $\<f_{\e},f\> = \frac{1}{2^n} \sum_{i=0}^n a_i(A) \e^i (1-\e)^{n-i}$. On the other hand, Corollary~\ref{cor:hypercontractive}, gives $\<f_{\e},f\> \le
2^{\(\phi\(\sigma,\e\) + 1 - H(\sigma)\) \cdot n} \cdot \|f\|^2_2 = |A| \cdot 2^{\(\phi\(\sigma,\e\) - H(\sigma)\) \cdot n}$.

\noi Hence, for any $0 \le i \le n$ and for any $0 < \e \le \frac12$ we have
\[
a_i(A)  ~\le~ \frac{\sum_{k=0}^n a_k(A) \e^k (1-\e)^{n-k}}{\e^i (1-\e)^{n-i}} ~\le~ |A| \cdot \frac{2^{\(\phi\(\sigma,\e\) + 1 - H(\sigma)\) \cdot n}}{\e^i (1-\e)^{n-i}}.
\]
Minimizing over $\e$ gives
\[
a_i(A) ~~\le~~ |A| \cdot 2^{\min_{0 < \e \le \frac12} \left\{\phi\(\sigma,\e\) + 1 - H(\sigma) - \frac in \log_2(\e) - \(1 - \frac in\) \log_2(1-\e) \right\} \cdot n} ~\le~
\]
\[
|A| \cdot 2^{\(\sigma H\(\frac{i}{2\sigma n}\) + (1-\sigma)H\(\frac{i}{2(n-\sigma n)}\)\) \cdot n},
\]
where the last step is via Lemma~\ref{lem:ref:phi:edge-isop}. \eprf

\subsubsection*{Proof of near tightness of (\ref{ineq:edge-isop}) for spheres or unions of spheres}

\noi For an even $i = 2j$, let $A$ be a sphere of radius $s$. Then $a_i(A) = |A| \cdot {s \choose j} {{n-s} \choose j}$, while (\ref{ineq:edge-isop}) gives $a_i(A) \le |A| \cdot 2^{sH\(\frac js\) + (n-s) H\(\frac{j}{n-s}\)}$. By (\ref{binomial-H}), the upper bound provided by (\ref{ineq:edge-isop}) is larger than $a_i(A)$ by at most a factor of $\Theta\(j \cdot \sqrt{\frac{(s-j)(n-s-j)}{s(n-s)}}\) \le O\(i\)$.

\noi For an odd $i = 2j-1$, let $A$ be the union of two adjacent spheres of dimension $n-1$ and radii $s-1$ and $s$. Note that $|A| = {n \choose s}$, and that $a_i(A) = 2|A|{s \choose j} {{n-s-1} \choose {j-1}}$. As above, this loses a factor of at most $O(i)$ to $|A| \cdot 2^{sH\(\frac{j}{s}\) + (n-s-1)H\(\frac{j-1}{n-s-1}\)}$. Next, a simple analysis, which we omit, shows that $sH\(\frac{i}{2s}\) + (n-s)H\(\frac{i}{2(n-s)}\) = sH\(\frac{2j-1}{2s}\) + (n-s)H\(\frac{2j-1}{2(n-s)}\)$ is larger than $sH\(\frac{j}{s}\) + (n-s-1)H\(\frac{j-1}{n-s-1}\)$ by at most $\frac12 \log_2\(\frac{n-s}{s}\) + O(1)$, and hence the upper bound on $a_i(A)$ provided by (\ref{ineq:edge-isop}) is tight, up to a factor of $O\(\sqrt{\frac{n-s}{s}} \cdot i\)$.
\eprf

\subsubsection*{Proof of Corollary~\ref{cor:Kleitman-West}}

\noi First, note that for any $t \ge 1$ holds $tH\(\frac 1t\) = \log_2(t) + (t-1) \log_2\(\frac{t}{t-1}\) \le \log_2(t) + \frac{1}{\ln 2}$, where in the last step we have used the fact that $\ln(1 + x) \le x$ for any $x > -1$.

\noi Let $1 \le s \le \frac n2$ be an integer. Using the observation above and (\ref{ineq:edge-isop}) with $\sigma = \frac sn$ and $i=2$, gives
\[
a_2(A) ~\le~ |A| \cdot 2^{sH\(\frac 1s\) + (n-s)H\(\frac{1}{n-s}\)} ~\le~ e^2 s(n-s) \cdot |A|.
\]
\eprf

\subsubsection*{Proof of (\ref{worst error exponent})}
\label{subsubsec:8}

\noi If $R = 1$, the claim becomes $P_{\mathrm{ue}}(1,\e) = 0$, which is correct, since in this case the undetected error probability is $1 - (1-\e)^n \approx 1$. So, we may assume $0 < R < 1$. Let $\sigma = H^{-1}(R)$. Let $n$ be a large integer. We will assume, somewhat unaccurately, that $s = \sigma n$ is integer, whenever required (to do this with full accuracy we would sandwich $P_{\mathrm{ue}}(R,\e)$ between $P_{\mathrm{ue}}\(H\(\frac{\lfloor \sigma n\rfloor}{n}\) ,\e\)$ and  $P_{\mathrm{ue}}\(H\(\frac{\lceil \sigma n\rceil}{n}\) ,\e\)$ and proceed similarly).

\noi Let $A(\sigma,n)$ be the union spheres of dimension $n-1$ and radii $s-1$ and $s$ around zero. We will show below that $A$ maximizes the undetected error probability for all codes in $\H$ of cardinality at most $2^{H\(\sigma\) n} = 2^{Rn}$, up to a factor of at most $O\(s^{\frac12} n^{\frac32}\)$. This will imply that
\[
P_{\mathrm{ue}}(R,\e) ~=~ \limsup_{n \rarrow \infty} \(\frac 1n \log_2\big(P_{\mathrm{ue}}\(A(\sigma,n),\e\)\big)\).
\]

\noi Let $A = A(\sigma,n)$, and let $f = 1_A$. Then $f = g + h$, where $g$ is the characteristic function of the sphere of radius $s$ around zero and $h$ is the characteristic function of the sphere of radius $s-1$ around zero. Note that $P_{\mathrm{ue}}\(A,\e\) =  \frac{2^n}{|A|} \cdot \<f_{\e},f\> - (1 - \e)^n$. We claim that for $\e > 0$ the first of these terms is exponentially in $n$ larger than the second one, and consequently $P_{\mathrm{ue}}\(A,\e\) \approx  \frac{2^n}{|A|} \cdot \<f_{\e},f\>$. Indeed, we have
\[
\frac 1n \log_2\(\frac{2^n}{|A|} \cdot \<f_{\e},f\>\)  ~\ge~ \frac 1n \log_2\(\frac{2^n}{|A|} \cdot \<g_{\e},g\>\) ~\in~ \phi(\sigma,\e) + 1 - H(\sigma) \pm o_n(1) ~=~
\]
\[
\max_{0 \le x \le \sigma} \Big\{\alpha_{\sigma,\e}(x)\Big\} \pm o_n(1).
\]
For the second step, see Section~\ref{subsec:hamming}, and for the third step Section~\ref{subsubsec:ref:phi}. As stated in Section~\ref{subsubsec:ref:phi}, the value of $x$ for which the maximum is attained is $x^{\ast} = x(\sigma, \e) = \frac{-\e^2 + \e \sqrt{\e^2 + 4(1-2\e) \sigma(1-\sigma)}}{2(1-2\e)}$, which is strictly positive for $\sigma, \e > 0$. And it is easy to see that $\alpha_{\sigma,\e}\(x^{\ast}\) > \alpha_{\sigma,\e}(0) = \log_2(1-\e)$.

\noi Next, note that $\<g_{\e},g\> \le \<f_{\e},f\> = \<g_{\e} + h_{\e}, g + h\> \le 2 \cdot \(\<g_{\e},g\> + \<h_{\e},h\>\)$. For the last step recall that the noise operator is a positive semidefinite linear operator, and hence we have the Cauchy-Schwarz inequality $\<g_{\e}, h\> \le \<g_{\e},g\>^{\frac12} \cdot \<h_{\e},h\>^{\frac12}$, and similarly for $\<g, h_{\e}\>$. Since both $\frac 1n \log_2\(\frac{2^n}{|A|} \cdot \<g_{\e},g\>\)$ and  $\frac 1n \log_2\(\frac{2^n}{|A|} \cdot \<h_{\e},h\>\)$ are in $\alpha_{\sigma,\e}\(x^{\ast}\) \pm o_n(1)$, this implies that also $\frac 1n \log_2\(P_{\mathrm{ue}}\(A,\e\)\) \in \alpha_{\sigma,\e}\(x^{\ast}\) \pm o_n(1)$. Hence
\[
P_{\mathrm{ue}}(R,\e) ~=~ \limsup_{n \rarrow \infty} \(\frac 1n \log_2\big(P_{\mathrm{ue}}\(A(\sigma,n),\e\)\big)\) ~=~ \alpha_{\sigma,\e}\(x^{\ast}\),
\]
proving (\ref{worst error exponent}).

\noi To conclude the proof it remains to show that $A$ maximizes the undetected error probability up to a polynomial factor. Let $t = \frac{2s(n-s)}{n}$. We first claim that the number of pairs of points in $A$ at distances between $1$ and $t$ from each other is not negligible. More precisely, $\sum_{i=1}^t  a_i(A) \ge \Omega\(n^{-\frac32}\) \cdot |A|^2$. This can be shown by observing that distance distributions inside a Hamming sphere and between two distinct Hamming spheres are closely related to hypergeometric distributions with appropriate parameters, and by applying straightforward first moment estimates for hypergeometric distributions. We omit the details.

\noi Let now $C \subseteq \H$ with $|C| \le 2^{H(\sigma)n}$. On one hand we have by Theorem~\ref{thm:edge-isop} that
\[
\frac{1}{|C|} \cdot \sum_{i=1}^t a_i(C) \e^i (1-\e)^{n-i} ~\le~ O\(s^{-\frac12} n^{\frac32}\) \frac{1}{|A|} \cdot \sum_{i=1}^t a_i(A) \e^i (1-\e)^{n-i} ~\le~ O\(s^{-\frac12} n^{\frac32}\) \cdot P_{\mathrm{ue}}(A,\e).
\]
On the other hand, we have
\[
\frac{1}{|C|} \cdot \sum_{i=t+1}^n a_i(C) \e^i (1-\e)^{n-i} ~\le~ \e^{t+1} (1-\e)^{n-t-1} \cdot \frac{1}{|C|} \sum_{i=t+1}^n a_i(C) ~\le~  \e^{t+1} (1-\e)^{n-t-1} \cdot |C| ~\le
\]
\[
\e^t (1-\e)^{n-t} \cdot O\(s^{\frac12}\) |A| ~\le~ O\(s^{\frac12}n^{\frac32}\) \cdot \frac{1}{|A|}  \sum_{i=1}^t a_i(A) \e^i (1-\e)^{n-i} ~\le~ O\(s^{\frac12}n^{\frac32}\) \cdot P_{\mathrm{ue}}(A,\e).
\]
We use (\ref{binomial-H}) in the third step, and the inequality $\sum_{i=1}^t a_i(A) \ge \Omega\(n^{-\frac32}\) \cdot |A|^2$ in the fourth step.

\noi Combining these two inequalities gives
\[
P_{\mathrm{ue}}(C,\e) ~=~ \frac{1}{|C|} \cdot \sum_{i=1}^n a_i(C) \e^i (1-\e)^{n-i} ~\le~  O\(s^{\frac12}n^{\frac32}\) \cdot P_{\mathrm{ue}}(A,\e).
\]

\eprf

\subsection*{Proof of Theorem~\ref{thm:NHC}}

\subsubsection*{Proof of (\ref{ineq:NHC})}

\noi Note that it would suffice to show a somewhat weaker statement:
\beqn
\label{thm-NHC-weak}
\|f_{\e}\|_2 ~\le~ 2^{o(n)} \cdot 2^{\eta\(\frac 1n \log_2\(\frac{\|f\|_{1 + (1-2\e)^2}}{\|f\|_1}\),~\e\) \cdot n} \cdot \|f\|_{1 + (1-2\e)^2},
\eeqn
since the $2^{o(n)}$ term by can be removed by a tensorization argument, like in the proof of Corollary~\ref{cor:hypercontractive}. We proceed to show (\ref{thm-NHC-weak}), with an error term which is polynomial in $n$.

\noi We may assume that $f \ge 0$, since replacing $f$ with $|f|$ increases the LHS of (\ref{thm-NHC-weak}) and does not change the RHS. We may also assume, by homogeneilty, that $\|f\|_1 = 1$. This means that $\|f\|_{\infty} \le 2^n$, and that the points at which $f < 2^{-n}$, say, contribute little to both sides ot the inequality, so we may ignore them for the sake of discussion (that is, we may and will assume that $f$ vanishes on these points). All the remaining points can be partitioned into $O(n)$ level sets $A_1,...A_r$ such that $f$ varies by a factor of $2$ at most in each level set. Let $\alpha_i = \frac 1n \log_2\(\frac{|A_i|}{2^n}\)$, and let $\nu_i = \frac 1n \log_2\(v_i\)$, where $v_i$ is, say, the median value of $f$ on $A_i$. Then, up to an additive error term of $O\(\frac{\log(n)}{n}\)$, we have, for any $p \ge 1$, that
\[
\frac 1n \log_2 \|f\|_p ~~\approx~~ \frac 1n \log_2\(\(\sum_{i=1}^r \frac{|A_i|}{2^n} \cdot v_i^p\)^{1/p}\) ~~\approx~~ \max_{1 \le i \le r} \left\{\frac{\alpha_i - 1}{p} + \nu_i\right\}.
\]
Here we use the approximate equality sign "$\approx$" to register that the equality holds up to a negligible error.

\noi Next, we estimate the LHS of (\ref{thm-NHC-weak}) in terms of $\{\alpha_i\}$ and $\{\nu_i\}$. Let $f_i$ be the restriction of $f$ to $A_i$. Then $f_{\e} = \sum_{i=1}^r \(f_i\)_{\e}$, and
we have that,up to an additive error term of $O\(\frac{\log(n)}{n}\)$,
\[
\frac 1n \log_2 \<f_{\e}, f_{\e}\> ~~\approx~~ \frac 1n \log_2\(\sum_{i=1}^r \<\(f_i\)_{\e}, \(f_i\)_{\e}\>\) ~~=~~ \frac 1n \log_2\(\sum_{i=1}^r \<\(f_i\)_{2\e(1-\e)}, f_i\>\)  ~~\approx~~
\]
\[
\max_{1 \le i \le r} \left\{\frac 1n \log_2 \<\(f_i\)_{2\e(1-\e)}, f_i\>\right\} ~~\approx~~ \max_{1 \le i \le r} \left\{\frac 1n \log_2\(v_i^2 \cdot \<\(1_{A_i}\)_{2\e(1-\e)}, 1_{A_i}\>\)\right\} ~~\le~~
\]
\[
\max_{1 \le i \le r} \Big\{\tilde{\phi}\(\alpha_i, 2\e(1-\e)\) + 2 \nu_i\Big\},
\]
The first step follows from the Cauchy-Schwarz inequality, the second step uses the semigroup property of noise operators, and the last step follows from Corollary~\ref{cor:hypercontractive}, and the definition of $\tilde{\phi}$.

\noi We will show (\ref{thm-NHC-weak}) to be a simple corollary of the following lemma. Given $0 \le \alpha_1,...\alpha_r \le 1$ and $0 \le \nu_1,...\nu_r$, we write $N(p)$ for $\max_{1 \le i \le r} \left\{\frac{\alpha_i - 1}{p} + \nu_i\right\}$.

\lem
\label{lem:NHC}
Let $0 \le \delta \le 1/2$, and let $p \ge 2-2\delta$. Then, for any $0 \le \alpha_1,...\alpha_r \le 1$ and $0 \le \nu_1,...\nu_r$ holds
\[
\max_{1 \le i \le r} \Big\{\tilde{\phi}\(\alpha_i, \delta\) +  2\nu_i\Big\} ~~\le~~ \tilde{\phi}\(\alpha^{\ast}, \delta\) + 2\nu^{\ast},
\]
where $\alpha^{\ast} = 1 - \frac{p}{p-1} \cdot \(N(p) - N(1)\)$, and $\nu^{\ast} = \frac{p}{p-1} \cdot N(p) - \frac{1}{p-1} \cdot N(1)$.
\elem

\prf
We need to show that for all $1 \le i \le r$ holds $\tilde{\phi}\(\alpha_i, \delta\) + 2\nu_i \le \phi\(\alpha^{\ast},\delta\) + 2 \nu^{\ast}$. Fix $i$. Recall that $\frac{1}{1-\delta} \le \frac{\partial \tilde{\phi}(\alpha,\delta)}{\partial \alpha} \le 2$. There are two cases to consider.

\begin{enumerate}

\item $\nu_i \le \nu^{\ast}$.

\noi If also $\alpha^{\ast} \ge \alpha_i$, the claim follows from the monotonicity of $\tilde{\phi}$ in $\alpha$. If, on the other hand, $\alpha^{\ast} < \alpha_i$, then $\tilde{\phi}\(\alpha_i, \delta\) - \tilde{\phi}\(\alpha^{\ast},\delta\) \le 2\alpha_i - 2\alpha^{\ast}$, and hence it only remains to verify that $\alpha_i + \nu_i \le \alpha^{\ast} + \nu^{\ast}$. To see that note that
\[
\alpha_i + \nu_i ~~\le~~ 1 + N(1) ~~=~~ \alpha^{\ast} + \nu^{\ast}.
\]

\item $\nu_i > \nu^{\ast}$.

\noi Note that $\frac{\alpha_i}{p} + \nu_i \le \frac 1p + N(p) = \frac{\alpha^{\ast}}{p} + \nu^{\ast}$. In particular, we have that $\alpha^{\ast} > \alpha_i$. And hence
\[
\tilde{\phi}\(\alpha^{\ast},\delta\) - \tilde{\phi}\(\alpha_i, \delta\) ~~\ge~~ \frac{1}{1-\delta} \cdot \(\alpha^{\ast} - \alpha_i\) ~~\ge~~ 2 \cdot \frac{\alpha^{\ast} - \alpha_i}{p} ~~\ge~~ 2\nu_i - 2\nu^{\ast}.
\]

\end{enumerate}

\eprf

\noi We can now conclude the proof of (\ref{thm-NHC-weak}). Given a function $f$, define, as above, the partition of $\H$ into level sets $A_1,...,A_r$ of $f$, and define $\{\alpha_i\}$ and $\{\nu_i\}$ correspondingly. Apply the lemma with $p = 1 + (1-2\e)^2$ and $\delta = 2\e(1-\e)$. By the discussion above, and bearing in mind that $\Big |\frac{\partial \tilde{\phi}(\alpha,\e)}{\partial \alpha} \Big |$ is bounded by a constant (in fact by $2$), we have, up to an additive error term of $O\(\frac{\log(n)}{n}\)$, that
\[
\frac 1n \log_2 \|f_{\e}\|^2_2 ~=~ \frac 1n \log_2 \<f_{\e}, f_{\e}\> ~\lesssim~ \max_{1 \le i \le r} \Big\{\tilde{\phi}\(\alpha_i, 2\e(1-\e)\) + 2 \nu_i\Big\} ~\le~ \tilde{\phi}\(\alpha^{\ast}, 2\e(1-\e)\) + 2\nu^{\ast} ~=
\]
\[
\tilde{\phi}\(1 - \frac{p}{p-1} \cdot \(N(p) - N(1)\), ~2\e(1-\e)\) + 2 \cdot \frac{N(p) - N(1)}{p-1} + 2N(p) ~\approx
\]
\[
\tilde{\phi}\(1 - \frac{p}{p-1} \cdot \frac 1n \log_2\(\frac{\|f\|_p}{\|f\|_1}\), ~2\e(1-\e)\) + 2 \cdot \frac{\frac 1n \log_2\(\frac{\|f\|_p}{\|f\|_1}\)}{p-1} + \frac 2n \log_2\(\|f\|_p\) ~=
\]
\[
2\eta_p\(\frac 1n \log_2\(\frac{\|f\|_p}{\|f\|_1}\),~\e\) + \frac 2n \log_2\(\|f\|_p\) ~=~ 2\eta\(\frac 1n \log_2\(\frac{\|f\|_{1 + (1-2\e)^2}}{\|f\|_1}\),~\e\) + \frac 2n \log_2\(\|f\|_{1 + (1-2e)^2}\).
\]
\eprf

\subsubsection*{Proof of near tightness of (\ref{ineq:NHC}) for spheres}

\pro
\label{pro:tight-edge-isop}
Let $f$ be proportional to $1_S$, where $S$ is a Hamming sphere of radius $s$. Then (\ref{ineq:NHC}) is tight for $f$ up to a factor of $O\(s^{3/4}\)$.
\epro

\prf

\noi We write $p$ for $1 + (1-2\e)^2$. We may and will assume, for simplicity, that $\|f\|_1 = 1$ (which means $f = \frac{2^n}{{n \choose s}} \cdot 1_S$). This reduces (\ref{ineq:NHC}), after taking binary logarithms of both sides, and expanding the definition of $\eta$, into
\[
\frac 1n \log_2\(\|f_{\e}\|_2\) ~~\le~~ \frac12 \tilde{\phi}\(1 - \frac{p}{p-1} \cdot \frac 1n \log_2\(\|f\|_p\), ~2\e(1-\e)\) + \frac{p}{p-1} \cdot \frac 1n \log_2\(\|f\|_p\).
\]
We want to show this is nearly tight for $f$. Let $\sigma = \frac sn$. We proceed by comparing both sides of this inequality to $M = \frac12 \phi(\sigma, 2\e(1-\e)) + 1 - H(\sigma)$.

\noi First, consider the RHS. Let $g(x) = \frac12 \tilde{\phi}\(1 - \frac{p}{p-1} \cdot x, ~2\e(1-\e)\) + \frac{p}{p-1} \cdot x$. Then the RHS is $g\(\frac 1n \log_2\(\|f\|_p\)\)$. Observe that $M = g\(\frac{p-1}{p} (1 - H(\sigma))\)$. Note also that $\frac 1n \log_2\(\|f\|_p\) = \frac{p-1}{p} \(1 - \frac 1n \log_2\({n \choose s}\)\)$. Let $x = \frac{p-1}{p} (1 - H(\sigma))$ and let $y = \frac{p-1}{p} \(1 - \frac 1n \log_2\({n \choose s}\)\)$. By (\ref{binomial-H}), we have that $y \ge x$ and that $y - x \le \frac{p-1}{p} \cdot \(\frac{1}{2n} \log_2(s) + O\(\frac 1n\)\)$, where the asymptotic notation hides absolute constants. Since $\tilde{\phi}$ increases and $2 \ge \tilde{\phi}' \ge \frac{1}{1-2\e(1-\e)} \ge 1$, we have that $g$ increases and $g' \le \frac{p}{2p-2}$. This implies that $g(y) \ge g(x)$ and $|g(y) - g(x)| \le \frac12 \cdot \frac 1n \log_2(s) + O\(\frac{1}{n}\)$. In other words, $|\mathrm{RHS} - M| \le \frac{1}{2n} \log_2(s) + O\(\frac{1}{n}\)$.

\noi Next, consider the LHS. We have that it equals (writing $\delta$ for $2\e(1-\e)$):
\[
\frac{1}{2n} \log_2 \(\<f_{\e}, f_{\e}\>\) ~=~ \frac{1}{2n} \log_2 \(\<f_{\delta}, f\>\) ~=~ \frac{1}{2n} \log_2 \(\frac{2^n}{{n \choose s}} \sum_{i=0}^{s} {s \choose i}{{n-s} \choose i} \delta^{2i} (1-\delta)^{n-2i}\) ~\ge
\]
\[
\frac{1}{2n} \log_2 \(\frac{2^n}{{n \choose s}}\) + \frac{1}{2n} \max_{0 \le i \le s} \log_2\({s \choose i}{{n-s} \choose i} \delta^{2i} (1-\delta)^{n-2i}\).
\]
In the first step we have used the semigroup property of noise operators. For the second step, see Section~\ref{subsec:hamming}.

\noi By (\ref{binomial-H}) we have  $\frac{1}{2n} \log_2 \(\frac{2^n}{{n \choose s}}\) \ge \frac{1 - H(\sigma)}{2} + \frac{1}{4n} \log_2(s) - O\(\frac{1}{n}\)$. Similarly, by (\ref{binomial-H}):
\[
\frac{1}{2n} \max_{0 \le i \le s} \log_2\({s \choose i}{{n-s} \choose i} \delta^{2i} (1-\delta)^{n-2i}\) ~\ge~
\]
\[
\frac{1}{2} \max_{0 \le i \le s} \left\{\sigma H\(\frac{i/n}{\sigma}\) + (1-\sigma) H\(\frac{i/n}{1-\sigma}\) + 2\frac in\log_2(\delta) + \(1-2\frac in\) \log_2(1-\delta)\right\} - \frac{1}{2n} \log_2(s) - O\(\frac{1}{n}\)
\]
\[
\ge~ \frac{1}{2} \max_{0 \le x \le \sigma} \left\{\sigma H\(\frac{x}{\sigma}\) + (1-\sigma) H\(\frac{x}{1-\sigma}\) + 2x\log_2(\delta) + \(1-2x\) \log_2(1-\delta)\right\} - \frac{1}{2n} \log_2(s) - O\(\frac{1}{n}\).
\]
The second inequality is by Lemma~\ref{lem:ref:phi:disc-cont}. Summing up, we have that the LHS is bounded from below by
\[
\frac{1}{2} \max_{0 \le x \le \sigma} \left\{\sigma H\(\frac{x}{\sigma}\) + (1-\sigma) H\(\frac{x}{1-\sigma}\) + 2x\log_2(\delta) + \(1-2x\) \log_2(1-\delta)\right\} +
\]
\[
\frac{1 - H(\sigma)}{2}  - \frac{1}{4n} \log_2(s) - O\(\frac{1}{n}\) ~=~
\]
\[
=~ \frac12 \phi\(\sigma, ~2\e(1-\e)\) + 1 - H(\sigma) - \frac{1}{4n} \log_2(s) - O\(\frac{1}{n}\) ~=~ M - \frac{1}{4n} \log_2(s) - O\(\frac{1}{n}\).
\]

\noi The first step follows from the definition of $\phi$. Wrapping everything up, we have that $\mathrm{LHS} - \mathrm{RHS} \le \frac{3}{4n} \log_2(s) + O\(\frac{1}{n}\)$ and hence the hypercontractive inequality is tight for $f$ up to a multiplicative factor of $O\(s^{3/4}\)$.

\eprf

\subsubsection{Proof of Corollary~\ref{cor:sphere-stable}}

\noi We may and will assume, by homogeneity, that $\|f\|_1 = 1$. Let $F(p)$ be the (normalized) binary logarithm of the RHS of the inequality in the claim of the corollary. That is $F(p)  = \eta_p\(\frac 1n \log_2\(\|f\|_p\),~\e\) + \frac 1n \log_2\(\|f\|_p\)$. We will show that $F(p)$ increases in $p$, and hence the claim of the corollary for $p \ge 1 + (1-2\e)^2$ follows from the claim for $p = 1 + (1-2\e)^2$, proved in (\ref{ineq:NHC}). Expanding the definition of $\eta_p$, we have, as in the proof of Proposition~\ref{pro:tight-edge-isop}, that
\[
F(p) ~=~ \frac12 \tilde{\phi}\(1 - \frac{p}{p-1} \cdot \frac 1n \log_2\(\|f\|_p\), ~2\e(1-\e)\) + \frac{p}{p-1} \cdot \frac 1n \log_2\(\|f\|_p\).
\]
Since the derivative of $\tilde{\phi}$ with respect to its first argument is bounded from above by $2$, it suffices to show that $\frac{p}{p-1} \cdot \log_2\(\|f\|_p\)$ is increasing in $p$ to infer that $F$ is increasing. Let $G(t) = \log_2\(\|f\|_{1/t}\)$, for $0 < t \le 1$. The function $G$ is decreasing and convex (this is a consequence of H\"older's inequality, see \cite{HLP}, theorems 196-197). Moreover, $G(1) = 0$, since $\|f\|_1 = 1$. It is easy to see that this implies that $\frac{p}{p-1} \cdot \log_2\(\|f\|_p\) = \frac{G\(\frac 1p\) - 1}{1-\frac 1p}$ is increasing in $p$.
\eprf

\subsection*{Proof of Theorem~\ref{thm:max-proj}}

\subsubsection*{Proof of (\ref{ineq:max-proj})}

\noi First, as is observed at the beginning of the proof of Proposition~\ref{pro:max-proj}, it suffices to consider the case $0 \le k \le \frac n2$.

\noi Let $0 \le \e < \frac12$ and consider the action of the noise operator $T_{\e}$ on $f$. Since $T_{\e} = \sum_{k=0}^n (1-2\e)^k \Pi_k$ (see Section~\ref{subsubsec:Fourier}), we have that $f_{\e} = \sum_{k=0}^n (1-2\e)^k f_k$, and therefore $\<f_{\e}, f_{\e}\> = \sum_{k=0}^n (1-2\e)^{2k} \<f_k, f_k\>$. This implies that
$\<f_k, f_k\>$ is upperbounded by $\frac{\<f_{\e}, f_{\e}\>}{(1-2\e)^{2k}}$. Taking logarithms of both sides of this inequality, and using Corollary~\ref{cor:sphere-stable} in the second step, we get that
\[
\frac 1n \log_2\(\|f_k\|_2\) ~\le~ \frac 1n \log_2\(\|f_{\e}\|_2\) - \frac kn \log_2(1-2\e) ~\le~
\eta_p\(r(p),\e\) + \frac 1n \log_2\(\|f\|_p\) - \frac kn \log_2(1-2\e).
\]

\noi Hence (\ref{ineq:max-proj}) would follow if we verify the identity
\[
\min_{0 \le \e \le \frac12} \left\{\eta_p\(r(p),\e\) - \frac kn \log_2(1-2\e)\right\} ~=~ \pi\(\frac kn \wedge \frac{n-k}{n},H^{-1}\(1 - \frac{p}{p-1} \cdot r(p)\)\) - \frac{p-2}{2p - 2} \cdot r(p).
\]

\noi Writing $\sigma$ for $H^{-1}\(1 - \frac{p}{p-1} \cdot r(p)\)$, $\kappa$ for $\frac kn$, $\delta$ for $2\e(1-\e)$, and expanding all the definitions, this reduces to verifying that for all $0 \le \kappa,\sigma, \delta \le \frac12$ holds, writing $x = x(\sigma,\delta)$ for  $\frac{-\delta^2 + \delta \sqrt{\delta^2 + 4(1-2\delta) \sigma(1-\sigma)}}{2(1-2\delta)}$, that $\pi(\kappa,\sigma) = \pi(\sigma, \kappa)$ is given by
\[
\frac12 \min_{0 \le \delta \le \frac12} \left\{\sigma H\(\frac{x}{\sigma}\) + \(1-\sigma\) H\(\frac{x}{1-\sigma}\) + 2x\log_2(\delta) + (1-2x) \log_2(1-\delta) - \kappa \log_2(1-2\delta)\right\}.
\]
This is shown in Lemma~\ref{lem:ref:pi-min}, and (\ref{ineq:max-proj}) follows.

\subsubsection*{Proof of near tightness of (\ref{ineq:max-proj}) for spheres}

\pro
\label{pro:tight-max-proj}
Let $f$ be proportional to $1_S$, where $S$ is a Hamming sphere of radius $s$. Then (\ref{ineq:max-proj}) is tight for $f$ in the following sense: If $k \le \frac n2 - \sqrt{s(n-s)}$, then (\ref{ineq:max-proj}) is tight for $f$up to a factor of $O\(k^{1/4}\)$. Moreover, (\ref{ineq:max-proj}) is tight for $f$ up to a factor of $O(t)$, provided $k$ is a point at which the $\ell_2$ norm of $K_s$ is attained, up to a factor of $t$.
\epro

\prf

\noi Let $S$ be the Hamming sphere of radius $s$ around $0$. Let $f = 1_S$. Recall (see Section~\ref{subsubsec:Fourier}) that $\widehat{f} = \frac{1}{2^n} \cdot K_s$, and therefore for any $0 \le k \le n$ holds $f_k = \frac{1}{2^n} K_s(k) \cdot K_k$. Hence $\|f_k\|_2 = \frac{\sqrt{{n \choose k}}}{2^n} |K_s(k)|$. In particular, $\|f_k\|_2 = \|f_{n-k}\|_2$, and in the following argument it suffices to consider the case $0 \le k \le \frac n2$.

\noi Let $R$ be the RHS of (\ref{ineq:max-proj}). Then $\frac 1n \log_2\(\|f\|_p\) = 1 - \frac 1n \log_2\({n \choose s}\)$ and $r(p) = \frac{p-1}{p} \cdot \(1 - \frac 1n \log_2\({n \choose s}\)\)$. Substituting, we get $\frac 1n \log_2(R) = \pi\(\frac kn, H^{-1}\(\frac 1n \log_2\({n \choose s}\)\)\) - \frac12 \(1 - \frac 1n \log_2\({n \choose s}\)\)$.

\noi We consider two cases.

\begin{enumerate}

\item $0 \le k \le \frac n2 - \sqrt{s(n-s)}$.

\noi In this case, see (\ref{Krawchouk-tail}), we have $K_s(k) \ge \frac{{n \choose s}}{2^{H\(\frac sn\) \cdot n}} \cdot 2^{\tau\(\frac sn, \frac kn\) \cdot n}$. Hence, recalling that $\tau(x,y) = \pi(y,x) - \frac{H(y) - H(x) - 1}{2}$, we have, after some rearrangement, that
\[
\frac 1n \log_2(\|f_k\|_2) ~\ge~ \pi\(\frac kn, \frac sn\) - \frac12 \cdot \(1 - \frac 1n \log_2\({n \choose s}\)\) + \frac 12 \cdot \(\frac 1n \log_2\({n \choose s}\) - H\(\frac sn\)\) +
\]
\[
\frac 12 \cdot \(\frac 1n \log_2\({n \choose k}\) - H\(\frac kn\)\).
\]
By the monotonicity of $\pi$ we have that $\pi\(\frac kn, H^{-1}\(\frac 1n \log_2\({n \choose s}\)\)\) \le \pi\(\frac kn, \frac sn\)$ and hence
\[
\frac 1n \log_2(R) - \frac 1n \log_2(\|f_k\|_2) ~\le~ \frac 12 \cdot \(H\(\frac sn\) - \frac 1n \log_2\({n \choose s}\)\) +
\]
\[
\frac 12 \cdot \(H\(\frac kn\) - \frac 1n \log_2\({n \choose k}\)\) ~\le~ \frac14 \log_2(ks) + O(1),
\]
where in the last inequality we have used (\ref{binomial-H}). This proves the first part of the proposition.

\item $K_s$ attains its $\ell_2$ norm on $k$ up to a factor of $t$.

\noi This means that $\frac{1}{2^n} {n \choose k} K_s^2(k) \ge \Omega\(\frac{1}{t^2}\) \cdot \|K_s\|^2_2 = \Omega\(\frac{1}{t^2}\) \cdot {n \choose s}$, which implies $|K_s(k)| \ge \Omega\(\frac{1}{t}\) \cdot \sqrt{\frac{{n \choose s} 2^n}{{n \choose k}}}$. Hence $\|f_k\|_2 \ge \Omega\(\frac{1}{t}\) \cdot \sqrt{\frac{{n \choose s}}{2^n}}$ and, recalling that $\pi$ is non-positive,
\[
\frac 1n \log_2(R) - \frac 1n \log_2(\|f_k\|_2) ~\le~ \pi\(\frac kn, \frac 1n \log_2\({n \choose s}\)\) - \log_2(t) + O(1) ~\le~ \log_2(t) + O(1).
\]
This proves the second part of the proposition.
\end{enumerate}

\eprf

\noi To complete the proof of the tightness of (\ref{ineq:max-proj}) for spheres recall that, by Corollary~\ref{cor:norm-any-two-roots}, between any two consecutive roots of $K_s$ there is a point on which $K_s$ attains its $\ell_2$ norm up to a factor of $O\(n^{5/2}\)$.

\eprf

\section{Proof of Theorem~\ref{thm:norms} and related results}
\label{sec:norms}

\noi We first deduce Corollary~\ref{cor:also balls} from Theorem~\ref{thm:norms}.

\prf (of Corollary~\ref{cor:also balls}):

\noi Let $g$ be a polynomial of degree $s$. Write $g = \sum_{r=0}^s a_r f_r$, where $f_r$ is a homogeneous polynomial of degree $r$, $0 \le r \le s$. By the triangle inequality for the $\ell_p$ norm we have $\|f\|_p \le \sum_{r = 0}^s |a_r| \|f_r\|_p$. On the other hand, the Parseval identity gives $\|f\|_2 = \sqrt{\sum_{r=0}^s a_r^2 \|f\|_r^2} \ge \max_{0 \le r \le s} |a_r| \|f_r\|_2$. Note also that (\ref{ineq:norms}) is equivalent to $\frac{\|f\|_p}{\|f\|_2} \le 2^{\frac{\psi\(p, \frac sn\)}{p} \cdot n}$.

\noi Taking all of the above into account, and recalling that $\psi(p,x)$ increases in $x$, we have that
\[
\frac{\|g\|_p}{\|g\|_2} ~\le~ \frac{\sum_{r = 0}^s |a_r| \|f_r\|_p}{\max_{0 \le r \le s} |a_r| \|f_r\|_2} ~\le~ \sum_{r = 0}^s \frac{\|f_r\|_p}{\|f_r\|_2} ~\le~ \sum_{r = 0}^s 2^{\frac{\psi\(p, \frac rn\)}{p} \cdot n} ~\le~ n \cdot 2^{\frac{\psi\(p, \frac sn\)}{p} \cdot n}.
\]

\noi We proceed with a tensorization argument. For an integer $m \ge 1$, let $G_m = g^{\otimes m}$. Note that $G_m$ is a polynomial of degree at most $sm$ on $\{0,1\}^{nm}$. By the above,
\[
\frac{\|g\|_p}{\|g\|_2} ~=~ \(\frac{\|G_m\|_p}{\|G_m\|_2}\)^{\frac 1m} ~\le~ \(nm \cdot 2^{\frac{\psi\(p, \frac sn\)}{p} \cdot nm}\)^{\frac 1m} ~=~ (nm)^{\frac 1m} \cdot 2^{\frac{\psi\(p, \frac sn\)}{p} \cdot n}.
\]

\noi Taking $m$ to infinity gives $\frac{\|g\|_p}{\|g\|_2} \le 2^{\frac{\psi\(p, \frac sn\)}{p} \cdot n}$.

\eprf

\noi We proceed with the proof of Theorem~\ref{thm:norms}. First, we introduce some notation.  Let $R(n,s,p)$ be the maximum of the ratio $\frac{||f||^p_p}{||f||^p_2}$ over all homogeneous polynomials of degree $s$ on $\H$. Let $r(n,s,p) = \frac{||K_s||^p_p}{||K_s||^p_2}$. Then (\ref{ineq:norms}) becomes
\[
R(n,s,p) ~\le~ 2^{\psi\(p,\frac sn\) \cdot n}.
\]

\noi The key step required to show this is the following claim.
\thm
\label{thm:key-step}
Let $p \ge 2$ be fixed. Then, for all $0 \le s \le n/2$ holds
\[
R(n,s,p) ~\le~ 2^{O\(\frac{n}{\log(n)}\)} \cdot r(n,s,p).
\]
Here the constant in the asymptotic notation may depend on $p$.
\ethm

\noi The inequality (\ref{ineq:norms}) will follow from Theorem~\ref{thm:key-step} and the following limit estimate.
\lem
\label{lem:psi}
For any integers $n \ge 1$ and $0 \le s \le \frac n2$, and for any $p \ge 2$ holds
\[
\lim_{m \rarrow \infty} \Big(r(nm,sm,p)\Big)^{\frac 1m} ~=~ 2^{\psi\(p,\frac sn\) \cdot n}.
\]
\elem

\noi In fact, assume Theorem~\ref{thm:key-step} and Lemma~\ref{lem:psi} to hold.
Let $f$ be a homogeneous polynomial of degree $s$ on $\H$, such that $R(n,s,p) = \frac{\E f^p}{\E^{p/2} f^2}$. We proceed with a tensorization argument. For an integer $m \ge 1$, let $F_m = f^{\otimes m}$. Then $F_m$ is a homogeneous polynomial of degree $sm$ on $\{0,1\}^{nm}$.
\noi Hence,
\[
R(n,s,p) ~=~ \frac{\E f^p}{\E^{p/2} f^2} ~=~ \(\frac{\E F_m^p}{\E^{p/2} F_m^2}\)^{\frac 1m} ~\le~ R(nm,sm,p)^{\frac 1m}  ~\le~ \(2^{O\(\frac{nm}{\log(nm)}\)} \cdot r(nm,sm,p)\)^{\frac 1m},
\]
where the second inequality follows from Theorem~\ref{thm:key-step}.  Taking $m$ to infinity, and using Lemma~\ref{lem:psi}, gives $R(n,s,p) \le 2^{\psi\(p,\frac sn\) \cdot n}$, establishing the first claim of the theorem.

\noi The second claim of Theorem~\ref{thm:norms} will be dealt with in the following proposition.
\pro
\label{pro:gap}
There is an absolute constant $C$ such that for any integers $n \ge 1$ and $0 \le s \le \frac n2$, and for any $p \ge 2$ holds
\[
2^{n \cdot\psi\(p,\frac sn\)} ~\le~ n \cdot C^p \cdot s^{\frac p4} \cdot r(n,s,p).
\]
\epro

\noi In the remainder of this section we prove Theorem~\ref{thm:key-step}. Lemma~\ref{lem:psi}~and~Proposition~\ref{pro:gap} will be proved in Section~\ref{sec:lemmas}.

\noi {\it Notation:} For the duration of this section let $s_0 = s_0(n) = \frac{n}{\ln n}$. Let $\e = \e(n) = \frac{n^{\frac{11}{2}} \ln^{\frac12} n}{s_0^6}$. Note that $\e(n)$ behaves like $\frac{1}{\sqrt{n}}$, up to polylogarithmic factors. The proof of Theorem~\ref{thm:key-step} will rely on the following four claims.

\lem
\label{lem:small-s}
Theorem~\ref{thm:key-step} holds for all $s \le s_0$.
\elem
\prf
We use (\ref{ineq:moments-HC}). Since $r(n,s,p) \ge 1$, we have
\[
R(n,s,p) ~\le~ (p-1)^{\frac{ps}{2}} ~\le~ (p-1)^{\frac{ps_0}{2}} ~\le~ 2^{\frac{p \log_2(p-1)}{2} \cdot \frac{n}{\ln n}} ~\le~  2^{\frac{p \log_2(p-1)}{2} \cdot \frac{n}{\ln n}}  \cdot r(n,s,p).
\]

\eprf

\lem
\label{lem:large-s}
Theorem~\ref{thm:key-step} holds for all $\frac n2 - s_0 \le s \le \frac n2$.
\elem
\prf
Let $\delta_0$ be the characteristic function of $0$. Clearly,
\[
R(n,s,p) ~\le~ \frac{\E \delta^p_0}{\E^{p/2} \delta^2_0} ~=~ 2^{\(\frac p2 - 1\)n}.
\]

\noi On the other hand, recall that $K_s(0) = ||K_s||^2_2 = {n \choose s}$. Note also that, by (\ref{binomial-H}), we have ${n \choose s} \ge {n \choose \frac n2 - s_0} \ge \Omega\(\frac{1}{\sqrt n}\) \cdot 2^{H\(\frac 12  - \frac{s_0}{n}\) \cdot n} \ge 2^{n - O\(\frac{s_0^2}{n}\)}$. Hence
\[
r(n,s,p) ~\ge~ \frac{\frac{1}{2^n} \cdot K^p_s(0)}{||K_s||^p_2} ~=~ \frac{{n \choose s}^{p/2}}{2^n} ~\ge~
\frac{2^{\(\frac p2 -1 \) n}}{2^{O\(\frac{p s_0^2}{n}\)}} ~\ge~ 2^{-O\(\frac{n}{\ln^2 n}\)} \cdot R(n,s,p) .
\]
\eprf

\noi The proofs of the next two claims are harder. We will first state the claims and show how to deduce Theorem~\ref{thm:key-step} from the preceding two lemmas and these two claims and then prove the claims.

\pro
\label{pro:induct}
There exists an explicitly defined (see (\ref{F-def})) function $F = F_p$ of two nonnegative variables $x$ and $y$ such that
\begin{enumerate}
\item
The function $F$ is increasing in both $x$ and $y$ is $1$-homogeneous.
\item
For any $1 \le s \le (n+1)/2$ the following inductive relation holds
\[
R(n+1,s,p) \quad \le \quad F\(R(n,s,p), R(n,s-1,p)\).
\]
\end{enumerate}
\epro

\pro
\label{pro:F}
There exists a sufficiently large constant $n_0$ such that for all $n \ge n_0$ and for all $s_0(n+1) \le s \le (n+1)/2 - s_0(n+1)$ holds
\[
r(n+1,s,p) \quad \in \quad \(1 \pm O\(\e\)\)^{2p} \cdot  F\Big(r(n,s,p), r(n,s-1,p)\Big).
\]
\epro

\noi We now prove Theorem~\ref{thm:key-step}, assuming the four claims above to hold, and proceeding similarly to \cite{KS}. We will show by induction on $n$ that for all $n$ and for all $1 \le s \le \frac n2$ holds $R(n,s,p) \le 2^{c \frac{n}{\log(n)}} \cdot r(n,s,p)$, for some constant $c$ which may depend on $p$.

\noi For any fixed $n_0$, we may assume, by choosing $c$ to be sufficiently large, that the claim holds for $n \le n_0$, which takes care of the base step. We pass to the induction step. Assume the claim holds for $n$ and we will show that it holds for $n+1$ as well. Let $1 \le s \le (n+1)/2$ be given. We may and will assume that $n \ge n_0$, for a sufficiently large $n_0$. By Lemmas~\ref{lem:small-s}~and~\ref{lem:large-s} the claim holds for $s \le s_0 = s_0(n+1)$ and for $s \ge (n+1)/2 - s_0$. So we may assume $s_0 < s < (n+1)/2 - s_0$.

\noi Let $R_0 = R(n,s,p)$ and $R_1 = R(n,s-1,p)$. By Proposition~\ref{pro:induct} $R(n+1,s,p) \le F\(R_0,R_1\)$. Let $\rho = \max\left\{\frac{R_0}{r(n,s,p)}, \frac{R_1}{r(n,s-1,p)}\right\}$. By the induction hypothesis $\rho \le 2^{c \frac{n}{\log(n)}}$. By the monotonicity and $1$-homogeneity of $F$ given in Proposition~\ref{pro:induct}, and by Proposition~\ref{pro:F}, we have that
\[
F\(R_0,R_1\) ~\le~ F\Big(\rho \cdot r(n,s,p), \rho \cdot r(n,s-1,p)\Big) ~=~ \rho \cdot F\Big(r(n,s,p), r(n,s-1,p)\Big) ~\le~
\]
\[
2^{c \frac{n}{\log(n)}} \cdot F\Big(r(n,s,p), r(n,s-1,p)\Big) ~\le~
\]
\[
2^{c \frac{n}{\log(n)}} \cdot \(1 + O\(\e\)\)^{2p} \cdot r(n+1,s,p) ~\le~ 2^{c \frac{n+1}{\log(n+1)}} \cdot r(n+1,s,p),
\]

\noi completing the proof of Theorem~\ref{thm:key-step}. Note that the last inequality holds, for a sufficiently large $n$, since $\e = \tilde{O}\(\frac{1}{\sqrt{n}}\)$.

\eprf

\subsection{Proof of Proposition~\ref{pro:induct}}

\noi Let $p$ be given. We start with defining the function $F = F_p$ at a point $(x,y)$ where $x, y \ge 0$. If $y = 0$, let $F(x,y) = x$. If $y \not = 0$, let $\rho = \(\frac xy\)^{2/p}$. Let $P(z) = \frac{\(\sqrt{z} + 1\)^p + \Big | \sqrt{z} - 1 \Big |^p}{2}$.
We define $F(x,y)$ by
\beqn
\label{F-def}
F(x,y) ~~=~~ y \cdot \mbox{sup}_{\beta \in [0, \infty)} \frac{P(\rho \beta)}{(\beta+1)^{p/2}}.
\eeqn

\noi By definition, $F$ is clearly $1$-homogeneous. Since, as is easy to see, $P$ increases in $z$, for $z \ge 0$, we also have that $F$ is increasing in $x$. To see that $F$ increases in $y$, substitute $\alpha = \rho \beta$ and note that $F(x,y) = x \cdot \mbox{sup}_{\alpha \in [0, \infty)} \frac{P(\alpha)}{(\alpha+\rho)^{p/2}}$.

\noi We now proceed similarly to the proof of Proposition~4.5 in \cite{KS}.

\noi Let $f$ be a homogeneous polynomial of degree $s$ over $\{0,1\}^{n+1}$, such that $\frac{\E f^p}{\E^{p/2} f^2} = R(n+1,s,p)$. For $i = 0,1$ let $f_i$ be the restriction of $f$ to the $n$-dimensional subcube $\{x: x_{n+1} = i\}$. We view both of these subcubes as isomorphic to $\H$. Note that there is a homogeneous polynomial $g_0$ of degree $s$ over $\H$ and a homogeneous polynomial $g_1$ of degree $s-1$ over $\H$, such that $f_0 = g_0 + g_1$ and $f_1 = g_0 - g_1$. We sum up the above by writing $f \leftrightarrow \(g_0 + g_1, g_0 - g_1\)$.

\noi Let us first deal with the case in which one of the functions $g_i$ vanishes. If $g_1 = 0$ then $f \leftrightarrow \(g_0, g_0\)$, and hence $R(n+1,s,p) = \frac{\E f^p}{\E^{p/2} f^2} = \frac{\E g_0^p}{\E^{p/2} g_0^2} \le R(n,s,p)$. To see that the claim of the proposition holds it remains to verify that $x \le F(x,y)$. This however is true, since $F(x,0) = y$ and $F$ increases in $y$. Similarly, if $g_0 = 0$, we have $R(n+1,s,p) \le R(n,s-1,p)$. In this case we need to verify $y \le F(x,y)$. This is true, since $F(x,y) \ge y P(0) = y$.

\noi From now on we assume that both $g_i$ do not vanish. Let $R_0 =\frac{\E g_0^p}{\E^{p/2} g_0^2}$, and let $R_1 =\frac{\E g_1^p}{\E^{p/2} g_1^2}$. Note that $R_0 \le R(n,s,p)$ and $R_1 \le R(n,s-1,p)$. We use Hanner's inequality \cite{Lieb-Loss}: For $p \ge 2$ holds
\[
\|g_0+g_1\|_p^p + \|g_0-g_1\|_p^p ~\le~ \(\|g_0\|_p + \|g_1\|_p\)^p + \Big | \|g_0\|_p - \|g_1\|_p \Big |^p.
\]

\noi This implies that
\[
R(n+1,s,p) ~=~ \frac{\E f^p}{\E^{p/2} f^2} = \frac{\frac12 \cdot \(\|g_0+g_1\|_p^p + \|g_0-g_1\|_p^p\)}{\(\E g^2_0 + \E g^2_1\)^{p/2}} ~\le
\]
\[
\frac{\frac12 \cdot \(\(\|g_0\|_p +\|g_1\|_p\)^p + \Big | \|g_0\|_p - \|g_1\|_p\Big |^p\)}{\(\E g^2_0 + \E g^2_1\)^{p/2}} ~=~ R_1 \cdot \frac{\frac12 \cdot \(\(\frac{\|g_0\|_p}{\|g_1\|_p} + 1\)^p + \Big | \frac{\|g_0\|_p}{\|g_1\|_p} - 1\Big |^p\)}{\(\frac{\E g^2_0}{\E g^2_1} + 1\)^{p/2}}.
\]

\noi Let $\rho = \(\frac{R_0}{R_1}\)^{2/p}$, and $\beta = \frac{\E g^2_0}{\E g^2_1}$. Then $\frac{\|g_0\|_p}{\|g_1\|_p} = \sqrt{\rho \beta}$, and the last expression can be written as

\[
R_1 \cdot  \frac{\frac12 \cdot \(\(\sqrt{\rho \beta} + 1\)^p + \Big | \sqrt{\rho \beta} - 1\Big |^p\)}{\(\beta + 1\)^{p/2}} ~=~ R_1 \cdot \frac{P(\rho \beta)}{\(\beta + 1\)^{p/2}} ~\le~ F\(R_0,R_1\),
\]
where the last inequality follows from the definition of $F$.

\eprf

\subsection{Proof of Proposition~\ref{pro:F}}

\noi There are two functions on $\left[2, \infty\right] \times \left[0, \frac12\right]$ which will play an important role in the following argument. The first of these functions is the function $h(p,x)$ defined in Section~\ref{subsubsec:ref:h}. We define the second function to be $g(p,x) = x^{\frac{1}{p}} (1-x)^{\frac{p-1}{p}} - x^{\frac{p-1}{p}} (1-x)^{\frac{1}{p}}$. Note that $g$ is nonnegative. For fixed $p$ we will frequently omit the first variable and view $h$ and $g$ as functions of $x$ only.

\noi Given $n$, $s \le n/2$, and $p$, we define $i_0 = i_0(n,s,p)$ to be the unique real number in the interval $\left[0,n/2\right]$ satisfying
\beqn
\label{i_0}
1 - \frac{2s}{n} ~=~ h\(p,\frac{i_0}{n}\).
\eeqn


\noi We now define several quantities depending on $n, s, p$ and $i_0$. Assume $s > 0$. Let $t = t(n,s,p) = \frac{\(n-2i_0\) + \sqrt{\(n-2i_0\)^2 - 4s(n-s)}}{2(n-s)}$. Let
\beqn
\label{rho}
\rho(n,s,p) \queq \frac{n-2i_0}{s} \cdot t - 1,
\eeqn
and let
\beqn
\label{F}
\Phi(n,s,p) \queq \frac{n}{2(n-i_0)} \cdot \(\frac{s}{n}\)^{p/2} \cdot \(1 + \frac{n-s}{s} \cdot t\)^p.
\eeqn

\noi The claim of Proposition~\ref{pro:F} will be based on the following two claims.

\pro
\label{pro:exact}
Let $F$ be the function defined in (\ref{F-def}). Then, assuming $0 < s < n/2$, we have
\[
\Phi(n,s,p) \queq F\(\rho^{p/2}(n,s,p),1\)
\]
\epro

\pro
\label{pro:approximate}
There exists a sufficiently large constant $n_0$ such that for all $n \ge n_0$ and for all $s_0 \le s \le \frac n2 - s_0$ holds
\begin{enumerate}

\item
\[
\(\frac{r(n,s,p)}{r(n,s-1,p)}\)^{2/p} ~\in~ \(1 \pm O\(\e\)\)^2 \cdot \rho(n,s,p).
\]

\item
\[
\frac{r(n+1,s,p)}{r(n,s-1,p)} ~\in~ \(1 \pm O\(\e\)\)^p \cdot \Phi(n,s,p).
\]

\end{enumerate}
\epro

\noi We first derive Proposition~\ref{pro:F} from these two claims and then prove the claims. By $1$-homogeneity of $F$, the claim of the proposition is equivalent to
\[
\frac{r(n+1,s,p)}{r(n,s-1,p)} \quad \in \quad \(1 \pm O\(\e\)\)^{2p} \cdot  F\(\frac{r(n,s,p)}{r(n,s-1,p)}, 1\)
\]

\noi Assume Propositions~\ref{pro:exact}~and~\ref{pro:approximate} to hold. By the monotonicity of $F$ in both coordinates and by its $1$-homogeneity, we have that
\[
\frac{r(n+1,s,p)}{r(n,s-1,p)} ~\in~ \(1 \pm O\(\e\)\)^p \cdot \Phi(n,s,p) ~=~ \(1 \pm O\(\e\)\)^p \cdot F\(\rho^{p/2}(n,s,p),1\) ~\subseteq~
\]
\[
\(1 \pm O\(\e\)\)^{2p} \cdot F\(\frac{r(n,s,p)}{r(n,s-1,p)}, 1\).
\]

\eprf

\subsection{Proof of Proposition~\ref{pro:exact}}

\noi First, we observe that $\rho(n,s,p)$ lies between $1$ and $p-1$. This will be the contents of the following lemma.

\lem
\label{lem:rho-bounds}
For all $0 < s < n/2$ holds
\[
1 ~<~ \rho(n,s,p) ~<~ p - 1.
\]
\elem
\prf
\noi Let $x = \frac{i_0}{n}$. Then by (\ref{i_0}) we have $0 < x < 1/2$, and $h(x) = 1 - \frac sn$. Hence $\frac sn = \frac{1-h(x)}{2}$. In particular, $s$ is a function of $x$, and hence so is $t$. In fact,
$t = \frac{1-2x+g(x)}{1+h(x)}$. To see this, observe that a simple calculation gives
\[
t ~=~ \frac{\(n-2i_0\) + \sqrt{\(n-2i_0\)^2 - 4s(n-s)}}{2(n-s)} ~=~ \frac{(1-2x) + \sqrt{(1-2x)^2 - \(1-h^2(x)\)}}{1+h(x)}.
\]
Note that $h^2(x) - g^2(x) = 4x(1-x)$, and hence the last expression is indeed $\frac{1-2x+g(x)}{1+h(x)}$.

\noi Next, we write $\rho = \rho(n,s,p)$ as a function of $x$ as well:
$\rho = \frac{1- 2x + g(x)}{1- 2x - g(x)}$. This can be verified by a simple calculation, using again the identity $h^2(x) = g^2(x) + 4x(1-x)$:
\[
\rho ~=~ \frac{n-2i}{s} \cdot t - 1 ~=~ \frac{2(1-2x)(1-2x+g(x))}{1-h^2(x)} - 1 ~=~
\]
\[
\frac{\Big((1-2x) + g(x)\Big)^2}{(1-2x)^2 - g^2(x)} ~=~ \frac{1- 2x + g(x)}{1- 2x - g(x)}.
\]
\noi Since $g > 0$ for $0 < x  < 1/2$, this implies that $\rho > 1$.

\noi Next, we argue that $\rho < p - 1$. This is equivalent to $g < \frac{p-2}{p} \cdot (1-2x)$. Since both sides of this putative inequality vanish at $1/2$, it suffices to show that $g' > -\frac{2p-4}{p}$. Computing $g'$ and rearranging, we have that
\[
g' ~=~ \frac{1}{p} \cdot \(\(\frac{1-x}{x}\)^{\frac{p-1}{p}} + \(\frac{x}{1-x}\)^{\frac{p-1}{p}}\) - \frac{p-1}{p} \cdot \(\(\frac{1-x}{x}\)^{\frac{1}{p}} + \(\frac{x}{1-x}\)^{\frac{1}{p}}\).
\]

\noi Let $\gamma = \frac12 \cdot \(\(\frac{1-x}{x}\)^{\frac{1}{p}} + \(\frac{x}{1-x}\)^{\frac{1}{p}}\)$. Then $\gamma > 1$ and by convexity of the function $t \rarrow t^{p-1}$ we have $g' \ge \frac{2}{p} \cdot \gamma^{p-1} - \frac{2p-2}{p} \cdot \gamma$. Using this, it remains to verify the inequality $\gamma^{p-1} - (p-1) \gamma > -(p-2)$, for $\gamma > 1$, and this is true, since it holds with equality for $\gamma = 1$, and the derivative of the LHS is positive for $\gamma > 1$.

\eprf

\noi Next, we consider $F\(\rho^{p/2}(n,s,p),1\)$. Let $x = \rho^{p/2}(n,s,p)$ and $y = 1$. By the preceding lemma $\rho = \rho(n,s,p) = \(\frac{x}{y}\)^{2/p}$ lies between $1$ and $p-1$. Recall that $F(x,y) = y \cdot \mbox{sup}_{\beta \in [0, \infty)} \frac{P(\rho \beta)}{(\beta+1)^{p/2}}$. Let $f(\beta) = \frac{P(\rho \beta)}{(\beta+1)^{p/2}}$. The following lemma describes the behavior of $f$ when $1 < \rho < p-1$.

\lem
\label{lem:h for x > y}
Assume $1 < \rho < p - 1$. Then $f$ increases from $0$ to some point $\frac{1}{\rho} < \beta^{\ast} < \infty$ and decreases from $\beta^{\ast}$ on. In particular,
$\textup{sup}_{\beta \in [0, \infty)} \frac{P(\rho \beta)}{(\beta+1)^{p/2}} ~=~ \frac{P\(\rho \beta^{\ast}\)}{(\beta^{\ast}+1)^{p/2}}$.
\elem

\noi We will prove the lemma below. Here we proceed assuming that it holds. By the lemma, we may restrict our attention to the behavior of $f$ on $\left(\frac{1}{\rho}, \infty\right]$. In this interval $P(z) = \frac{(\sqrt{z} + 1)^p + (\sqrt{z} - 1)^p}{2}$. It will be convenient to make the one-to-one substitution $u = \frac{1}{\rho(\beta+1)}$. Then $0 < u < \frac{1}{\rho + 1}$ and $f(\beta) = \rho^{p/2} \cdot Q(u)$, where $Q(u) = \frac12 \cdot \(\(\sqrt{1 - \rho u} + \sqrt{u}\)^p + \(\sqrt{1 - \rho u} - \sqrt{u}\)^p\)$. In particular, $F(x,y) = y \cdot \rho^{p/2} \cdot Q\(u^{\ast}\) = \rho^{p/2} \cdot Q\(u^{\ast}\)$, where $u^{\ast}$ is the unique point in the interval $\(0,\frac12\)$ in which the derivative $Q'$ vanishes. A simple calculation gives that $u^{\ast}$ is implicitly given by the following identity (writing $u$ for $u^{\ast}$ for simplicity):
\beqn
\label{der-zer}
\(\sqrt{1 - \rho u} + \sqrt{u}\)^{2k-1} \cdot \(\sqrt{1 - \rho u} - \rho\sqrt{u}\)  \queq \(\sqrt{1 - \rho u} - \sqrt{u}\)^{2k-1} \cdot \(\sqrt{1 - \rho u} + \rho\sqrt{u}\).
\eeqn

\noi Given this, the claim of the proposition is immediately implied by the following two lemmas.

\lem
\label{lem:u-ast}
For $\rho = \rho(n,s,p)$ holds
\[
u^{\ast} \queq \frac{s}{\rho n}.
\]
\elem

\noi And

\lem
\label{lem:F-eq}
\[
\Phi(n,s,p) \queq \rho^{p/2} \cdot Q\(u^{\ast}\)
\]
\elem

\noi It remains to prove Lemmas~\ref{lem:h for x > y}~-~\ref{lem:F-eq}.

\noi {\bf Proof of Lemma~\ref{lem:h for x > y}}.

\noi We partition $[0,\infty]$ into two subintervals $\left[0,\frac{1}{\rho}\right]$ and $\left[\frac{1}{\rho}, \infty\)$. The claim of the lemma is implied by the following two claims.

\begin{enumerate}

\item

\noi On $\frac{1}{\rho} \le \beta < \infty$ the function $f$ increases up to some point $\frac{1}{\rho} < \beta^{\ast} < \infty$ and decreases from $\beta^{\ast}$ on.

\item

\noi On $0 \le \beta \le \frac{1}{\rho}$ the unction $f$ increases.

\end{enumerate}

\noi {\bf The case $\frac{1}{\rho} \le \beta < \infty$}: On this interval $f(\beta) = \frac{\(\sqrt{\rho \beta} + 1\)^p + \(\sqrt{\rho \beta} - 1\)^p}{(2\beta+1)^{p/2}}$ and, after some rearrangement and simplification, $f'$ is proportional to
\[
\frac{\beta+1}{\beta} \cdot \sqrt{\rho \beta} \cdot \(\(\sqrt{\rho \beta} + 1\)^{p-1} + \(\sqrt{\rho \beta} - 1\)^{p-1}\) - \(\(\sqrt{\rho \beta} + 1\)^p + \(\sqrt{\rho \beta} - 1\)^p\).
\]

\noi Set $z = \sqrt{\rho \beta}$. Then $z \ge 1$ and the above is proportional to
$\frac{z^2 + \rho}{z} - \frac{(z+1)^p + (z-1)^p}{(z+1)^{p-1} + (z-1)^{p-1}}$.

\noi Let
\[
t(z) ~=~  \frac{(z+1)^p + (z-1)^p}{(z+1)^{p-1} + (z-1)^{p-1}} \cdot z - z^2 ~=~ \frac{z \cdot \((z+1)^{p-1} - (z-1)^{p-1}\)}{(z+1)^{p-1} + (z-1)^{p-1}}.
\]

\noi Note that the sign of $f'$ is the same of that of $\rho - t(z)$. Hence, recalling that $1 < \rho < p-1$, the claim will follow if we show that the function $t(z)$ strictly increases from $1$ to $p-1$ on $[1, \infty)$.

\noi First, it is easy to see that $t(1) = 1$ and that $t(z) \rarrow_{z \rarrow \infty} p-1$. Next, we claim that $t' > 0$, which is the same as
\[
\(\((z+1)^{p-1} - (z-1)^{p-1}\) + (p-1)z \cdot \((z+1)^{p-2} - (z-1)^{p-2}\)\) \cdot \((z+1)^{p-1} + (z-1)^{p-1}\) >
\]
\[
(p-1)z \cdot \(\((z+1)^{p-1} - (z-1)^{p-1}\) \cdot \((z+1)^{p-2} + (z-1)^{p-2}\)\).
\]

\noi Rearranging and simplifying, this is the same as
\[
(z+1)^{2p-2} - (z-1)^{2p-2} ~>~ 4(p-1)z \cdot \(z^2-1\)^{p-2}.
\]

\noi Consider the function $r(z) = z^{2p-2}$. Since $2p-2 \ge 2$ we have that $r''' \ge 0$ and hence (by developing $r$ into Taylor series around $z$, up to the second term) that $r(z+1) - r(z-1) \ge 2 r'(z) = 4(p-1)z^{2p-3}$. Hence it suffices to show
\[
4(p-1)z^{2p-3} ~>~ 4(p-1)z \cdot \(z^2-1\)^{p-2},
\]
which is evidently true for $z \ge 1$.

\noi {\bf The case $0 \le \beta \le \frac{1}{\rho}$}: On this interval $f(\beta) = \frac{\(\sqrt{\rho \beta} + 1\)^p + \(1 - \sqrt{\rho \beta}\)^p}{(2\beta+1)^{p/2}}$. We have that $f'(0) = \infty$. Next, we consider $f'(\beta)$ for $\beta > 0$. We have, similarly to the above, that $f'$ is proportional to
\[
\frac{\beta+1}{\beta} \cdot \sqrt{\rho \beta} \cdot \(\(\sqrt{\rho \beta} + 1\)^{p-1} - \(\sqrt{1 - \rho \beta}\)^{p-1}\) - \(\(\sqrt{\rho \beta} + 1\)^p + \(\sqrt{1 - \rho \beta}\)^p\).
\]

\noi Set $z = \sqrt{\rho \beta}$. Then $0 < z \le 1$ and the above is proportional to
$\frac{z^2 + \rho}{z} - \frac{(z+1)^p + (1-z)^p}{(z+1)^{p-1} - (1-z)^{p-1}}$.

\noi Let
\[
t(z) ~=~  \frac{(z+1)^p + (1-z)^p}{(z+1)^{p-1} - (1-z)^{p-1}} \cdot z - z^2 ~=~ \frac{z \cdot \((z+1)^{p-1} + (1-z)^{p-1}\)}{(z+1)^{p-1} - (1-z)^{p-1}}.
\]

\noi The sign of $f'$ is the same of that of $\rho - t(z)$. Hence the claim will follow if we show that the function $t(z)$ strictly increases from $\frac{1}{p-1}$ to $1$ on $\left[0, \frac{1}{\rho}\right]$.

\noi First, it is easy to see that $\lim_{z \rarrow 0} t(0) = \frac{1}{p-1}$ and that $t(1) = 1$.

\noi Next, we claim that $t' > 0$, which, similarly to the discussion above is the same as
\[
(z+1)^{2p-2} - (1-z)^{2p-2} ~>~4(p-1)z \cdot \(z^2-1\)^{p-2}.
\]

\noi Consider the function $r(u) = u^{2p-2}$. Since $2p-2 \ge 2$ we have that $r''' \ge 0$ and hence that $r(1+z) - r(1-z) \ge 2z r'(1) = 4(p-1)z$, which is evidently larger then the RHS above.

\eprf

\noi {\bf Proof of Lemma~\ref{lem:u-ast}}.

\noi Let $u = \frac{s}{\rho n}$. Clearly $0 < u \le \frac{1}{2\rho} < \frac{1}{\rho + 1}$. So, it remains to verify that $u$ satisfies (\ref{der-zer}). As in the proof of Lemma~\ref{lem:rho-bounds}, we will write  everything as a function of $x = \frac{i_0}{n}$, reducing to an identity involving the functions $g = g(x)$ and $h = h(x)$, which we then proceed to verify.

\noi First, we observe that $u = \frac{1}{\rho} \cdot \frac sn = \frac{1-2x-g}{1-2x+g} \cdot \frac{1-h}{2}$. We also have $1 - \rho u = 1 - \frac sn = \frac{1 + h}{2}$, and that
$\rho^2 u = \rho \cdot \frac sn = \frac{1-2x+g}{1-2x-g} \cdot \frac{1-h}{2}$.

\noi Using this, substituting in (\ref{der-zer}) and simplifying, we need to show
\[
\(\sqrt{(1+h)(1-2x+g)} + \sqrt{(1-h)(1-2x-g)}\)^{p-1} \cdot \(\sqrt{(1+h)(1-2x-g)} - \sqrt{(1-h)(1-2x+g)}\) =
\]
\[
\(\sqrt{(1+h)(1-2x+g)} - \sqrt{(1-h)(1-2x-g)}\)^{p-1} \cdot \(\sqrt{(1+h)(1-2x-g)} + \sqrt{(1-h)(1-2x+g)}\).
\]

\noi Next, we multiply both sides by
\[
\(\sqrt{(1+h)(1-2x+g)} + \sqrt{(1-h)(1-2x-g)}\)^{p-1} \cdot \(\sqrt{(1+h)(1-2x-g)} + \sqrt{(1-h)(1-2x+g)}\).
\]

\noi We observe that
\[
(1+h)(1-2x+g) \cdot (1-h)(1-2x-g) = (1+h)(1-2x-g) \cdot (1-h)(1-2x+g) = \(1-h^2\) \cdot \((1-2x)^2 - g^2\) = \(1-h^2\)^2
\]
and hence, after some simplification,
\[
\(\sqrt{(1+h)(1-2x+g)} + \sqrt{(1-h)(1-2x-g)}\)^2 ~=~ 2\cdot \Big((1-2x) + gh + \(1-h^2\)\Big),
\]
and
\[
\(\sqrt{(1+h)(1-2x-g)} + \sqrt{(1-h)(1-2x+g)}\)^2 ~=~ 2\cdot \Big((1-2x) - gh + \(1-h^2\)\Big).
\]

\noi In addition, after some simplification, we have
\[
(1+h)(1-2x+g) - (1-h)(1-2x-g) = 2 \cdot \Big((1-2x)h + g\Big)
\]
and
\[
(1+h)(1-2x-g) - (1-h)(1-2x+g) = 2 \cdot \Big((1-2x)h - g\Big).
\]

\noi Taking all this into account, we need to show that
\beqn
\label{der-zero-simple}
\Big((1-2x) + gh + \(1-h^2\)\Big)^{p-1} \cdot \Big((1-2x)h - g\Big) = \Big((1-2x) - gh + \(1-h^2\)\Big)\cdot \Big((1-2x)h + g\Big)^{p-1}
\eeqn

\noi It's not hard to verify that
\[
(1-2z) + gh + \(1-h^2\) \queq 2 \cdot \Big((1-x)^2 - x^{\frac{2p-2}{p}} (1-x)^{\frac{2}{p}}\Big),
\]
that
\[
(1-2z) - gh + \(1-h^2\) \queq 2 \cdot \Big((1-x)^2 - x^{\frac{2}{p}} (1-x)^{\frac{2p-2}{p}}\Big),
\]
that
\[
(1-2z)h + g \queq 2x(1-x) \cdot \(\(\frac{1-x}{x}\)^{\frac{p-1}{p}} - \(\frac{x}{1-x}\)^{\frac{p-1}{p}}\),
\]
and that
\[
(1-2z)h - g \queq 2x(1-x) \cdot \(\(\frac{1-x}{x}\)^{\frac{1}{p}} - \(\frac{x}{1-x}\)^{\frac{1}{p}}\).
\]

\noi We can now complete the proof of the fact that (\ref{der-zer}) holds. In fact, substituting the above and simplifying, it is not hard to see that both sides of (\ref{der-zero-simple}) are equal to
\[
2^{p/2}\cdot x^{\frac{2p-2}{p}} (1-x)^{\frac{3p-3}{p}} \cdot \((1-x)^{\frac{2p-2}{p}} - x^{\frac{2p-2}{p}}\)^{p-1} \cdot \((1-x)^{\frac{2}{p}} - x^{\frac{2}{p}}\).
\]

\noi This completes the proof of Lemma~\ref{lem:u-ast}

\eprf

\noi {\bf Proof of Lemma~\ref{lem:F-eq}}.

\noi We again reduce the claim to an algebraic identity involving the functions $g = g(z)$ and $h = h(z)$, which we proceed to verify. Recalling that $Q(u) = \frac12 \cdot \(\(\sqrt{1 - \rho u} + \sqrt{u}\)^p + \(\sqrt{1 - \rho u} - \sqrt{u}\)^p\)$, taking $u = u^{\ast} = \frac{s}{\rho n} = \frac{1-2z-g}{1-2z+g} \cdot \frac{1-h}{2}$, and recalling the definition of $\Phi(n,s,p)$, we need to verify the identity
\beqn
\label{F-is-good}
\rho^{p/2} \cdot \(\(\sqrt{1 - \rho u} + \sqrt{u}\)^p + \(\sqrt{1 - \rho u} - \sqrt{u}\)^p\) \queq \frac{n}{n-i} \cdot \(\frac{s}{n}\)^{p/2}\cdot \(1 + \frac{n-s}{s} \cdot t\)^p.
\eeqn

\noi We proceed by expressing everything via the functions $g$ and $h$, as above. It is not hard to see that the LHS of (\ref{F-is-good}) is equal to
\[
\frac{\Big(\sqrt{(1+h)(1-2z+g)} + \sqrt{(1-h)(1-2z-g)}\Big)^p + \Big(\sqrt{(1+h)(1-2z+g)} - \sqrt{(1-h)(1-2z-g)}\Big)^p}{2^{p/2} \cdot (1-2z-g)^{p/2}},
\]
and the RHS of (\ref{F-is-good}) equals to
\[
\frac{1}{1-z} \cdot \(\frac{1-h}{2}\)^{p/2} \cdot \(\frac{2 - 2z + g -h}{1-h}\)^p \queq \frac{1}{1-z} \cdot \frac{1}{2^{p/2} (1-h)^{p/2}} \cdot \(2 - 2z + g -h\)^p.
\]

\noi Rearranging, we need to show that
\[
\Big(\sqrt{(1+h)(1-2z+g)} + \sqrt{(1-h)(1-2z-g)}\Big)^p + \Big(\sqrt{(1+h)(1-2z+g)} - \sqrt{(1-h)(1-2z-g)}\Big)^p ~=
\]
\[
\frac{1}{1-z} \cdot \(\frac{1-2z-g}{1-h}\)^{p/2} \cdot \(2 - 2z + g -h\)^p.
\]

\noi We start with simplifying the LHS of this putative identity. Recall that
\[
\(\sqrt{(1+h)(1-2z+g)} + \sqrt{(1-h)(1-2z-g)}\)^2 ~=~ 2\cdot \Big((1-2z) + gh + \(1-h^2\)\Big).
\]

\noi Similarly, it is easy to see that
\[
\(\sqrt{(1+h)(1-2z+g)} - \sqrt{(1-h)(1-2z-g)}\)^2 ~=~ 2\cdot \Big((1-2z) + gh - \(1-h^2\)\Big).
\]

\noi Hence, the LHS is $2^{p/2}$ times
\[
\Big((1-2z) + gh + \(1-h^2\)\Big)^{p/2} + \Big((1-2z) + gh - \(1-h^2\)\Big)^{p/2}.
\]

\noi Recall that we have
\[
(1-2z) + gh + \(1-h^2\) \queq 2 \cdot \Big((1-z)^2 - z^{\frac{2p-2}{p}} (1-z)^{\frac{2}{p}}\Big)
\]

\noi Similarly, it is easy to see that
\[
(1-2z) + gh - \(1-h^2\) \queq 2 \cdot \Big(z(1-z) \(\frac{1-z}{z}\)^{\frac{p-2}{p}} - z^2\Big).
\]

\noi Substituting and simplifying, it's not hard to verify that the LHS is
\[
2^p \cdot \((1-z)^{\frac{2p-2}{p}} - z^{\frac{2p-2}{p}}\)^{p/2}.
\]

\noi The RHS is harder to simplify, but we can write it as
\[
\frac{1}{1-z} \cdot 2^p \cdot \(\frac{(1-2z-g)\(1 - z + \frac{g-h}{2}\)^2}{1-h}\)^{p/2}.
\]

\noi Simplifying and rearranging, veryfying that these two expressions are equal amounts to verifying that
\[
\frac{(1-2z-g)\(1 - z + \frac{g-h}{2}\)^2}{1-h} \queq (1-z)^2 - z^{\frac{2p-2}{p}} (1-z)^{\frac{2}{p}}.
\]

\noi Substituting the definitions of $g$ and $h$, this is equivalent to
\[
\((1-2z) - \(z^{\frac{1}{p}}(1-z)^{\frac{p-1}{p}} - z^{\frac{p-1}{p}}(1-z)^{\frac{1}{p}}\)\) \cdot \((1-z) - z^{\frac{p-1}{p}}(1-z)^{\frac{1}{p}}\)^2 ~=
\]
\[
\(1 - \(z^{\frac{1}{2k}}(1-z)^{\frac{p-1}{p}} + z^{\frac{p-1}{p}}(1-z)^{\frac{1}{p}}\)\) \cdot \((1-z)^2 - z^{\frac{2p-2}{p}} (1-z)^{\frac{2}{p}}\).
\]

\noi Writing $a = 1-z$ and $b = z^{\frac{p-1}{p}}(1-z)^{\frac{1}{p}}$, the second term on the LHS is $(a-b)^2$, while the second term on the RHS is $a^2 - b^2$. So we can divide out by $a-b$, and have to show that
\[
\((1-2z) - \(z^{\frac{1}{p}}(1-z)^{\frac{p-1}{p}} - z^{\frac{p-1}{p}}(1-z)^{\frac{1}{p}}\)\) \cdot \((1-z) - z^{\frac{p-1}{p}}(1-z)^{\frac{1}{p}}\) ~=
\]
\[
\(1 - \(z^{\frac{1}{p}}(1-z)^{\frac{p-1}{p}} + z^{\frac{p-1}{p}}(1-z)^{\frac{1}{p}}\)\) \cdot \((1-z) + z^{\frac{p-1}{p}} (1-z)^{\frac{1}{p}}\)
\]

\noi Write $A := z^{\frac{1}{p}}(1-z)^{\frac{p-1}{p}}$ and $B := z^{\frac{p-1}{p}}(1-z)^{\frac{1}{p}}$. Then the last identity is
\[
\Big((1-2z) - (A-B)\Big) \cdot \Big((1-z) - B\Big) ~=~ \Big(1 - (A+B)\Big) \cdot \Big((1-z) + B \Big)
\]

\noi Simplifying and rearranging, this is the same as $AB = z(1-z)$, which is true.

\eprf

\subsection{Proof of Proposition~\ref{pro:approximate}}

\noi Recall that we assume that $n$ is large and that $s_0 \le s \le \frac n2 - s_0$, where $s_0 = s_0(n) = \frac{n}{\ln n}$.

\noi We first observe that under this assumption, the value of $i_0$ given by (\ref{i_0}) is bounded away from $0$ and from $\frac n2$.

\lem
\label{lem:i-star}
Let $0 \le i_0 \le \frac n2$ be given by (\ref{i_0}). Then
\[
\(\frac{s_0}{n}\)^p \quad \le \quad \frac{i_0}{n} \quad \le \quad \(\frac12 - \sqrt{\frac{s}{n}\(1-\frac{s}{n}\)}\) - \Omega\(\(\frac{p-2}{p}\)^2 \cdot \(\frac{s_0}{n}\)^4\)
\]
Here the asymptotic notation hides absolute factors.
\elem

\prf

\noi Let $x = \frac{i_0}{n}$, $y = \frac12 - \sqrt{\frac{s}{n}\(1-\frac{s}{n}\)}$. We start with the first inequality, since it is easier. The derivative of $h$ is computed in the proof of Lemma~\ref{lem:ref:h} below, and it is easy to see that for all $0 < z \le 1/2$  holds $h'(z) \le \frac{p-1}{p} \cdot z^{-\frac{1}{p}} + \frac{1}{p} \cdot z^{-\frac{p-1}{p}} \le z^{-\frac{1}{p}} + z^{-\frac{p-1}{p}}$. Since $h(0) = 0$, it follows that $h(z) \le  z^{\frac{p-1}{p}} + z^{\frac{1}{p}} \le 2 z^{\frac{1}{p}}$. Hence $h(x) = 1 - \frac{2s}{n} \ge \frac{2s_0}{n}$ implies $x \ge \(\frac{s_0}{n}\)^p$, completing the first inequality.

\noi We pass to the second inequality. Recall that $h^2(z) = g^2(z) + 4z(1-z)$. Hence we have, observing that $4y(1-y) = \(1 - \frac{2s}{n}\)^2$, that
\[
h(y) ~=~ \sqrt{4y(1-y) + g^2(y)} ~=~ \sqrt{\(1 - \frac{2s}{n}\)^2 + g^2(y)} ~\ge~ \(1 - \frac{2s}{n}\) + \frac{g^2(y)}{4},
\]
where the last inequality follows from the following easily verifiable clam: Let $0 \le a, \e \le 1$. Then $\sqrt{a^2 + \e} \ge a + \frac{\e}{4}$.

\noi Next, we have
\[
g(y) = \(y(1-y)\)^{\frac{1}{p}} \cdot \((1-y)^{\frac{p-2}{p}} - y^{\frac{p-2}{p}}\) ~\ge~ \(\frac 14 \cdot \(1 - \frac{2s}{n}\)^2\)^{\frac{1}{p}} \cdot \frac{2p-4}{p} \cdot \sqrt{\frac{s}{n}\(1-\frac{s}{n}\)} ~\ge
\]
\[
\(\frac 14 \cdot \(1 - \frac{2s}{n}\)^2\)^{\frac{1}{2}} \cdot \frac{2p-4}{p} \cdot \sqrt{\frac{s}{n}\(1-\frac{s}{n}\)} ~=~ \frac{p-2}{p} \cdot \(1 - \frac{2s}{n}\) \cdot \sqrt{\frac{s}{n}\(1-\frac{s}{n}\)}.
\]
To see the first inequality, observe that for $0 \le t \le 1$ holds $\(t^{\frac{p-2}{p}}\)' \ge \frac{p-2}{p}$, and hence that $(1-y)^{\frac{p-2}{p}} - y^{\frac{p-2}{p}} \ge \frac{p-2}{p} \cdot (1-2y)$. It follows that
\[
h(y) ~\ge~ \(1 - \frac{2s}{n}\) + \(\frac{p-2}{2p}\)^2 \cdot \(1 - \frac{2s}{n}\)^2 \cdot \frac{s}{n}\(1-\frac{s}{n}\).
\]

\noi Let $\delta = \(\frac{p-2}{2p}\)^2 \cdot \(1 - \frac{2s}{n}\)^2 \cdot \frac{s}{n}\(1-\frac{s}{n}\)$. Recall that by definition $h(x) = 1 - \frac{2s}{n}$. Hence we get $h(y) \ge h(x) + \delta$. Computing $h'$ as in the proof of Lemma~\ref{lem:ref:h}, it is easy to see that for all $0 < z < 1$ holds $0 < h'(z) < \frac{1}{z}$. From this, $\delta \le h(y) - h(x) \le \frac{y-x}{x}$, which implies to $x \le y -  \frac{\delta}{1+\delta} \cdot y \le y - \frac{\delta}{2} \cdot y$ (the last inequality follows since clearly $\delta \le \frac12$).

\noi Recalling that $y = \frac12 - \sqrt{\frac{s}{n}\(1-\frac{s}{n}\)}$, and that by assumption $s_0 \le s \le \frac n2 - s_0$, it is easy to check that $\frac{\delta}{2} \cdot y \ge \Omega\(\(\frac{p-2}{p}\)^2 \cdot \(\frac{s_0}{n}\)^4\)$, and the claim of the lemma holds.

\eprf

\rem
\label{rem:i_0-integer}
We will assume from now on, to avoid complications in notation arising from replacing $i_0$ by the nearest integer, that $i_0$ is integer, whenever it is convenient for us to do so. It is easy to see that the error this introduces is negligible.
\erem

\noi The key step in the proof of Proposition~\ref{pro:approximate} is the following claim, which may be of independent interest.

\pro
\label{pro:2k-norm}
Let $p > 2$ be fixed. Let $0 < i_0 < \frac n2$ satisfy (\ref{i_0}). Then the $\ell_p$ norm of $K_s$ is attained, up to a small error, in a union of intervals of length $O\(\sqrt{n \log n}\)$ around $i_0$ nad $n - i_0$. More precisely, there is an absolute constant $C$ such that if $I$ is the interval of length $C \sqrt{n \log n}$ around $i_0$ then, for a sufficiently large $n$, depending on $p$, and for any $\sigma \le s \le n/2 - \sigma$ holds
\[
\frac{1}{2^n} \sum_{i \in I \cup (n - I)} {n \choose i} \(|K_s(i)|\)^p ~\ge~ \(1 - O\(\frac{1}{n^2}\)\) \cdot ||K_s||^p_p.
\]
\epro

\noi This proposition and the argument leading towards its proof will have the following corollary as in easy implication. We write a superscript $K_s^{(n)}$ for the Krawchouk polynomial on the $n$-dimensional cube, when we consider functions on cubes of different dimensions (in the second claim of the corollary).

\cor
\label{cor:i_0-approx}
Let $i_0$ be given by (\ref{i_0}). Then
\begin{enumerate}
\item
\[
\frac{||K_s||^p_p}{||K_{s-1}||^p_p} ~\in~ \(1 \pm O(\e)\)^p \cdot \frac{K^{p}_s\(i_0\)}{K^{p}_{s-1}\(i_0\)}
\]

\item
\[
\frac{||K^{(n+1)}_s||^p_p}{||K^{(n)}_{s-1}||^p_p} ~\in~ \(1 \pm O(\e)\)^p \cdot \frac{\frac{{{n+1} \choose i_0}}{2^{n+1}} \cdot \(K^{(n+1)}_s\(i_0\)\)^p}{\frac{{{n} \choose i_0}}{2^{n}} \cdot \(K^{(n)}_{s-1}\(i_0\)\)^p}
\]

\end{enumerate}

\ecor

\noi Looking ahead, the first claim of Proposition~\ref{pro:approximate} will be a simple consequence of the first claim of this corollary, and the second claim of the proposition  will follow easily from the second claim of the corollary. We will prove Proposition~\ref{pro:2k-norm} and Corollary~\ref{cor:i_0-approx} and, following this, complete the proof of Proposition~\ref{pro:approximate}.

\noi We proceed with the proof of Proposition~\ref{pro:2k-norm}.

\lem
\label{lem:ratios-k-norm}
Let $0 < i_0 < \frac n2$ be given by (\ref{i_0}). Let $n^{\frac23} \ll \Delta \ll \frac n2 - \sqrt{s(n-s)} - i_0$. Let $i_1 = \frac n2 - \sqrt{s(n-s)} - \Delta$. Then for any $0 \le i \le i_0$ holds
\[
\frac{{n \choose {i+1}} K^p_s(i+1)}{{n \choose {i}} K^p_s(i)} ~\ge~ \frac{i_0}{i+1} \cdot \(1 - O\(\frac{s}{\Delta^2}\)\)^p,
\]
for any $i_0 < i < i_1$ holds
\[
\frac{{n \choose {i+1}} K^p_s(i+1)}{{n \choose {i}} K^p_s(i)} ~\le~ \frac{i_0}{i+1} \cdot \(1 +O\(\frac{s}{\Delta^2}\)\)^p,
\]
and for any $i_1 \le i < x_s - 1$ holds
\[
\frac{{n \choose {i+1}} K^p_s(i+1)}{{n \choose {i}} K^p_s(i)} ~\le~ \frac{i_0}{i_1+1} \cdot \(1 +O\(\frac{s}{\Delta^2}\)\)^p.
\]
\elem
\prf

\noi We will need the two following facts: The location of the first root $x_s$ of the Krawchouk polynomial $K_s$ and the behaviour of the values of $K_s$ in the interval $\(0, x_s\)$. Recall that (see e.g., \cite{KZ}) we have
\beqn
\label{first-root}
\frac n2 - \sqrt{s(n-s)}  ~\le~ x_s ~\le~ \frac n2 - \sqrt{s(n-s)} + O\(n^{\frac23}\).
\eeqn

\noi The following fact has been shown in \cite{KL, MRRW} (see also proof of Lemma~5.1 in \cite{S} for a more detailed calculation). Let $0 \le i \le x_s - \Delta$ for some $\Delta \gg \sqrt{s}$. Then
\beqn
\label{ratio-estimate}
\frac{K_s(i+1)}{K_s(i)} ~\in~\(1, 1 \pm O\(\frac{s}{\Delta^2}\)\) \cdot \frac{(n-2s) + \sqrt{(n-2s)^2 - 4i(n-i)}}{2(n-i)}
\eeqn

\noi We proceed with the proof. We will prove the first inequality. The second is similar. The third follows from the second immediately, since the ratio $\frac{{n \choose {i+1}} K^p_s(i+1)}{{n \choose {i}} K^p_s(i)}$ decreases in $i$ on $0 \le i \le x_s - 1$ (to see this note that clearly the ratio of the binomial coefficients decreases, and as observed in the proof of the preceding lemma, the ratio $\frac{K_s(i+1)}{K_s(i)}$ decreases as well).

\noi Observe that (\ref{i_0}) means that $\frac{(n-2s) + \sqrt{(n-2s)^2 - 4i_0\(n-i_0\)}}{2\(n-i_0\)} = \(\frac{i_0}{n-i_0}\)^{\frac{1}{p}}$. To see this, note that, as in the proof of Lemmas~\ref{lem:u-ast}~and~\ref{lem:F-eq} above, we can rewrite this equaity in terms of the variable $x = \frac{i_0}{n}$ and the functions $g$ and $h$ of this variable. It transforms into $\frac{h(x) + \sqrt{h^2(x) - 4x(1-x)}}{2(1-x)} = \(\frac{x}{1-x}\)^{\frac{1}{p}}$, which is the same as $h + g = 2x^{\frac{1}{p}} (1-x)^{\frac{p-1}{p}}$, and this follows directly from the definition of $g$ and of $h$.

\noi Hence, we have (using (\ref{ratio-estimate}) and Lemma~\ref{lem:i-star}) that
\[
\frac{{n \choose {i+1}} K^p_s(i+1)}{{n \choose {i}} K^p_s(i)} ~=~ \frac{n-i}{i+1} \(\frac{K_s(i+1)}{K_s(i)}\)^{p} ~\ge~ \frac{n-i}{i+1} \(\frac{K_s\(i_0+1\)}{K_s\(i_0\)}\)^{p} ~\ge~
\]
\[
\(1 - O\(\frac{s}{\Delta^2}\)\)^p \cdot \frac{n-i}{i+1} \cdot \(\frac{(n-2s) + \sqrt{(n-2s)^2 - 4i_0\(n-i_0\)}}{2\(n-i_0\)}\)^p ~=~
\]
\[
\(1 - O\(\frac{s}{\Delta^2}\)\)^p \cdot \frac{n-i}{i+1} \cdot \frac{i_0}{n-i_0} ~\ge~ \(1 -O\(\frac{s}{\Delta^2}\)\)^p \cdot \frac{i_0}{i+1}.
\]
\eprf

\lem
\label{lem:2-increase}
Let $0 \le i \le x_s - \Delta$ for some $\Delta \gg \frac{n^{\frac56}}{s_0^{\frac16}}$. Then
\[
\frac{{n \choose {i+1}} K^2_s(i+1)}{{n \choose i} K^2_s(i)} ~\ge~ 1 + \Omega\(\frac{\Delta s_0^{\frac12}}{n^{\frac32}}\).
\]
\elem
\prf
We have, by (\ref{first-root}) and (\ref{ratio-estimate}) that
\[
\frac{{n \choose {i+1}} K^2_s(i+1)}{{n \choose i} K^2_s(i)} ~=~ \frac{n-i}{i+1} \cdot \(\frac{K_s(i+1)}{K_s(i)}\)^2 ~\ge~
\]
\[
\(1 - O\(\frac{s}{\Delta^2}\)\)^2 \cdot \frac{n-i}{i+1} \cdot \(\frac{(n-2s) + \sqrt{(n-2s)^2 - 4i(n-i)}}{2(n-i)}\)^2 ~\ge
\]
\[
\(1 - O\(\frac{s}{\Delta^2}\)\)^2 \cdot \frac{(n-2s)^2}{4(n-i)(i+1)}
\]
The quadratic $Q(x) = 4x(n-x)$ equals $(n-2s)^2$ at $x = \frac n2 - \sqrt{s(n-s)}$. It is easy to see that this means that for $i < \frac n2 - \sqrt{s(n-s)} - \Omega(\Delta)$ (this estimate on $i$ is valid by (\ref{first-root}) and by our assumption on $\Delta$), we have $4(i+1)(n-i) \le (n-2s)^2 - \Omega\(\Delta \cdot (n-2i)\) \le (n-2s)^2 - \Omega\(\Delta \cdot \sqrt{s(n-s)}\)$. Hence
\[
\frac{(n-2s)^2}{4(n-i)(i+1)} ~\ge~ 1 + \Omega\(\frac{\Delta \cdot \sqrt{s(n-s)}}{(n-2s)^2}\) ~\ge~ 1 + \Omega\(\frac{\Delta s_0^{\frac12}}{n^{\frac32}}\).
\]

\eprf

\noi We can now complete the proof of Proposition~\ref{pro:2k-norm}, using the three auxiliary claims above.

\noi {\bf Proof of Proposition~\ref{pro:2k-norm}}.

\noi Let $\Delta = n^{\frac45}$. By Corollary~\ref{cor:norm-any-two-roots} and by Lemma~\ref{lem:2-increase}, there exists $x_s - \Delta \le i' \le x_s$ such that, say, ${n \choose i'} K^2_s\(i'\) \ge \frac{1}{n^3} \cdot 2^n {n \choose s}$, which means that $K_s\(i'\) \ge \sqrt{\frac{1}{n^3} \cdot \frac{2^n {n \choose s}}{{n \choose i'}}}$. Let $i_1 = \lfloor x_s - \Delta \rfloor$. Let $D = i_1 - i_0$. Note that by Lemma~\ref{lem:i-star} we have $D \ge \Omega\(\frac{s_0^4}{n^3}\)$, and hence by the third claim of Lemma~\ref{lem:ratios-k-norm} we have that
\[
{n \choose {i_1}} K_s^p\(i_1\) ~\ge~ {n \choose {i'}} K_s^p\(i'\) ~\ge~ \frac{1}{n^{3p/2}} \cdot \frac{2^{pn/2} {n \choose s}^{p/2}}{{n \choose i_1}^{p/2-1}}.
\]

\noi By Lemma~\ref{lem:ratios-k-norm},
\[
\frac{{n \choose i_0} K^p_s\(i_0\)}{{n \choose i_1} K^p_s\(i_1\)} ~\ge~ \(1 +O\(\frac{s}{\Delta^2}\)\)^{-pD} \cdot \frac{\(i_0+1\)\(i_0 + 2\)...\(i_0+D\)}{i^D_0} ~\ge~
\]
\[
e^{-\frac{pDn}{\Delta^2}} \cdot e^{\Omega\(\frac{D^2}{n}\)} ~\ge~ e^{\Omega\(\frac{D^2}{n}\)} ~\ge~ e^{\Omega\(\frac{s_0^8}{n^7}\)}.
\]
This implies that
\[
{n \choose i_0} K^p_s\(i_0\) ~\ge~ e^{\Omega\(\frac{\sigma^8}{n^7}\)} \cdot {n \choose i_1} K^p_s\(i_1\) ~\ge~
\]
\[
e^{\Omega\(\frac{s_0^8}{n^7}\)} \cdot \frac{1}{n^{3p/2}} \cdot \frac{2^{pn/2} {n \choose s}^{p/2}}{{n \choose i_1}^{p/2-1}} ~\ge~
e^{\Omega\(\frac{s_0^8}{n^7}\)} \cdot  \frac{2^{pn/2} {n \choose s}^{p/2}}{{n \choose i_1}^{p/2-1}}.
\]

\noi Next, for any $i_1 \le i \le n/2$ holds ${n \choose i} K^2_s(i) \le 2^n ||K_s||^2_2 = 2^n {n \choose s}$, and hence ${n \choose i} \(|K_s|(i)\)^p \le \frac{2^{pn/2} {n \choose s}^{p/2}}{{n \choose i}^{p/2-1}} \le \frac{2^{pn/2} {n \choose s}^{p/2}}{{n \choose i_1}^{p/2-1}}$. This implies that
\[
{n \choose i_0} K^p_s\(i_0\) ~\ge~ e^{\Omega\(\frac{s_0^8}{n^7}\)} \cdot \(\sum_{i_1 \le i \le n/2} {n \choose i} \(|K_s|(i)\)^p\).
\]

\noi In addition, by Lemma~\ref{lem:ratios-k-norm}, similarly to the above, we have that for a sufficiently large constant $C$ holds that if $0 \le i < i_1$ and $|i-i_0| \ge C \cdot \sqrt{n \log n}$ then $\frac{{n \choose i_0} K^p_s\(i_0\)}{{n \choose i} K^p_s\(i\)} \ge n^3$. Let $I$ be the interval $\left[i_0 -  C \cdot \sqrt{n \log n}, i_0 + C \cdot \sqrt{n \log n} \right]$. Then by the above, we have $\sum_{i \not \in I, i \le n/2} {n \choose i} \(|K_s|(i)\)^p \le O\(\frac{1}{n^2}\) \cdot {n \choose i_0} K^p_s\(i_0\)$.

\noi Finally, recall that $K_s$ is symmetric around $n/2$ if $s$ is even, and antisymmetric if $s$ is odd. Taking all of this into account, we have
\[
\sum_{i \in I \cup (n - I)} {n \choose i} \(|K_s(i)|\)^p ~\ge~ \(1 - O\(\frac{1}{n^2}\)\) \cdot ||K_s||^p_p,
\]
completing the proof of the proposition.

\eprf

\noi {\bf Proof of Corollary~\ref{cor:i_0-approx}}.

\noi We start with the first part of the corollary. By Proposition~\ref{pro:2k-norm} applied to both $K_{s-1}$ and $K_s$ there is an interval $I$ of length $O\(\sqrt{n \log n}\)$ around $i_0$, where $i_0$ is defined by (\ref{i_0}) such that both the $\ell_p$ norms of $K_{s-1}$ and $K_s$ are attained, up to a factor of $1-O\(\frac{1}{n^2}\)$ on $I$ and on $n - I$. Taking into account the symmetry (or anti-symmetry) of the Krawchouk polynomials around $\frac n2$, the claim of the corollary will follow if we show that for all $i$ in $I$ holds $\frac{K_s(i)}{K_{s-1}(i)} \in \(1 \pm O(\e)\) \cdot \frac{K_s\(i_0\)}{K_{s-1}\(i_0\)}$.

\noi Recall that ${n \choose i} K_s(i) = {n \choose s}K_i(s)$. Hence $\frac{K_s(i)}{K_{s-1}(i)} = \frac{{n \choose s}}{{n \choose {s-1}}} \cdot \frac{K_i(s)}{K_{i}(s-1)}$.

\noi By the above discussion we know that each $i$ in $I$ satisfies $i \le x_s - \Omega\(\frac{s_0^4}{n^3}\)$. An easy calculation\footnote{We omit the details.} using (\ref{first-root}) shows that $i \le x_s - \Omega\(\frac{\sigma^4}{n^3}\)$ implies $s \le x_{i} - \Omega\(\frac{s_0^{4.5}}{n^{3.5}}\)$. Hence we may apply (\ref{ratio-estimate}), with roles of $s$ and $i$ reversed, to obtain $\frac{K_i(s)}{K_{i}(s-1)} \in \(1 \pm O\(\frac{n^8}{s_0^9}\)\) \cdot \frac{(n-2i) + \sqrt{(n-2i)^2 - 4s(n-s)}}{2(n-s)}$.

\noi This means that
\[
\frac{K_s(i)}{K_{s-1}(i)} ~/~ \frac{K_s\(i_0\)}{K_{s-1}\(i_0\)} ~\in~ \(1 \pm O\(\frac{n^8}{s_0^9}\)\) \cdot \frac{(n-2i) + \sqrt{(n-2i)^2 - 4s(n-s)}}{\(n-2i_0\) + \sqrt{\(n-2i_0\)^2 - 4s(n-s)}}.
\]
By Lemma~\ref{lem:i-star}, both $n - 2i$ and $n - 2i_0$ are lowerbounded by $\Omega\(\frac{s_0^5}{n^4}\)$. Hence $\frac{n-2i}{n - 2i_0} \in 1 \pm \frac{n^{9/2} \log^{1/2}(n)}{s_0^5}$. In addition, since both $i$ and $i_0$ are upperbounded by $x_s - \Omega\(\frac{s_0^4}{n^3}\) = \frac n2 - \sqrt{s(n-s)} - \Omega\(\frac{s_0^4}{n^3}\)$, it is easy to see that $\frac{(n-2i)^2 - 4s(n-s)}{\(n-2i_0\)^2 - 4s(n-s)} ~\in~ 1 \pm \frac{n^4 \log^{1/2}(n)}{s_0^{9/2}}$. Recalling the definition of $\e = \e(n)$ we see that both ratios lie in $1 \pm \e$ and hence
\[
\frac{K_s(i)}{K_{s-1}(i)} ~/~ \frac{K_s\(i_0\)}{K_{s-1}\(i_0\)} ~\in~ \(1 \pm O\(\frac{n^7}{\sigma^8}\)\) \cdot \(1 \pm O(\e)\) ~\subseteq~ 1 \pm O(\e).
\]

\noi We pass to the second part of the corollary. Again, we may focus our attention on an interval $I$ of length $O\(\sqrt{n \log n}\)$ around $i_0$ given by (\ref{i_0}) in which (half of) the $\ell_p$ norms of $K^{(n)}_{s-1}$, $K^{(n)}_{s}$, and $K^{(n+1)}_s$ are attained, up to a factor of $1-\frac{1}{n^2}$. Moreover, following the argument above, for all $i \in I$ the ratio $\frac{K^{(n)}_s(i)}{K^{(n)}_{s-1}(i)}$ is in $\(1 \pm O(\e)\) \cdot \frac{K^{(n)}_s\(i_0\)}{K^{(n)}_{s-1}\(i_0\)}$. Recall the identity (see \cite{Lev}, but also easy to verify directly) $K^{(n+1)}_s(i) = K^{(n)}_s(i) + K^{(n)}_{s-1}(i)$. This means that for all $i \in I$ holds
\[
\frac{K^{(n+1)}_s(i)}{K^{(n)}_{s-1}(i)} ~=~ \frac{K^{(n)}_s(i)}{K^{(n)}_{s-1}(i)} + 1 ~\in~ \(1 \pm O(\e)\) \cdot \(\frac{K^{(n)}_s\(i_0\)}{K^{(n)}_{s-1}\(i_0\)} + 1\) ~=~
\(1 \pm O(\e)\) \cdot  \frac{K^{(n+1)}_s\(i_0\)}{K^{(n)}_{s-1}\(i_0\)}.
\]

\eprf

\noi {\bf Proof of Proposition~\ref{pro:approximate}}

\noi We start with the first claim of the proposition. Let $\rho_1 = \(\frac{r(n,s,p)}{r(n,s-1,p)}\)^{2/p} = \frac{||K_s||^2_p}{\E \(K_{s}\)^2} ~/~ \frac{||K_{s-1}||^2_p}{\E \(K_{s-1}\)^2}$. We need to show that $\rho_1 \in \(1 \pm O(\e)\) \cdot \rho$, where $\rho = \frac{n-2i_0}{s} \cdot t - 1$, and $t = \frac{\(n-2i_0\) + \sqrt{\(n-2i_0\)^2 - 4s(n-s)}}{2(n-s)}$.

\noi Recall that $\E K^2_s = {n \choose s}$ and that ${n \choose i} K_s(i) = {n \choose s}K_i(s)$. Applying the first claim of Corollary~\ref{cor:i_0-approx} and using (\ref{ratio-estimate}) (with roles of $s$ and $i$ reversed) we have that (estimating $\frac{K_{i_0}\(s\)}{K_{i_0}\(s-1\)}$ as in the proof of Corollary~\ref{cor:i_0-approx}):
\[
\rho_1 ~\in~ \(1 \pm O\(\e\)\)^{2} \cdot \frac{{n \choose {s-1}}}{{n \choose s}} \cdot \(\frac{K_s\(i_0\)}{K_{s-1}\(i_0\)}\)^2 ~=~
\(1 \pm O\(\e\)\)^{2} \cdot \frac{{n \choose s}} {{n \choose {s-1}}} \cdot \(\frac{K_{i_0}\(s\)}{K_{i_0}\(s-1\)}\)^2 ~\subseteq~
\]
\[
\(1 \pm O\(\e\)\)^{2}\cdot \(1 \pm O\(\frac{n^7}{\sigma^8}\)\) \cdot \frac{n-s+1}{s} \cdot t^2(n,s) ~\subseteq~
\(1 \pm O\(\e\)\) \cdot \frac{n-s}{s} \cdot t^2(n,s)
\]

\noi Finally, recall that $t = t(n,s)$ is a root of the quadratic $(n-s)t^2 - \(n-2i_0\)t + s = 0$. Hence $t^2 = \frac{\(n-2i_0\)t-s}{n-s}$. Substituting this in the above expression and simplifying gives the first claim of the proposition.

\noi We proceed to the second claim of the proposition. We need to show that $\frac{r(n+1,s,p)}{r(n-1,s-1,p)} ~\in~ \(1 \pm O\(\e\)\)^p \cdot \Phi(n,s,p)$, where $\Phi(n,s,p) = \frac{n}{2\(n-i_0\)} \cdot \(\frac{s}{n}\)^{p/2} \cdot \(1 + \frac{n-s}{s} \cdot t\)^p$.

\noi We have that
$\frac{r(n+1,s,p)}{r(n-1,s-1,p)} = \frac{||K^{(n+1)}_s||^p_p}{\E^{p/2} \(K^{(n+1)}_{s}\)^2} ~/~ \frac{||K^{(n)}_{s-1}||^p_p}{\E^{p/2} \(K^{(n)}_{s-1}\)^2}$. Applying the second claim of Corollary~\ref{cor:i_0-approx} (and replacing, within negligible error, $\frac{n+1}{2\(n+1-i_0\)}$ by $\frac{n}{2\(n-i_0\)}$) we have that the RHS of this expression is in
\[
\(1 \pm O\(\e\)\)^p \cdot  \frac{n}{2\(n-i_0\)}  \(\frac{{n \choose {s-1}}}{{{n+1} \choose s}}\)^{p/2}
\(\frac{K^{(n+1)}_s\(i_0\)}{K^{(n)}_{s-1}\(i_0\)}\)^p ~\subseteq~
\]
\[
\(1 \pm O\(\e\)\)^p \cdot  \frac{n}{2\(n-i_0\)}  \(\frac{s}{n+1}\)^{p/2} \(1 + \frac{K^{(n)}_s\(i_0\)}{K^{(n)}_{s-1}\(i_0\)}\)^p.
\]

\noi Recalling that
\[
\frac{K^{(n)}_s\(i_0\)}{K^{(n)}_{s-1}\(i_0\)} = \frac{{n \choose s}}{{n \choose {s-1}}} \cdot \frac{K^{(n)}_{i_0}(s)}{K^{(n)}_{i_0}(s-1)} = \frac{n-s+1}{s} \cdot \frac{K^{(n)}_{i_0}(s)}{K^{(n)}_{i_0}(s-1)} \in \(1 \pm O\(\frac{n^7}{s_0^8}\)\) \cdot \frac{n-s}{s} \cdot t,
\]
and replacing, within negligible error, $\frac{s}{n+1}$ by $\frac{s}{n}$, we obtain the second claim of the proposition.

\eprf

\subsection{Proofs of Lemma~\ref{lem:psi} and Proposition~\ref{pro:gap}}
\label{sec:lemmas}

\subsubsection*{Proof of Lemma~\ref{lem:psi}}

\noi Let $n$ and $p \ge 2$ be fixed. Let $0 \le s \le \frac n2$. Let $m$ be an integer, and let $N = nm$ and $S = sm$. We need to show that
\[
\lim_{m \rarrow \infty} \Big(r(N,S,p)\Big)^{\frac 1m} ~=~ \lim_{m \rarrow \infty} \(\frac{\E \(K^{(N)}_S\)^p}{\E^{p/2} \(K^{(N)}_S\)^2}\)^{\frac 1m} ~=~ 2^{n \cdot\psi\(p,\frac sn\)},
\]
For $s = 0$ the claim of the lemma reduces to verifying that $\psi(p,0) = 0$, and for $s = \frac n2$ to verifying that $\psi\(p,\frac12\) = \frac{p-2}{2}$. Both these facts follow easily from the definition of $\psi$. So we may assume from now on that $0 < s < \frac n2$.

\noi Consider first the denominator. Recalling that $\E \(K^{(N)}_S\)^2 = {N \choose S}$ and using (\ref{binomial-H}), we have that
\[
\lim_{m \rarrow \infty} \({\E}^{p/2} \(K^{(N)}_S\)^2\)^{\frac 1m} ~=~ \lim_{m \rarrow \infty} {N \choose S}^\frac{p}{2m}  ~=~ \lim_{m \rarrow \infty} 2^{\frac{pN}{2m} H\(\frac SN\)} ~=~   2^{\frac12 pn H\(\frac sn\)}.
\]

\noi Next, consider the numerator. Recall that $s_0 = s_0(N) = \frac{N}{\ln N}$. For a sufficiently large $m$, $S = sm$ satisfies $s_0 < S < \frac N2 - s_0$. Hence, by Proposition~\ref{pro:2k-norm}, we have that, up to a constant factor, the value of $\E \(K^{(N)}_S\)^p$ is given by $\frac{1}{2^N} \sum_{i \in I} {N \choose i} \(K^{(N)}_S(i)\)^p$, where $I$ is an interval of length $O\(\sqrt{N \log N}\)$ around $i_0$, and $i_0$ is determined by  $h\(p,\frac{i_0}{N}\) = 1 - \frac{2S}{N} = 1 - \frac{2s}{n}$. Let $i_1$ be the leftmost integer point of $I$ and let $i_2$ be the rightmost point. Then for any $i \in I$ holds ${N \choose i} \(K^{(N)}_S(i)\)^p \le {N \choose i_2} \(K^{(N)}_S\(i_1\)\)^p$, since the binomial coefficients increase as $i$ increases in $I$, while the value of $K^{(N)}_S(i)$ decreases. Next, by the Lemma~\ref{lem:2-increase} (see also Remark~\ref{rem:i_0-integer}), we have that ${N \choose {i_0}} \(K^{(N)}_S\)^2\(i_0\) \ge {N \choose {i_1}} \(K^{(N)}_S\)^2\(i_1\)$. Hence $K^{(N)}_S\(i_1\) \le \sqrt{\frac{{N \choose {i_0}}}{{N \choose {i_1}}}} K^{(N)}_S\(i_0\)$. Altogether, we have,
\[
\frac{{N \choose {i_0}}}{2^N} \(K^{(N)}_S\(i_0\)\)^p ~~\le~~ \E \(K^{(N)}_S\)^p ~~\le~~ O\(|I| \cdot \frac{{N \choose i_2}}{2^N} \cdot \(\frac{{N \choose {\lceil i_0 \rceil}}}{{N \choose {i_1}}}\)^{p/2} \(K^{(N)}_S\(i_0\)\)^p\).
\]

\noi Taking the limit as $m$ goes to infinity, and using the approximation of the binomial coefficient ${b \choose a}$ by $2^{b H\(a/b\)}$, we have
\[
\lim_{m \rarrow \infty} \(\E \(K^{(N)}_S\)^p\)^{\frac 1m} ~=~ \lim_{m \rarrow \infty} \(\frac{{N \choose {i_0}}}{2^N} \(K^{(N)}_S\(i_0\)\)^p\)^{\frac 1m} ~=~
2^{(H(y)-1)n} \cdot \lim_{m \rarrow \infty} \(K^{(N)}_S\(i_0\)\)^{\frac{p}{m}},
\]
where we write $y$ for $\frac{i_0}{N}$, remembering that $y$ is determined by $h(p,y) = 1 - \frac{2s}{n}$.

\noi Note that $i_0 \le \frac N2 - \sqrt{S(N-S)}$, and hence we may apply (\ref{Krawchouk-tail}) to obtain $K^{(N)}_S\(i_0\) ~~\in~~ 2^{\(\tau\(\frac SN,y\) \pm o(1)\) \cdot N} = 2^{\(\tau\(\frac sn, y\) \pm o_N(1)\) \cdot N}$. This gives $\lim_{m \rarrow \infty} \(K^{(N)}_S\(i_0\)\)^{\frac{p}{m}} = 2^{\(p \cdot \tau\(\frac sn,y\)\) \cdot n}$, and $\lim_{m \rarrow \infty} \(\E \(K^{(N)}_S\)^p\)^{\frac 1m} = 2^{\(H(y) - 1 + p \cdot \tau\(\frac sn,y\)\) \cdot n}$.

\noi Summing up, we have
\[
\lim_{m \rarrow \infty} \Big(r(N,S,p)\Big)^{\frac 1m} ~~=~~ 2^{\(H(y) - 1 + p \cdot \tau\(\frac sn,y\) - \frac p2 H\(\frac sn\)\) \cdot n} = 2^{\psi\(p,\frac sn\) \cdot n},
\]
where in the last step we use the first definition of $\psi$ in Section~\ref{subsubsec:ref:psi}.

\eprf

\subsubsection*{Proof of Proposition~\ref{pro:gap}}

\noi We will need a simple technical lemma.

\lem
\label{lem:H}
\[
\max_{0 \le i \le n/2} \frac{2^{n H\(\frac{i}{n}\)}}{{n \choose {\lfloor i \rfloor}}} ~\le~ O(n).
\]
\elem
\prf

\noi For $0 \le i < 1$ we have $\frac{2^{n H\(\frac{i}{n}\)}}{{n \choose {\lfloor i \rfloor}}} \le 2^{n H\(\frac{1}{n}\)} \le 2^{n \cdot \(\frac{1}{n} \log_2(n) + O\(\frac{1}{n}\)\)} \le O(n)$.

\noi Using (\ref{binomial-H}) we have that for $1 \le i \le n/2$ holds
\[
\frac{2^{n H\(\frac{i}{n}\)}}{{n \choose {\lfloor i \rfloor}}} ~\le~ O\(\sqrt{i} \cdot 2^{n \cdot \(H\(\frac{i}{n}\) - H\(\frac{\lfloor i \rfloor}{n}\)\)}\) ~\le~ O\(\sqrt{i} \cdot 2^{n \cdot \frac{i - \lfloor i \rfloor}{n} \cdot H'\(\frac{\lfloor i \rfloor}{n}\)}\),
\]
where the inequality follows from the concavity of $H$. Recalling that $H'(x) = \log_2\(\frac{1-x}{x}\) \le \log_2\(\frac{1}{x}\)$, the last expression is at most $O\(\sqrt{i} \cdot \frac{n}{\lfloor i \rfloor}\) \le O(n)$.

\eprf

\noi We proceed with the proof of the proposition. We have that $\psi\(p,\frac sn\) = H(y) - 1 + p \cdot \tau\(\frac sn,y\) - \frac p2 H\(\frac sn\)$, where $y$ is determined by $h(p,y) = 1 - 2\frac{s}{n}$. Set $i_0 = \lfloor ny \rfloor$ and observe that $i_0  \le ny \le \frac n2 - \sqrt{s(n-s)}$. Hence, by (\ref{Krawchouk-tail}), we have
\[
K_s\(i_0\) ~\ge~ \frac{{n \choose s}}{2^{H\(\frac sn\) \cdot n}} \cdot 2^{\tau\(\frac sn, \frac{i_0}{n}\) \cdot n} ~\ge~ \frac{{n \choose s}}{2^{H\(\frac sn\) \cdot n}} \cdot 2^{\tau\(\frac sn,y\) \cdot n},
\]
where in the second step we have used the fact that $\tau$ is decreasing in the second variable.

\noi Hence, using Lemma~\ref{lem:H} in the third inequality and (\ref{binomial-H}) in the last inequality, we have
\[
r(n,s,p) ~=~ \frac{\E |K_s|^p}{\E^{p/2} K^2_s} ~=~ \frac{\E |K_s|^p}{{n \choose s}^{p/2}} ~\ge~
\frac{\frac{1}{2^n} \cdot {n \choose {i_0}} K^p_s\(i_0\)}{{n \choose s}^{p/2}} ~\ge~
\]
\[
\(\frac{{n \choose s}}{2^{H\(\frac sn\) \cdot n}}\)^{\frac p2} \cdot \frac{1}{2^n} \cdot {n \choose {i_0}} 2^{\(p\tau\(\frac sn,y\) - \frac p2 H\(\frac sn\)\) \cdot n} ~\ge~
\]
\[
\Omega\(\frac{1}{n}\) \(\frac{{n \choose s}}{2^{H\(\frac sn\) \cdot n}}\)^{\frac p2} \cdot 2^{\(H(y) - 1 + p\tau\(\frac sn,y\) - \frac p2 H\(\frac sn\)\) \cdot n} ~=~ \Omega\(\frac{1}{n}\) \(\frac{{n \choose s}}{2^{H\(\frac sn\) \cdot n}}\)^{\frac p2} \cdot 2^{\psi\(p,\frac sn\) \cdot n} ~\ge~
\]
\[
\Omega\(\frac{1}{n}\) \cdot C^{-p} \cdot s^{-\frac p4} \cdot 2^{\psi\(p,\frac sn\) \cdot n}.
\]
\eprf

\section{Appendix: Proofs of claims about univariate and bivariate functions}

\subsubsection*{Proof of Lemma~\ref{lem:ref:r}}

\noi The derivative $\frac{\partial r}{\partial y}$ is easily seen to be proportional, up to a positive factor, to $(1-2x)^2 + (1-2x)\sqrt{(1-2x)^2 - 4y(1-y)} -  2(1-y)$. For a fixed $y$, this is maximized at $x = 0$, in which case this is $0$.
\eprf

\subsubsection*{Proof of Lemma~\ref{lem:ref:tau}}

\noi For $y \ge \frac12 - \sqrt{x(1-x)}$ the claim of the lemma follows immediately from the definition of $\tau$. It should also be possible to verify the claim directly for $y < \frac12 - \sqrt{x(1-x)}$ , but we proceed by observing that in this range the claim follows immediately from the reciprocity of Krawchouk polynomials (property 1 in Section~\ref{subsec:krawchouk}), from (\ref{Krawchouk-tail}), from (\ref{binomial-H}), and from the continuity of the function $\tau$.
\eprf

\subsubsection*{Proof of Lemma~\ref{lem:ref:h}}

\noi For the first claim of the lemma, the values of $h$ at the endpoints of $x$ are easy to verify. And, it is easy to see that for $0 < x < 1/2$ holds
\[
\frac{\partial h}{\partial x} ~=~ \frac{p-1}{p} \cdot \(\(\frac{1-x}{x}\)^{\frac{1}{p}} - \(\frac{x}{1-x}\)^{\frac{1}{p}}\) + \frac{1}{p} \cdot \(\(\frac{1-x}{x}\)^{\frac{p-1}{p}} - \(\frac{x}{1-x}\)^{\frac{p-1}{p}}\) ~ > ~ 0.
\]

\noi For the second claim of the lemma, writing $h(p,x) = \sqrt{x(1-x)} \cdot \(\(\frac{1-x}{x}\)^{\frac{p-2}{2p}} + \(\frac{x}{1-x}\)^{\frac{p-2}{2p}}\)$, it is easy to see that for a fixed $0 < x < 1/2$ this is a strongly increasing function in $p$. It is also easy to see that  $h(2,x) = 2 \sqrt{x(1-x)}$ and that $h(p,x) \rarrow_{p \rarrow \infty} 1$.

\noi For the third claim of the lemma, let $z$ be such that $h(2,z) = 1 - 2x$. Since $h(2,z) = 2 \sqrt{z(1-z)}$,
this is equivalent (after rearranging) to $x = \frac12 - \sqrt{z(1-z)}$, which is the same as $z = \frac12 - \sqrt{x(1-x)}$. Since $h(p,z)$ increases in $p$ this means that $h(p,z) > 1 - 2x$. Since $h(p,u)$ increases in $u$, and $h(p,y) = 1 - 2x$, this implies $y < z = \frac12 - \sqrt{x(1-x)}$.
\eprf

\subsubsection*{Proof of Lemma~\ref{lem:ref:psi-aux}}

\noi We view $a$ as a function of $\delta$ for a fixed $p$. The boundary values of $a$ are easy to verify. It remains to check that $a$ decreases. We will show that $b = 1 - 2a = \delta(1-\delta) \cdot \frac{(1-\delta)^{p-2} + \delta^{p-2}}{(1-\delta)^p + \delta^p}$ increases. Let $g(\delta) = (1 - \delta)^{p-2} + \delta^{p-2}$, and $h(\delta) = (1 - \delta)^p + \delta^p$. Then $a = \delta(1-\delta) \cdot \frac gh$, and $a' = \frac{(1 - 2\delta) gh + \delta(1-\delta) \(g'h - gh'\)}{h^2}$. We will show that the numerator is positive, which will imply $a' > 0$.

\noi Computing and simplifying, we have that $g'h - gh' = (1-2\delta)\big(\delta(1-\delta)\big)^{p-3} \cdot \big(2\delta(1-\delta) + (p-2)\big)$, and that $gh = (1-\delta)^{2p-2} + \delta^{2p-2} + \((1-\delta)^2 + \delta^2\) \big(\delta(1-\delta)\big)^{p-2}$. Substituting and simplifying, we get that
\[
(1 - 2\delta) gh + \delta(1-\delta) \(g'h - gh'\) ~=~ (1 - \delta)^{2p-2} - \delta^{2p-2} + (p-1) (1-2\delta) \big(\delta(1-\delta)\big)^{p-2},
\]
which is positive for all  $p \ge 2$ and $0 \le \delta < \frac12$.

\eprf

\subsubsection*{Proof of Proposition~\ref{pro:psi-second rep}}

\noi We need to show that
\beqn
\label{inden-psi-one-two}
H(y) - 1 + p\tau(x,y) - \frac p2 H(x) ~=~ (p-1) + \log_2\Big((1-\delta)^p + \delta^p\Big) - \frac{p}{2} H(x) - px\log_2(1-2\delta).
\eeqn

\noi We fix $p$ and view both sides as functions of a free variable $0 \le \delta \le \frac12$. Recall that $x = x(\delta) = \(\frac12 - \delta\) \cdot \frac{(1-\delta)^{p-1} - \delta^{p-1}}{(1-\delta)^p + \delta^p}$ and that $y = y(\delta)$ is determined by $h(p,y) = 1 - 2x$. In fact, we claim that $y(\delta) = \frac{\delta^p} {(1-\delta)^p + \delta^p}$. To see this, one has to verify $h(p,y) = 1 - 2x$, and this is easy to do.

\noi We start with an auxiliary lemma.
\lem
\label{lem:psi-I-aux}
With our definitions of $x$ and $y$, we have that
\[
I(0,x) - I(y,x) ~=~ x \log_2\(1-2\delta\) - 1 + H(x) - y \log_2(\delta) - (1-y)\log_2(1-\delta).
\]
\elem
\prf

\noi It is not hard to verify directly that we have

\begin{enumerate}

\item $1 - 2x = \delta(1-\delta) \cdot \frac{(1-\delta)^{p-2} + \delta^{p-2}}{(1-\delta)^p + \delta^p}$.

\item $\sqrt{(1-2x)^2 - 4y(1-y)} = \delta(1-\delta) \cdot \frac{(1-\delta)^{p-2} - \delta^{p-2}}{(1-\delta)^p + \delta^p}$.

\item $1 - 2y - \sqrt{(1-2x)^2 - 4y(1-y)} = (1-2\delta) \cdot \frac{(1-\delta)^{p-1} + \delta^{p-1}}{(1-\delta)^p + \delta^p}$.

\item $1 - 2x + \sqrt{(1-2x)^2 - 4y(1-y)} = 2\delta \cdot \frac{(1-\delta)^{p-1}} {(1-\delta)^p + \delta^p} = \frac{\delta}{1-\delta} \cdot (2 - 2y)$.

    This implies that $\frac{1 - 2x + \sqrt{(1-2x)^2 - 4y(1-y)}}{2 - 2y} = \frac{\delta}{1-\delta}$.

\item $2 - 2y - (1-2x)^2 - (1-2x) \cdot \sqrt{(1-2x)^2 - 4y(1-y)} = \frac{2-4\delta}{(1 - \delta)^2} \cdot (1-y)^2$.

\end{enumerate}

\noi Substituting this in the definition of $I$ leads, after some simplification, to
\[
I(0,x) - I(y,x) ~=~ \frac{1-2x}{2} \log_2\(\frac{(1-\delta)^{p-1} - \delta^{p-1}}{(1-\delta)^{p-1} + \delta^{p-1}}\) - y \log_2\(\frac{\delta}{1-\delta}\) - \frac12 \log_2\(8x(1-x)\) + \frac12 \log_2\(\frac{2-4\delta}{(1-\delta)^2}\).
\]

\noi We claim that the RHS of this expression can be further simplified to the RHS in the claim of the lemma. Indeed, simplifying, we need to verify
\[
x \log_2\(1-2\delta\) + H(x) ~=~ \(\frac12 - x\) \cdot \log_2\(\frac{(1-\delta)^{p-1} - \delta^{p-1}}{(1-\delta)^{p-1} + \delta^{p-1}}\) - \frac12 \log_2\(x(1-x)\) + \frac12 \log_2\(1 - 2\delta\).
\]

\noi Expanding the enropy and rearranging, this is the same as
\[
\(\frac12 - x\) \cdot \log_2\((1-2\delta) \cdot \frac{(1-\delta)^{p-1} - \delta^{p-1}}{(1-\delta)^{p-1} + \delta^{p-1}}\) + \(x - \frac12\) \log_2(x) + \(\frac12 - x\) \log_2(1-x) = 0,
\]
which is equivalent to $(1-2\delta) \cdot \frac{(1-\delta)^{p-1} - \delta^{p-1}}{(1-\delta)^{p-1} + \delta^{p-1}} = \frac{x}{1-x}$, and this is easy to verify directly.
\eprf

\noi We continue with the proof of the proposition. Recalling the definition of $\tau$, and substituting the identity proved in the previous lemma in the LHS of (\ref{inden-psi-one-two}), we get that it equals to
\[
(p-1) + H(y) - \frac p2 H(x) + p\cdot\Big(y\log_2(\delta) + (1-y) \log_2(1-\delta) - x \log_2(1-2\delta)\Big).
\]

\noi Finally, it is easy to see that $H(y) + p\cdot \(y \log_2(\delta) + (1-y) \log_2(1-\delta)\) = \log_2 \((1-\delta)^p + \delta^p\)$. Substituting this in the last expression gives the RHS of (\ref{inden-psi-one-two}).

\subsubsection*{Proof of Proposition~\ref{pro:psi-concave convex}}

\noi We start with the first claim of the proposition. For $x = 0$ we use the second definition of $\psi$. We have that $\delta = \frac12$, and hence $\psi(p,0) = p - 1 + \log_2\((1-\delta)^p + \delta^p\) = 0$. For $p = 2$ we use the first definition of $\psi$. Since $h\(2, y\) = 2\sqrt{y(1-y)}$, we have that $x = \frac12 - \sqrt{y(1-y)}$, which is the same as $y = \frac12 - \sqrt{x(1-x)}$. Hence, by the definition of $\tau$, we have $\tau(x,y) = \frac{1+H(x)-H(y)}{2}$. Substituting this in the definition of $\psi$, we have $\psi(2,x) = 0$.

\noi We proceed with the second claim of the proposition, using the first definition of $\psi$. We view $x$ as fixed, and write $g(p) = \psi(p,x) =  H(y) - 1 + p\tau(x,y) - \frac p2 H(x)$, where $y = y(p)$ is determined by $h(p,y) = 1 - 2x$. Note that for $x > 0$ holds $y < \frac12 - \sqrt{x(1-x)}$. And hence $\frac{\partial \tau(x,y)}{\partial y} = \log_2\(r(x,y)\)$. Therefore we have that
\[
g' ~=~ \frac{\partial y}{\partial p} \cdot \(\log_2\(\frac{1-y}{y}\) + p \frac{\partial \tau}{\partial y}\) + \tau(x,y) - \frac12 H(x) ~=~
\]
\[
\frac{\partial y}{\partial p} \cdot \(\log_2\(\frac{1-y}{y}\) + p \log_2\(r(x,y)\)\) + \tau(x,y) - \frac12 H(x).
\]

\noi Next, we claim that the expression in brackets vanishes.
\lem
Let $0 \le x,y \le \frac12$ be such that $0 \le x < \frac12 - \sqrt{y(1-y)}$. Then there is a unique $p > 2$ such that $h(p,y) = 1 - 2x$ and this $p$ is given by
\[
p ~=~ - \frac{\log_2\(\frac{1-y}{y}\)}{\log_2\(r(x,y)\)} = \frac{\log_2\(\frac{y}{1-y}\)}{\log_2\(r(x,y)\)},
\]
where $r(x,y) = \frac{(1-2x) + \sqrt{(1-2x)^2 - 4y(1-y)}}{2 - 2y}$.
\elem
\prf
By the properties of the function $h$, there is a unique $p > 2$ such that $h(p,y) = 1 - 2x$. So it suffices verify the identity $h(p,y) = 1 - 2x$ for $p$ given in the claim of the lemma.

\noi Writing $M$ for $(1-2x) + \sqrt{(1-2x)^2 - 4y(1-y)}$ and $r = \frac{M}{2 - 2y}$ for $r(x,y)$, we have that
\[
h(p,y) ~=~ y^{\frac 1p} (1-y)^{\frac{p-1}{p}} + y^{\frac{p-1}{p}} (1-y)^{\frac 1p} ~=~ y^{\frac12} (1-y)^{\frac12} \cdot \(\(\frac{1-y}{y}\)^{\frac{p-2}{2p}} + \(\frac{y}{1-y}\)^{\frac{p-2}{2p}}\) ~=~
\]
\[
 y^{\frac12} (1-y)^{\frac12} \cdot \(\frac{y^{\frac12}}{(1-y)^{\frac12} r} + \frac{(1-y)^{\frac12} r}{y^{\frac12}}\) ~=~ \frac{y + (1-y) r^2}{r} ~=~
\]
\[
\frac{M^2 + 4y(1-y)}{2M} ~=~ 1-2x.
\]
\eprf

\noi Using the lemma gives
\[
g' ~=~ \tau(x,y) - \frac12 H(x).
\]

\noi We claim that this is positive for any $p > 2$ and hence $g$ is increasing. Recall that $\tau(x,y)$ is decreasing in $y$. Since $0 \le y < \frac12 - \sqrt{x(1-x)}$, we have that $\tau(x,y) \ge \tau\(x, \frac12 - \sqrt{x(1-x)}\) = \frac{1 + H(x) - H\(\frac12 - \sqrt{x(1-x)}\)}{2} > \frac12 H(x)$. Hence $g' > 0$.

\noi Next, we have $g'' = \frac{\partial \tau}{\partial y} \cdot \frac{\partial y}{\partial p} > 0$, since both terms in the product are negative ($\tau$ decreases in $y$ and $h$ increases in $p$, while $y(p)$ is determined by $h(p,y) = 1 - 2x$, and $x$ is fixed). This means that $g$ is strongly convex, completing the proof of the second claim of the proposition.

\noi We pass to the third claim of the proposition, using the second definition of $\psi$. We view $p$ as fixed, and write $g(x) = \psi(p,x) = (p-1) + \log_2\Big((1-\delta)^p + \delta^p\Big) - \frac{p}{2} H(x) - px\log_2(1-2\delta)$, where $\delta$ is determined by $x = \(\frac12 - \delta\) \cdot \frac{(1-\delta)^{p-1} - \delta^{p-1}}{(1-\delta)^p + \delta^p}$. We have
\[
g' ~=~ \frac{1}{ \ln 2} \cdot \frac{1}{(1-\delta)^p + \delta^p} \cdot \(p\delta^{p-1} - p(1-\delta)^{p-1}\) \cdot \delta' - \frac p2 \log_2\(\frac{1-x}{x}\) - p \log_2(1-2\delta) + px\frac{1}{\ln 2} \cdot \frac{1}{1 - 2\delta} \cdot 2 \delta'
\]
Note that the first and the fourth terms cancel out, by the definition of $\delta$, and hence we get
\[
g' ~=~ -\frac p2 \log_2\(\frac{1-x}{x}\)  - p \log_2(1-2\delta).
\]
To show that $g$ is increasing amounts to showing that $\log_2 \(\frac{1-x}{x}\) + \log_2\((1-2\delta)^2\) < 0$, which, recalling $\frac{x}{1-x} = (1-2\delta) \cdot \frac{(1-\delta)^{p-1} - \delta^{p-1}}{(1-\delta)^{p-1} + \delta^{p-1}}$, is easily simplifiable to
\[
\frac{(1-\delta)^{p-1} - \delta^{p-1}}{(1-\delta)^{p-1} + \delta^{p-1}} ~>~ 1 - 2\delta,
\]
for any $p > 2$ and $0 < \delta < \frac12$. Note that both sides of this inequality coincide for $p=2$. We claim that the LHS increases with $p$. Indeed, by a simple calculation, the derivative of the LHS w.r.t. $p$ is proportional to $\ln\(\frac{1-\delta}{\delta}\) \cdot \big(\delta(1-\delta)\big)^{p-1}$, which is clearly positive.

\noi Computing the derivative of $g$ at $0$ gives, by L'Hospital,
\[
g'(0) ~=~ -\frac p2 \cdot \log_2\(\lim_{x \rarrow 0} \frac{(1-x)(1-2\delta)^2}{x}\) ~=~ \frac p2 \cdot \log_2\(\lim_{\delta \rarrow \frac12} \frac{(1-2\delta)^2}{x}\),
\]
as claimed.

\noi It is easy to see that the limit is $\frac{1}{p-1}$, and we get
\[
g'(0) ~=~ \frac{p \log_2(p-1)}{2}.
\]

\noi We proceed to argue that $g'' < 0$ for $x > 0$, and hence $g$ is strongly concave. It is easy to see that $g''$ is proportional to $\frac{1}{x(1-x)} + \frac{4 \delta'}{1-2\delta}$.
It will be convenient to state the inequality $g'' < 0$ in terms of $\delta$. Recalling the definition of $\delta$, and rearranging, we need to show that
\[
4x(1-x) ~\ge~ -(1-2\delta) x'(\delta),
\]
with equality holding only at $\delta = \frac12$. Here we write $x(\delta) = \(\frac12 - \delta\) \cdot \frac{(1-\delta)^{p-1} - \delta^{p-1}}{(1-\delta)^p + \delta^p}$. Note that both sides vanish at $\delta = \frac12$, so we need only to show strict inequality for $0 \le \delta < \frac12$.

\noi We now introduce some notation to make the following calculations easier to write. Let $L_+(p) = L_+(p, \delta) = (1-\delta)^p + \delta^p$ and let $L_-(p) = L_-(p, \delta) = (1-\delta)^p - \delta^p$. Then the inequality above transforms into (after some simplification):
\[
(2 - 4\delta) L_-(p-1) \cdot \(2L_+(p) - (1-2\delta) L_-(p-1)\)  ~>~
\]
\[
(p-1) (1-2\delta)^2 L_+(p-2) L_+(p) + (2-4\delta) L_-(p-1) L_+(p) - p(1 - 2\delta)^2 L^2_-(p-1).
\]
Simplifying, the LHS is $(2-4\delta) L_-(2p-2)$.

\noi We also have the following identities:
$L_+(p-2) L_+(p) = L_+(2p-2) + \((1-\delta)^2 + \delta^2\) \cdot \big((\delta(1-\delta)\big)^{p-2}$; $L_-(p-1) L_+(p) = L_-(2p-1) - (1-2\delta) \big((\delta(1-\delta)\big)^{p-1}$; $L^2_-(p-1) = L_+(2p-2) - 2 \big((\delta(1-\delta)\big)^{p-1}$.

\noi Substituting, collecting similar terms together, and simplifying, we get to
\[
2L_-(2p-2) + (1-2\delta) L_+(2p-2) - 2L_-(2p-1) ~>~ (p-1)(1-2\delta)\big(\delta(1-\delta)\big)^{p-2}.
\]
Expaning the "$L$" notation, this is the same as
\beqn
\label{aux-psi-x-concave}
\frac{(1-\delta)^{2p-2} - \delta^{2p-2}}{1-2\delta} ~>~ (p-1)\big(\delta(1-\delta)\big)^{p-2}.
\eeqn

\noi To show this, we start with an auxiliary claim.
\lem
\label{lem:aux-a-psi}
Let $p > 2$. Then the function $f(\delta) = \frac{(1-\delta)^p - \delta^p}{1-2\delta}$ decreases on $\left[0, \frac12\)$.
\elem
\prf
Computing the derivative and simplifying, we have that $f'$ is proportional to $p \cdot \(\delta(1-\delta)^{p-1} -\delta^{p-1}(1-\delta)\) - (p-2) \cdot \((1-\delta)^p - \delta^p\)$. So, we need to show that the second term is greater than the first one. Both are equal to zero at $\delta = \frac12$, so it suffices to show that the derivative of the second term is smaller than that of the first term, which, after simplifying, amounts to $(1-\delta)^{p-1} + \delta^{p-1} > \delta(1-\delta)^{p-2} + \delta^{p-2}(1-\delta)$, which is true for $p > 2$, and for $0 \le \delta < \frac12$.
\eprf

\noi Now consider (\ref{aux-psi-x-concave}). Note that for $\delta = \frac12$ both sides (LHS at the limit for $\delta \rarrow \frac12$) equal $(p-1) \(\frac12\)^{2p - 4}$. In addition, by the preceding lemma, the LHS decreases, while it is easy to see that the RHS increases in $\delta$.
\eprf

\subsubsection*{Proof of Lemma~\ref{lem:ref:pi}}

\noi First,
\[
\pi(x,y) ~=~ \tau(x,y) - \frac{1 + H(x) - H(y)}{2} ~=~ \tau(y,x) - \frac{1 - H(x) + H(y)}{2} ~=~ \pi(y,x).
\]
The second equality is by Lemma~\ref{lem:ref:tau}.

\noi Next, by the proof of Lemma~\ref{lem:ref:pi-min}, we have that for $y < \frac12 - \sqrt{x(1-x)}$ the derivative $\frac{\partial \pi(x,y)}{\partial y} = -\log_2(1-2\delta)$, where $\delta = \delta(x,y)$ is easily seen to be strictly positive if $y > 0$. Hence $\pi$ is strongly increasing in $y$ and, by symmetry, also in $x$.

\noi Finally, recall that $\pi$ is continuous on $\left[0, \frac12\right] \times \left[0, \frac12\right]$ and that $\pi(x,y) = 0$ if $y \ge \frac12 - \sqrt{x(1-x)}$. This implies that $\pi(x,y) < 0$ for $y < \frac12 - \sqrt{x(1-x)}$.
\eprf

\subsubsection*{Proof of Lemma~\ref{lem:ref:pi-min}}

\noi Fix $\sigma$ and $\kappa$, and let $F(\delta) = \sigma H\(\frac{x}{\sigma}\) + \(1-\sigma\) H\(\frac{x}{1-\sigma}\) + 2x\log_2(\delta) + (1-2x) \log_2(1-\delta) - \kappa \log_2(1-2\delta)$. Then
\[
F'(\delta) ~~=~~ \frac{\partial x}{\partial \delta} \cdot \log_2\(\frac{\delta^2 (\sigma-x)(1-\sigma-x)}{(1-\delta)^2 x^2}\) + \frac{1}{\ln 2} \cdot \(\frac{2x}{\delta} - \frac{1-2x}{1-\delta} + \frac{2\kappa}{1-2\delta}\) ~~=
\]
\[
\frac{1}{\ln 2} \cdot \(\frac{2x}{\delta} - \frac{1-2x}{1-\delta} + \frac{2\kappa}{1-2\delta}\) ~~=~~ \frac{1}{(\ln 2) \delta(1-\delta)(1-2\delta)} \cdot \Big(2x(1-2\delta) + 2\kappa \delta(1-\delta) - \delta(1-2\delta)\Big),
\]
where the first equality follows since it is easy to check that the term multiplying $\frac{\partial x}{\partial \delta}$ vanishes by the definition of $x$.

\noi Recalling the definition of $x$ and simplifying, we have
\[
2x(1-2\delta) + 2\kappa \delta(1-\delta) - \delta(1-2\delta) ~~=~~ \delta \cdot \(\sqrt{\delta^2 + 4\sigma(1-\sigma)(1-2\delta)} - (1-2\kappa)(1-\delta)\).
\]
Now there are two cases to consider.
\begin{enumerate}

\item $\kappa \ge \frac12 - \sqrt{\sigma(1-\sigma)}$.

\noi In this case it is easy to see, by squaring both sides and analyzing the obtained quadratic inequality, that $\sqrt{\delta^2 + 4\sigma(1-\sigma)(1-2\delta)} \ge (1-2\kappa)(1-\delta)$ for all $0 \le \delta \le \frac12$, and hence $F'(\delta) \ge 0$ for all $\delta$. It follows that $F$ is increasing and its minimum is given by $F(0) = 0$. On the other hand, by the definition of $\pi$, in this case $\pi(\sigma,\kappa) = 0$, and the claim of the lemma holds.

\item $\kappa < \frac12 - \sqrt{\sigma(1-\sigma)}$. In this case, again analyzing the appropriate quadratic inequality, it is easy to see that the minimum of $F$ is attained at the only zero of $F'$ on $\left[0,\frac12\right]$, that is at
    \[
    \delta = \delta(\sigma,\kappa) ~~=~~ \frac{(1-2\sigma) \sqrt{(1-2\sigma)^2 - 4\kappa(1-\kappa)} - \Big((1-2\sigma)^2 - 4\kappa(1-\kappa)\Big)}{4\kappa (1-\kappa)}.
    \]

\noi So, we need to verify that for this value of $\delta$ holds
\[
\sigma H\(\frac{x}{\sigma}\) + \(1-\sigma\) H\(\frac{x}{1-\sigma}\) + 2x\log_2(\delta) + (1-2x) \log_2(1-\delta) - \kappa \log_2(1-2\delta)  = 2\pi\(\sigma,\kappa\).
\]

\noi Let $L(\sigma,\kappa)$ denote the LHS of the above. We need to verify $L = 2\pi$. First, note that for $\kappa = \frac12 - \sqrt{\delta(1-\delta)}$, we have $\delta(\sigma,\kappa) = 0$ and both $L$ and $\pi$ vanish. So we have an equality at an endpoint, and hence it suffices to show that the derivatives $\frac{\partial L}{\partial \kappa}$ and $\frac{\partial (2\pi)}{\partial \kappa}$ coincide. We have
\[
\frac{\partial L}{\partial \kappa} =  \frac{\partial x}{\partial \kappa} \cdot \log_2\(\frac{\delta^2 (\sigma-x)(1-\sigma-x)}{(1-\delta)^2 x^2}\) + \frac{\partial \delta}{\partial \kappa} \cdot \frac{1}{\ln 2} \cdot \(\frac{2x}{\delta} - \frac{1-2x}{1-\delta} + \frac{2\kappa}{1-2\delta}\) - \log_2(1-2\delta) =
\]
\[
- \log_2(1-2\delta).
\]
To see the equality, note that the first summand vanishes by the definition of $x$, and the second summand vanishes by the definition of $\delta$.

\noi On the other hand, by the definition of $\pi$, we have $2\pi(\sigma,\kappa) = H(\sigma) + H(\kappa) - 1 + 2I(\kappa,\sigma) - 2I(0,\sigma)$. Hence, using the notation of Section~\ref{subsubsec:ref:I},
\[
\frac{\partial (2\pi)}{\partial \kappa} ~=~ \log_2\(\frac{1-\kappa}{\kappa}\) + 2 \log_2\(r(\sigma,\kappa)\) ~=~
\]
\[
\log_2\(\frac{1-\kappa}{\kappa}\) + 2 \log_2\(\frac{(1-2\sigma) + \sqrt{(1-2\sigma)^2 - 4\kappa(1-\kappa)}}{2 - 2\kappa}\) ~=~
\]
\[
\log_2\(\frac{\((1-2\sigma) + \sqrt{(1-2\sigma)^2 - 4\kappa(1-\kappa)}\)^2}{4\kappa(1-\kappa)}\).
\]

\noi Let $C(\sigma,\kappa) = (1-2\sigma)^2 - 4\kappa(1-\kappa)$. Then the last expression can be written as $\log_2\(\frac{\((1-2\sigma) + \sqrt{C}\)^2}{4\kappa(1-\kappa)}\)$. We can also write $\delta = \frac{(1-2\sigma) \sqrt{C} - C}{4\kappa(1-\kappa)}$. It is easy to see that this implies $1 - 2\delta = \frac{\((1-2\sigma) - \sqrt{C}\)^2}{4\kappa(1-\kappa)}$.

\noi It remains to observe that
\[
(1-2\delta) \cdot \frac{\((1-2\sigma) + \sqrt{C}\)^2}{4\kappa(1-\kappa)}  ~~=~~ \frac{\((1-2\sigma) - \sqrt{C}\)^2}{4\kappa(1-\kappa)} \cdot \frac{\((1-2\sigma) + \sqrt{C}\)^2}{4\kappa(1-\kappa)}  ~~=~~ 1,
\]
which means that
\[
\frac{\partial L}{\partial \kappa} = \log_2\(\frac{1}{1-2\delta}\) = \log_2\(\frac{\((1-2\sigma) + \sqrt{C}\)^2}{4\kappa(1-\kappa)}\) = \frac{\partial (2\pi)}{\partial \kappa}.
\]
\end{enumerate}

\noi This completes the proof of the lemma.

\eprf

\subsubsection*{Proof of Lemma~\ref{lem:ref:phi}}

\noi For fixed $\sigma$ and $\e$, let $\alpha(x) =  \sigma H\(\frac{x}{\sigma}\) + (1-\sigma) H\(\frac{x}{1-\sigma}\) + 2x\log_2(\e) + (1-2x) \log_2(1-\e)$. As observed in \cite{ACKL}, the maximum of $\alpha$ is attained at the only zero of $\alpha'$, that is at $x = x(\sigma, \e) = \frac{-\e^2 + \e \sqrt{\e^2 + 4(1-2\e) \sigma(1-\sigma)}}{2(1-2\e)}$. For this value of $x$ holds $\phi(\sigma,\e) = H(\sigma) - 1 + \alpha(x)$.

\noi Similarly, for fixed $\sigma$ and $\e$, let $\beta(y) = y \log_2(1 - 2\e) + H(y) + 2\tau(\sigma,y)$. We are interested in the maximum of $\beta$ on $0 \le y \le \frac12$. First, note that, by the definition of $\tau$, we have $\beta(y) = \log_2(1-2\e) \cdot y + 1 + H(\sigma)$ for $y \ge \frac12 - \sqrt{\sigma(1-\sigma)}$, and hence $\beta$ decreases for $\frac12 - \sqrt{\sigma(1-\sigma)} \le y \le \frac12$. For $0 \le y < \frac12 - \sqrt{\sigma(1-\sigma)}$, we have $\beta' = \log_2(1-2\e) + \log_2\(\frac{1-y}{y}\) + 2\frac{\partial \tau}{\partial y} = (1-2\e) + \log_2\(\frac{1-y}{y}\) + 2\log_2\(r(\sigma,y)\)$. It is easy to see that the maximum of $\beta$ is attained at the only zero of $\beta'$, that is at $y = y(\sigma, \e) = \frac{(1-\e) - \sqrt{\e^2 + 4(1-2\e) \sigma(1-\sigma)}}{2 - 2\e}$.

\noi Write $\phi_2(\sigma,\e)$ for $\beta(y(\sigma,\e)) - 2$. We need to verify $\phi  = \phi_2$. First, we check the boundary conditions $\phi(\sigma,0) = \phi_2(\sigma,0)$ for all $0 \le \sigma \le \frac12$. We have $x(\sigma,0) = 0$ and hence $\phi(\sigma,0) = H(\sigma) - 1 + \alpha(0) = H(\sigma) - 1$. On the other hand, $y(\sigma,0) = \frac12 - \sqrt{\sigma(1-\sigma)}$. We have $\tau\(\sigma, \frac12 - \sqrt{\sigma(1-\sigma)}\) = \frac{1 + H(\sigma) - H\(\frac12 - \sqrt{\sigma(1-\sigma)}\)}{2}$, and hence $\phi_2(\sigma,0) = \beta\(\frac12 - \sqrt{\sigma(1-\sigma)}\) - 2 = H\(\frac12 - \sqrt{\sigma(1-\sigma)}\) + 2\tau\(\sigma, \frac12 - \sqrt{\sigma(1-\sigma)}\) - 2 = H(\sigma) - 1$ as well.

\noi Next, we verify that $\frac{\partial \phi}{\partial \e} = \frac{\partial \phi_2}{\partial \e}$, which will complete the proof. Writing $x$ for $x(\sigma,\e)$ and $y$ for $y(\sigma,\e)$, it is easy to see that $\frac{\partial \phi}{\partial \e} = \frac{\partial x}{\partial \e} \cdot \(\frac{\partial \alpha}{\partial x}_{|x(\sigma,\e)}\) + \frac{2x - \e}{\ln(2) \e(1-\e)} = \frac{2x - \e}{\ln(2)\e(1-\e)}$. Similarly, $\frac{\partial \phi_2}{\partial \e} = \frac{\partial y}{\partial \e} \cdot \(\frac{\partial \beta}{\partial y}_{|y = y(\sigma,\e)}\) - \frac{2y}{\ln(2) (1-2\e)} = - \frac{2y}{\ln(2) (1-2\e)}$. So, it remains to verify $\frac{2x - \e}{\e(1-\e)} = -\frac{2y}{1-2\e}$, which is easy to do directly.

\eprf

\subsubsection*{Proof of Lemma~\ref{lem:ref:phi:edge-isop}}

\noi Recall that $\phi(\sigma,\e) = H(\sigma) - 1 + \sigma H\(\frac{x}{\sigma}\) + (1-\sigma) H\(\frac{x}{1-\sigma}\) + 2x\log_2(\e) + (1-2x) \log_2(1-\e)$, where $x = x(\sigma, \e) = \frac{-\e^2 + \e \sqrt{\e^2 + 4(1-2\e) \sigma(1-\sigma)}}{2(1-2\e)}$. Substituting, we need to show that
\[
\min_{0 < \e \le \frac12} \left\{\sigma H\(\frac{x}{\sigma}\) + (1-\sigma) H\(\frac{x}{1-\sigma}\) + (2x-y)\log_2\(\frac{\e}{1-\e}\) \right\}  =
\sigma H\(\frac{y}{2\sigma}\) + (1-\sigma) H\(\frac{y}{2(1 - \sigma)}\).
\]

\noi Fix $\sigma$ and $y$, and let $F(\e) = \sigma H\(\frac{x}{\sigma}\) + (1-\sigma) H\(\frac{x}{1-\sigma}\) + (2x-y)\log_2\(\frac{\e}{1-\e}\)$. Then
\[
F'(\e) ~=~ \frac{\partial x}{\partial \e} \cdot \log_2\(\frac{\e^2 (\sigma-x)(1-\sigma-x)}{(1-\e)^2 x^2}\) + (2x-y) \cdot \frac{1}{\ln 2} \frac{1-2\e}{\e(1-\e)} ~=~ (2x-y) \cdot \frac{1}{\ln 2} \frac{1-2\e}{\e(1-\e)},
\]
since the term multiplying $\frac{\partial x}{\partial \e}$ vanishes by the definition of $x$.

\noi It is not hard to verify that $x$ strictly increases in $\e$ from $0$ to $\sigma(1-\sigma)$, and hence $F$ has a unique minimum at $\e$ for which $x = \frac y2$. Substituting this value of $x$ in $F$, gives the claim of the lemma.
\eprf

\subsubsection{Proof of Lemma~\ref{lem:ref:phi:disc-cont}}

\noi We write $\alpha$ for $\alpha_{\sigma,\e}$. Let $x = x^{\ast}(\sigma, \e) = \frac{-\e^2 + \e \sqrt{\e^2 + 4(1-2\e) \sigma(1-\sigma)}}{2(1-2\e)}$ be the point of maximum of $\alpha$, that is $A = \alpha\(x^{\ast}\)$. We distinguish between four cases.

\begin{enumerate}

\item $\sigma = 0$. In this case $x = 0$ as well, and $A = \alpha(0) = B$, and the claim holds. So in the remaining cases we may and will assume $\sigma \ge \frac 1n$.

\item $0 < x \le \frac{1}{5n}$. In this case we will compare $A$ with $B' = \alpha(0) = \log_2(1-\e)$. Clearly $A \ge B \ge B'$ and hence $A - B \le A - B'$. Recall that $x$ satisfies: $(1-\e)^2 x^2 = \e^2 (\sigma - x)(1 - \sigma - x)$. So, $\e^2 = \Theta\(\frac{x^2}{\sigma - x}\)$, with asymptotic notation hiding absolute constants. Expanding the entropies, and using the fact that $\sigma \ge \frac 1n \ge 5x$, it is easy to verify that $A - B' ~=~ O\(\frac 1n\)$ in this case.

\item $x > \frac{1}{5n}$ and $\sigma - x < \frac 2n$.

\noi In this case we will compare $A$ with $B' = \alpha(\sigma) = \alpha\(\frac{s}{n}\)$. It is not hard to verify (expanding the entropies and rearranging) that in this case
\[
A - B' ~\le~ x \log_2\(\frac{\sigma}{x}\) + (\sigma-x)\log_2\(\frac{\sigma^2}{(1-2\sigma)x^2}\) + O\(\frac 1n\).
\]
Observe that $x \log_2\(\frac{\sigma}{x}\) = O\(\frac 1n\)$, and that $\frac{\sigma^2}{(1-2\sigma)x^2}$ is bounded from above by an absolute constant. To see the second claim, note that $\frac{\sigma}{x}$ is bounded, and that the maximal value of $x$, attained for $\sigma = \frac12$ is $\frac{\e}{2}$, which is at most $\frac14$. So $\sigma$ cannot be close to  $\frac12$. Taking all of this into account, we get $A - B' ~=~ O\(\frac 1n\)$.

\item $x > \frac{1}{5n}$ and $\sigma - x > \frac 2n$.

\noi Recall that $\alpha'(y) = \log_2\(\frac{\e^2 (\sigma - y) (1 - \sigma - y)}{(1-\e)^2 y^2}\)$. We will choose $y$ to be the nearest fraction of the form $\frac in$ approximating $x$ from above, and set $B' = \alpha(y)$.  Then $|x - y| \le \frac 1n$ and for any $z$ between $x$ and $y$ holds
\[
\frac{\e^2 (\sigma - z) (1 - \sigma - z)}{(1-\e)^2 z^2} ~\le~ O\(\frac{\e^2 (\sigma - x) (1 - \sigma - x)}{(1-\e)^2 x^2}\) ~\le~ O(1).
\]
To see the first inequality note that all the terms in the first expression change by at most a costant factor compared to the second expression. For the second inequality, recall that $\frac{\e^2 (\sigma - x) (1 - \sigma - x)}{(1-\e)^2 x^2} = 1$. Hence $\alpha'(z) \le O(1)$ for all $y \le z \le x$ and hence $A - B \le A - B' = \alpha(x) - \alpha(y) \le O\(\frac 1n\)$.

\end{enumerate}

\eprf

\subsubsection{Proof of Lemma~\ref{lem:ref:tilde-phi}}

\noi Recall that $\phi(\sigma,\e) = H(\sigma) -1 + \sigma H\(\frac{x}{\sigma}\) + (1-\sigma) H\(\frac{x}{1-\sigma}\) + 2x\log_2(\e) + (1-2x) \log_2(1-\e)$, where $x = x(\sigma, \e) = \frac{-\e^2 + \e \sqrt{\e^2 + 4(1-2\e) \sigma(1-\sigma)}}{2(1-2\e)}$.

\noi We have that
\[
\frac{\partial \phi}{\partial \sigma} ~~=~~ \log_2\(\frac{(\sigma-x)(1-\sigma-x) \e^2}{(1-\e)^2 x}\) \cdot
\frac{\partial x^2}{\partial \sigma} + \log_2\(\frac{1-\sigma}{\sigma}\) + H\(\frac{x}{\sigma}\) - \frac{x}{\sigma} \log_2\(\frac{\sigma - x}{x}\) - H\(\frac{x}{1-\sigma}\) +
\]
\[
\frac{x}{1-\sigma} \log_2\(\frac{1-\sigma-x}{x}\).
\]
The first summand vanishes, since its first term vanishes by the definition of $x$ , and it is easy to see that the rest can be simplified to $\log_2\(\frac{1 - \sigma - x}{\sigma - x}\)$. Hence $\frac{\partial \phi}{\partial \sigma} ~~=~~ \log_2\(\frac{1 - \sigma - x}{\sigma - x}\)$. Since $\tilde{\phi}(\alpha,\e) = \phi\(H^{-1}(\alpha),\e\)$, we have
\[
\frac{\partial \tilde{\phi}}{\partial \alpha} ~~=~~ \(H^{-1}\)'(\alpha) \cdot \frac{\partial \phi}{\partial \sigma}_{|\sigma = H^{-1}(\alpha)} ~~=~~ \frac{\ln\(\frac{1 - \sigma - x}{\sigma - x}\)}{\ln\(\frac{1-\sigma}{\sigma}\)},
\]
where $\sigma = H^{-1}(\alpha)$. Computing the second derivative, we have that, similarly,
\[
\frac{\partial^2 \tilde{\phi}}{\partial \alpha^2} ~~=~~ \frac{1}{\log_2\(\frac{1-\sigma}{\sigma}\)} \cdot \frac{\partial}{\partial \sigma}_{|\sigma = H^{-1}(\alpha)} ~\frac{\ln(1 - \sigma - x) - \ln(\sigma - x)}{\ln(1-\sigma) - \ln(\sigma)},
\]
where $\sigma = H^{-1}(\alpha)$.

\noi In order to show that $\tilde{\phi}$ is concave, we need to show that $\frac{\partial}{\partial \sigma} \frac{\ln(1 - \sigma - x) - \ln(\sigma - x)}{\ln(1-\sigma) - \ln(\sigma)} \le 0$. Computing the derivative and rearranging, it is easy to see that this is equivalent to (writing $x'$ for $\frac{\partial x}{\partial \sigma}$):
\[
\sigma(1-\sigma) \ln\(\frac{1-\sigma}{\sigma}\) ~~\ge~~ \frac{(\sigma-x)(1-\sigma-x)}{(1-2x) - (1-2\sigma) x'} \cdot \ln\(\frac{1-\sigma-x}{\sigma-x}\).
\]
Next, note that $x' = \frac{1-2\sigma}{2 \frac{(1-\e)^2 - \e^2}{\e^2} \cdot x + 1}$. Substituting and using the fact that $(\sigma-x)(1-\sigma-x) \e^2 = (1-\e)^2 x^2$, the first term on the right can be simplified to $\sigma(1-\sigma) - \frac12 x$. Observing that, for fixed value of $\sigma$, the value of $x(\sigma, \e)$ increases from $0$ to $\sigma(1-\sigma)$, as $\e$ goes from $0$ to $\frac12$, it remains to verify that
\beqn
\label{phi-last}
\sigma(1-\sigma) \ln\(\frac{1-\sigma}{\sigma}\) ~~\ge~~ \(\sigma(1-\sigma) - \frac12 x\) \cdot \ln\(\frac{1-\sigma-x}{\sigma-x}\),
\eeqn
for any $0 \le \sigma \le \frac12$ and $0 \le x \le \sigma(1-\sigma)$. We proceed to show this. For a fixed $\sigma$, let $f(x)$ denote the RHS of this inequality. It is easy to see that $f(0) = f(\sigma(1-\sigma)) = \sigma(1-\sigma) \ln\(\frac{1-\sigma}{\sigma}\)$. We claim that $f$ is convex, which will imply (\ref{phi-last}). Indeed, direct calculation shows that
\[
f'' ~~=~~ \frac{(1 - 2\sigma) x}{(\sigma-x)^2 (1-\sigma-x)^2} \cdot \(\frac12 - 2\sigma(1-\sigma)\) ~~\ge~~ 0.
\]
This completes the proof of concavity of $\tilde{\phi}$.

\noi Let us also observe, for application in the proof of Lemma~\ref{lem:ref:eta-p} below,  that $f'' > 0$ for all $0 < x < \sigma(1-\sigma)$, which means that the inequality in (\ref{phi-last}) is strong for all $0 < x < \sigma(1-\sigma)$. This means that $\frac{\partial^2 \tilde{\phi}}{\partial \alpha^2} < 0$, for any $0 < \alpha < 1$, assuming $0 < \e < \frac12$.

\noi Next, we compute $\tilde{\phi}$ at $1$. Note that $x\(\frac12, \e\) = \frac{\e}{2}$, and hence $\tilde{\phi}(1,\e) = \phi\(\frac12, \e\) = H(\e) - H(\e) = 0$.

\noi We proceed to compute the right derivative of $\tilde{\phi}$ at $0$. Using the calculations above, and recalling that $x$ is between $0$ and $\sigma(1-\sigma)$, we have, by two applications of L'Hospital's rule, writing $\tilde{\phi}'(0,\e)$ for the right derivative at $0$:
\[
\tilde{\phi}'(0,\e) ~~=~~ \lim_{\sigma \rarrow 0} \frac{\ln\(\frac{1 - \sigma - x}{\sigma - x}\)}{\ln\(\frac{1-\sigma}{\sigma}\)} ~~=~~ \lim_{\sigma \rarrow 0} \left[\sigma(1-\sigma) \cdot \frac{(1-2x) - (1-2\sigma) x'}{(\sigma-x)(1-\sigma-x)}\right] ~~=~~
\]
\[
\lim_{\sigma \rarrow 0} \frac{\sigma(1-\sigma)}{\sigma(1-\sigma) - \frac12 x} ~~=~~ \lim_{\sigma \rarrow 0} \frac{1}{1 - \frac12 x'} ~~=~~ 2,
\]
where in the last equality we have used $\lim_{\sigma \rarrow 0} x' = \lim_{\sigma \rarrow 0} \frac{1-2\sigma}{2 \frac{(1-\e)^2 - \e^2}{\e^2} \cdot x + 1} = 1$.

\noi To compute the left derivative of $\tilde{\phi}$ at $\frac12$, recall that $x\(\frac12, \e\) = \frac{\e}{2}$, and hence, proceeding as in the preceding calculation,
\[
\tilde{\phi}'(1,\e) ~~=~~ \lim_{\sigma \rarrow \frac12} \frac{\ln\(\frac{1 - \sigma - x}{\sigma - x}\)}{\ln\(\frac{1-\sigma}{\sigma}\)} ~~=~~ \lim_{\sigma \rarrow \frac12} \frac{\sigma(1-\sigma)}{\sigma(1-\sigma) - \frac12 x} ~~=~~ \frac{1}{1-\e}.
\]

\noi Finally, we observe that since $\tilde{\phi}$ is concave, the value of its derivative on $(0,1)$ is bounded from below by  $\tilde{\phi}'(1,\e) =  \frac{1}{1-\e}$ and hence it is strongly increasing. This completes the proof of the lemma.
\eprf

\subsubsection{Proof of Lemma~\ref{lem:ref:eta-p}}

\noi We rely on the results in Lemma~\ref{lem:ref:tilde-phi}. The concavity of $\eta_p(x,\e)$ in $x$ follows immediately from that of $\tilde{\phi}$. Next, we compute the derivative of $\eta_p$ w.r.t. $x$. Since $\tilde{\phi}(y,\e)$ is concave in $y$, its derivative w.r.t. to $y$ is lower-bounded by the derivative at $y=1$, that is by $\frac{1}{1-\e}$. Hence
\[
\frac{\partial \eta_p}{\partial x} ~=~ -\frac{p}{2p-2} \frac{\partial \tilde{\phi(y,\e)}}{\partial y} \(1 - \frac{p}{p-1} \cdot x, ~2\e(1-\e)\) + \frac{1}{p-1} ~\le~ -\frac{p}{2p-2} \cdot \frac{1}{1 - 2\e(1-\e)} + \frac{1}{p-1} ~=
\]
\[
\frac{1}{p-1} \cdot \(-\frac{p}{1 + (1-2\e)^2} + 1\) ~\le~ 0.
\]

\noi Therefore, $\eta_p$ is decreasing in $x$.

\noi Let now $0 < \e < \frac12$.  the second derivative of $\tilde{\phi}(x,\e)$ w.r.t. $x$ is negative for all $0 < x < 1$, and hence the inequality in the above computation is sharp for all $x < \frac{p-1}{p}$. Hence $\eta_p$ is strongly decreasing in $x$. Observing that $\eta_p(0,\e) = \frac12 \tilde{\phi}\(1, ~2\e(1-\e)\) = 0$, this means that $\eta_p$ is strictly negative for $0 < x \le \frac{p-1}{p}$.

\eprf

\subsubsection*{Acknowledgement}

\noi We are grateful to Yuzhou Gu, Elchanan Mossel, Or Ordentlich, and Yury Polyanskiy for valuable remarks.

\end{document}